\documentclass{article}
\usepackage[utf8]{inputenc}
\usepackage{amsmath}
\usepackage{amssymb}
\usepackage{stmaryrd}
\usepackage{amsthm}
\usepackage{bbold}
\usepackage[left=3cm,right=3cm,top=3cm,bottom=3cm]{geometry}
\usepackage[linkcolor=blue, colorlinks=true]{hyperref}
\usepackage{graphicx}
\usepackage{appendix}
\usepackage{tikz}
\usepackage{subcaption}
\usetikzlibrary[patterns]

\newcommand{\footremember}[2]{
    \footnote{#2}
    \newcounter{#1}
    \setcounter{#1}{\value{footnote}}
}
\newcommand{\footrecall}[1]{
    \footnotemark[\value{#1}]
} 

\newcounter{count}
\newtheorem{hyp}{Assumption}

\newtheorem{theorem}{Theorem}
\newtheorem{prop}{Proposition}
\newtheorem{lemme}[prop]{Lemma}
\newtheorem{corr}[prop]{Corollary}
\newtheorem{conj}{Conjecture}

\title{An individual-based stochastic model reveals strong constraints on allometric relationships with minimal metabolic and ecological assumptions}
\author{
  Sylvain Billiard\footremember{alley}{Univ. Lille, CNRS, UMR 8198 – Evo-Eco-Paleo, F-59000 Lille, France}, Virgile Brodu\footremember{trailer}{Université de Lorraine, CNRS, Inria, IECL, F-54000 Nancy, France}, Nicolas Champagnat\footrecall{trailer}, Coralie Fritsch\footrecall{trailer}
  }

\begin{document}
\maketitle
\begin{abstract}

We design a stochastic individual-based model structured in energy, for single species consuming an external resource, where populations are characterized by a typical energy at birth in $\mathbb{R}^{*}_{+}$. The resource is maintained at a fixed amount, so we benefit from a branching property at the population level. Thus, we focus on individual trajectories, constructed as Piecewise Deterministic Markov Processes, with random jumps modelling births and deaths in the population; and a continuous and deterministic evolution of energy between jumps. We are mainly interested in the case where metabolic (\textit{i.e.} energy loss for maintenance), growth, birth and death rates depend on the individual energy over time, and follow allometric scalings (\textit{i.e.} power laws). Our goal is to determine in a bottom-up approach what are the possible allometric coefficients (\textit{i.e.} exponents of these power laws) under elementary --and ecologically relevant-- constraints, for our model to be valid for the whole spectrum of possible body sizes. We show in particular that assuming an allometric coefficient $\alpha$ related to metabolism strongly constrains the range of possible values for the allometric coefficients $\beta$, $\delta$, $\gamma$, respectively related to birth, death and growth rates. We further identify and discuss the precise and minimal ecological mechanisms that are involved in these strong constraints on allometric scalings.

\end{abstract}


\section{Introduction}

Since the synthesis work of Peters \cite{peters_1983}, a lot of papers in the field of evolutionary ecology \cite{yodzis_1992,brown_1993,denechere_2022} have stood for a very particular relationship between body mass or length and metabolism, across several orders of magnitude of species size. It is now well-known as an \textit{allometric} relationship, or simply \textit{allometry}, presented in the form $B \propto M^{\alpha} $, where $B$ is a measure of the metabolic rate, $M$ the mass and $\alpha$ the so-called \textit{allometric} coefficient. There is no clear consensus on a precise value for $\alpha$; some researchers argue that there exists a global coefficient based on a precise modelling of vascular systems \cite{savdee_08}, whereas others show that the coefficient may depend on specific features within species, and vary across the broad range of living organisms \cite{delong2010}. This kind of relationship is valid for metabolism, but also for a large number of ecological parameters (such as birth and death rates), with different allometric coefficients \cite{malerba_2019}. This allometric structure is very precise and often presented as
\begin{align}
\beta=\delta=\gamma-1=\alpha-1,
\label{eq:equalallo}
\end{align} 
with $\beta$, $\delta$, $\gamma$, $\alpha$ being allometric coefficients related to, respectively, the birth, death, growth and metabolic rates in the population. Even if these allometries seem to be a key ingredient for modelling the behavior of species and their ecological features, it is still mainly justified by raw data and phenomenological approaches \cite{niklas2001invariant}. Other studies are based on physical and physiological arguments for individuals, extrapolated at the population level \cite{brown_04}. The Metabolic Theory of Ecology justifies \textit{a posteriori} these relationships with energetic optimization and dimensional arguments \cite{savdee_08}, obliterating that evolutionary and ecological processes should be taken into account to fully explain this phenomenon. The principle of this theory is to understand ecological dynamics in terms of energy transfer and mass variation, with allometric scalings involved in the processes. These allometries are \textit{interspecific}, in the sense that a very broad range of living species apparently follows the same allometric laws, even if the fluctuation of their energies is on a different scale \cite{brown_04}. This is distinct from \textit{intraspecific} dynamics for fixed species, where numerous puzzling facts subsist \cite{koz_1997}, \cite{malerba_2019}.
\\\\
In this article, we will deliberately use the word \textit{energy} instead of \textit{mass} or \textit{size}. Indeed, biologists measure transfer of biomass rather than the masses themselves when they estimate metabolic rates of individuals \cite{yodzis_1992, denechere_2022}. This is why we prefer this generic term of \textit{energy}, associated to a fluctuating stock of resources meant to be invested by an individual for subsistence or reproduction. What interests us here is mostly the variety of reasonable assumptions we can make for the conception of a general interspecific model in the spirit of the Metabolic Theory of Ecology, and their consequences on the allometric scalings involved. Our goal is not to explain the emergence of allometries and their evolutionary aspects, but rather to identify admissible values for allometric coefficients, for a given fundamental coefficient $\alpha$ linking metabolism and individual energy. Also, we want to give probabilistic arguments and purely mathematical constraints, rather than the usual approach based on agreeing with experiment results.
\\\\
We design a structured individual-based model, and we study species characterized by a typical energy $x_{0}$. The energy $x_{0}$ is the energy of every new individual appearing in the population at a birth event. This specific reproduction mechanism makes our model closer to age-structured ones \cite{dieckmann22}, rather than growth-fragmentation patterns \cite{doumic2015statistical}, \cite{coranicofab}. We assume that the resource consumed by individuals is maintained at a fixed level $R$, thus neglecting individuals interactions through the resource consumption. Thanks to this assumption, we benefit from a \textit{branching property}. The whole system is ruled by two sorts of mechanisms: random jumps conforming to births and deaths in the population; and a continuous and deterministic evolution of individual energies between the jumps. Hence, our model gives rise to Piecewise Deterministic Markov Processes \cite{davis1984piecewise}. Our fundamental assumption is that all these mechanisms follow allometric scalings. Moreover, one of the new features of our model is that we will consider the broad range of possible characteristic energy $x_{0}$, to account for allometric properties valid from micro-organisms to mammals. Precisely, we highlight the necessary allometric scalings to obtain a biologically relevant model, allowing the survival of species (with a sufficient amount of resource), for every $x_{0}$. We will let $x_{0} \rightarrow 0$ and $+ \infty$ to infer mathematical constraints from our biological assumptions. Thus, our model intends to simulate a general interspecific structure, rather than a precise intraspecific explanation of allometries (where $x_{0}$ would be fixed once and for all). 
\\\\
In Section~\ref{sec:model}, we describe the general features of our model. At the individual level in Section~\ref{subsec:gensett}, we precise the main mechanisms involved: births, deaths, energy loss and gain. In this paper, individual trajectories $(\xi_{t})_{t \geq 0}$ are the main object of interest. In Section~\ref{subsec:cns}, we give a necessary and sufficient condition for these trajectories to be biologically relevant (this is Theorem~\ref{theo:cns}), and the extended generator (see Definition 5.2. in \cite{davis1984piecewise}) of the individual process under this condition. At the population level in Section~\ref{subsec:whole}, we gather individual trajectories into a measure-valued process $\mu$. In Section~\ref{subsec:generation}, we focus on a Galton-Watson process, describing the size of generations of $\mu$. In Section~\ref{subsec:allomsett}, we introduce allometric relationships in the model. In Section~\ref{sec:notations}, we present the main results of this article. First, assuming that the allometric coefficient for metabolism verifies $\alpha \leq 1$ (a lot of papers argue for $\alpha = 0.75$ \cite{peters_1983, brown_04, savdee_08}), we state our main theorem in Section~\ref{sec:main} (see Theorem~\ref{theo:short}). We also present our results when $\alpha > 1$ in Section~\ref{appendsix:swag}, where we obtain similar but weaker constraints on the allometric coefficients in Theorem~\ref{theo:sharpdeux}. In Section~\ref{sec:discussion}, we discuss the mathematical aspects and biological interpretations of our modelling and results, and give some perspectives. Section~\ref{proofs} is devoted to the technical intermediate results, leading to Theorem~\ref{theo:cns},~\ref{theo:short} and~\ref{theo:sharpdeux}. The proofs are based on three different ways to construct the individual process. The first one, given in Section~\ref{subsec:construc}, uses exponential random variables and is useful for Sections~\ref{subsec:delta} to~\ref{subsec:laprevue}. The second construction highlights a coupling between two processes, distinct by a different value for $x_{0}$, and is established in Section~\ref{proof:beplusde}. Also, we construct two other useful couplings between our process and similar ones with modified birth and death rates, in Section~\ref{subsec:laprevue}. The third construction, in Section~\ref{sec:pta}, uses a Poisson point measure, which brings useful martingale properties. In Section~\ref{subsec:simulations}, we present numerical simulations to illustrate our results and conjectures. In Appendix~\ref{app:construc}, we give some technical details about the definition of the population process of Section~\ref{subsec:whole}, and the embedded branching process of Section~\ref{subsec:generation}. In Appendix~\ref{appendix:pta}, we give details for the proof in Section~\ref{sec:pta}.

\section{Description of the model}
\label{sec:model}
\subsection{A general setting}
\label{subsec:gensette}

We design an individual-based model, structured by a positive trait called \textit{energy}. We study single species characterized by their energy at birth $x_{0} >0$, living in an environment described by a constant amount of resources $R$. First, we will describe in Section~\ref{subsec:gensett} the individual dynamics through the process $(\xi_{t})_{t \geq 0}$, which is the main process of interest in this article. We give in Theorem~\ref{theo:cns} of Section~\ref{subsec:cns} a necessary and sufficient condition for this process to be biologically relevant, and provide its extended generator under this condition in Proposition~\ref{prop:infinitesimal}. Then in Section~\ref{subsec:whole}, we gather individual trajectories into a population process. One standard way to look at it is to construct a measure-valued process \cite{four_04, champagnat2005individualbased}, such that at time $t \geq 0$, the population is described by a point measure $\mu_{t}$ on $\mathbb{R}^{*}_{+}$. In our context, this measure-valued process will benefit from a \textit{branching property}, which is the case for classical and similar approaches (see \cite{marguet2016} or Remark 2.2. in \cite{coranicofab}). Finally in Section~\ref{subsec:generation}, we define a process $(\Upsilon_{n})_{n \geq 0}$, representing the population size generation by generation in $(\mu_{t})_{t \geq 0}$, and such that $(\Upsilon_{n})_{n \geq 1}$ is a Galton-Watson process. This allows us to draw a link between a survival criterion for the population process $(\mu_{t})_{t \geq 0}$, and the mean number of birth jumps for the individual process $(\xi_{t})_{t \geq 0}$.

\subsubsection{Individual trajectories}
\label{subsec:gensett}

For $x_{0}, \xi_{0} >0$ and $R \geq 0$, we denote by $(\xi_{t,x_{0},R,\xi_{0}})_{t \geq 0}$ the process describing the evolution of energy over time, for an individual starting from energy $\xi_{0}$ at time 0, within a species characterized by $x_{0}$, with resources $R$. We simply write $(\xi_{t})_{t \geq 0}$ in the following if there is no ambiguity. An individual trajectory will be deterministic between some random jump times, corresponding to birth or death events. We define two cemetery states $\partial,\flat \notin \mathbb{R}$. If an individual energy reaches $\partial$ at time $t$, it means that the associated individual is dead and then $\xi_{s} = \partial$ for $s \geq t$. If an individual energy reaches $\flat$ at time $t$, it means that this energy exploded in finite time or reached the value 0 at time $t$, and then $\xi_{s} = \flat$ for $s \geq t$. The process $\xi$ takes its values in $\mathbb{R}^{*}_{+} \cup \{ \partial, \flat \}$. We precise now the general structure and parameters of our model. 
\begin{itemize}
\item \textbf{Birth}
\\ At time $t$, an individual with energy $\xi_{t-} \notin \{ \partial, \flat \}$ jumps from $\xi_{t-}$ to $\xi_{t} := \xi_{t-}-x_{0}$ at rate $b_{x_{0}}(\xi_{t-})$. This event corresponds to the transfer of a constant amount of energy $x_{0}$ to a single offspring born at time $t$.
\smallbreak
We assume that the birth rate $b_{x_{0}}$, defined on $\mathbb{R}^{*}_{+} $, is equal to 0 for every $x \leq x_{0}$, so that no individual with non-positive energy appears during a birth event. Also, we assume that it is equal to some positive and continuous function $\tilde{b}$ (not depending on $x_{0}$) for $x > x_{0}$. Finally, for $x >0$, $$b_{x_{0}}(x) =\mathbb{1}_{x > x_{0}}\tilde{b}(x).$$ 

\item \textbf{Death}
\\ At time $t$, every individual with energy $\xi_{t-} \notin \{ \partial, \flat \}$ dies at rate $d(\xi_{t-})$. The energy then jumps from $\xi_{t-}$ to $\xi_{t}:=\partial$, and we set $\xi_{s} := \partial$ for all $s \geq t$. The function $d$ is defined on $\mathbb{R}^{*}_{+} $, and assumed to be positive and continuous.
\item \textbf{Energy loss and resource consumption}
\begin{enumerate}
    \item Every individual with energy $\xi_{t} \notin \{ \partial, \flat \}$ loses energy over time at rate $\ell(\xi_{t})$. The function $\ell$ is defined on $\mathbb{R}^{*}_{+} $ and assumed to be positive, increasing and continuous. This is the function related to metabolism, that is the loss of energy for maintenance, in our model. The metabolism is commonly observed to increase when the mass of an individual increases \cite{wickman24}.
    \item In order to balance this energy loss, an individual with energy $\xi_{t} \notin \{ \partial, \flat \}$ consumes the resource at rate $f(\xi_{t},R)$. Importantly, we suppose that $f$ is of the form $f(\xi,R) := \phi(R)\psi(\xi),$ where $\phi$ is a continuous, non-negative, and increasing function on $\mathbb{R}^{+}$. Also, we assume that $\psi$ is a continuous, positive and increasing function on $\mathbb{R}^{*}_{+} $. 
\\\\
The fact that there is no resource consumption without resource imposes $\phi(0)=0$, and for a fixed energy $x$, we assume that the resource consumption is bounded even if there is an infinite amount of resource, which can always be reformulated as $\lim_{R \rightarrow + \infty}\phi(R)=1$. The shape of $\phi$ is typically the one for the functional response in evolutionary ecology: we can think of Holling type II or III functional responses \cite{yodzis_1992}. The function $\psi$ accounts for the conversion of resource into energy by the individual, which is why it is positive and increasing with energy, as for $\ell$. 
\end{enumerate}
For every $x>0$ and $R \geq 0$, we write $g(x,R) :=f(x,R) - \ell(x)$. We suppose that for every $R \geq 0$, $g(.,R)$ is $\mathcal{C}^{1}$ on $\mathbb{R}^{*}_{+}$, so in particular it is locally Lipschitz continuous. Between two random jump times (due to birth or death events), the energy evolves according to the following equation:
\begin{align}
    \dfrac{\mathrm{d}\xi_{t}}{\mathrm{d}t} =  g(\xi_{t},R).
    \label{eq:indivenergy}
\end{align}
The regularity of $g$ ensures that Equation~\eqref{eq:indivenergy} admits a unique positive local solution, starting from any positive energy $\xi_{0}$. This solution is denoted as $A_{\xi_{0},R} : t \mapsto A_{\xi_{0},R}(t)$ or simply $A_{\xi_{0}}(.)$ if there is no ambiguity on $R$. It is defined on an interval $[0,t_{\mathrm{max}}(\xi_{0},R)[$, where $t_{\mathrm{max}}(\xi_{0},R)$, or simply $t_{\mathrm{max}}(\xi_{0})$, is the deterministic time when it reaches 0 or $+ \infty$ ($t_{\mathrm{max}}(\xi_{0})$ is equal to $+ \infty$ if this never happens). If this happens between two random jumps at some time $t$, then the energy jumps from $\xi_{t-}$ to $ \xi_{t} := \flat$, and we set $\xi_{s} = \flat$ for all $s \geq t$. Classical arguments show that the flow $(\xi_{0},t) \mapsto A_{\xi_{0}}(t)$ is $\mathcal{C}^{1,1}$ on $[0,t_{\mathrm{max}}(\xi_{0})[$, meaning that it is differentiable with continuous derivatives in both variables $(\xi_{0},t)$.
\end{itemize}
\noindent In the following, we refer to all the previous assumptions for the model as the `general setting of Section~\ref{subsec:gensett}'. We define formally the process $(\xi_{t})_{t \geq 0}$ in Section~\ref{subsec:construc}, using an iterative construction of jump times. Also, we give constructions that are equivalent in law (in the sense that they have same distribution of sample paths) in Sections~\ref{proof:beplusde} and~\ref{sec:pta}. In the following, we denote by $\mathbb{P}_{x_{0}, R, \xi_{0}}$ the law associated to the individual process $(\xi_{t,x_{0},R,\xi_{0}})_{t \geq 0}$. If there is no ambiguity on some of these parameters, we will allow ourselves to lighten the notations into $\mathbb{P}_{\xi_{0}}$ for example. The associated expectations are naturally written as $\mathbb{E}_{x_{0}, R, \xi_{0}}$, or simply $\mathbb{E}$.

\subsubsection{Necessary and sufficient condition for a biologically relevant individual process}
\label{subsec:cns}

We consider the process $(\xi_{t,x_{0},R,\xi_{0}})_{t \geq 0}$ defined in Section~\ref{subsec:gensett}, and formally constructed in Section~\ref{subsec:construc}. We write respectively $T_{0,x_{0},R,\xi_{0}}$ and $T_{\infty,x_{0},R,\xi_{0}}$, or simply $T_{0}$ and $T_{\infty}$, for the random times when the process $(\xi_{t,x_{0},R,\xi_{0}})_{t \geq 0}$ reaches respectively 0 and $+ \infty$. Also, we denote by $T_{d,x_{0},R,\xi_{0}}$, or simply $T_{d}$, the random time when our process reaches $\partial$ (\textit{i.e.} the time of death of the individual). All these hitting times may be finite or infinite. We say that a trajectory $(\xi_{t})_{t \geq 0}$ is biologically relevant, if and only if the following event occurs:
$$ \{T_{d} < T_{0} \wedge T_{\infty} \}.$$ 
Remark that this event is equivalent to $ \{ T_{d}<+ \infty \}$, because if the time of death is finite, it necessarily means that the process did not reach 0 or $+ \infty$ before the death. In other terms, $(\xi_{t})_{t \geq 0}$ is not biologically relevant if it reaches $\flat$ before $\partial$, or if $(\xi_{t})_{t \geq 0}$ never reaches $\partial$ (\textit{i.e.} the individual never dies). Let us state a first situation where the trajectories are biologically relevant. For $R \geq 0$, we say that $R \in \mathfrak{R}_{0}$, if and only if
$$ \exists x >0, \forall y \in ]0,x], \quad g(y,R) < 0, $$
and we say that $R \in \mathfrak{R}_{\infty}$, if and only if
$$ \exists x >0, \forall y \geq x, \quad g(y,R) > 0. $$
Remark that under the general setting of Section~\ref{subsec:gensett}, we always have $0 \in \mathfrak{R}_{0}$, and both $\mathfrak{R}_{0}$ and $\mathfrak{R}_{\infty}$ are intervals when they are non-empty. We will see in Lemma~\ref{lemme:premierr} in Section~\ref{subsec:teun} that if $R \notin \mathfrak{R}_{0} \cup \mathfrak{R}_{\infty} $, then the process is almost surely biologically relevant for every $x_{0},\xi_{0}$. In the following, we make additional assumptions for our individual trajectories to be biologically relevant almost surely for the remaining values of $R$ (see Theorem~\ref{theo:cns} below). Let us introduce some notations for this purpose.
For $R\geq 0$ and $x_{0}>0$, we define the operator $K_{x_{0},R}$ such that for any bounded function $f : \mathbb{R}^{*}_{+} \rightarrow \mathbb{R}$,
\begin{align}
K_{x_{0},R}f : \xi_{0} \mapsto \displaystyle{\int_{0}^{t_{\max}(\xi_{0},R)}b_{x_{0}}(A_{\xi_{0},R}(u))e^{-\displaystyle{ \int_{0}^{u}}(b_{x_{0}}+d)(A_{\xi_{0},R}(\tau))\mathrm{d}\tau}f(A_{\xi_{0},R}(u)-x_{0}) \mathrm{d}u }.
\label{eq:kixzerooo}
\end{align} 
We write $\mathbf{1}$ for the constant function equal to 1. Remark that $K_{x_{0},R}$ is linear, positive and $K_{x_{0},R}\mathbf{1} \leq \mathbf{1}$. For any bounded function $f$ and $k \geq 1$, we write $K_{x_{0},R}^{k}f$ for the operator $K_{x_{0},R}$ applied to $K_{x_{0},R}^{k-1}f$, and $K_{x_{0},R}^{0}f=f$. Remark that $K_{x_{0},R}\mathbf{1}(\xi_{0})$ corresponds to the probability that, starting from energy $\xi_{0}$ with characteristic energy $x_{0}$ and resource $R$, the next jump of the process is a birth event (see Lemma~\ref{lemme:lemmatwo} of Section~\ref{subsec:teun} for a proof).
\begin{hyp}[\textbf{Individual energy avoids 0}]
For all $x_{0}>0$,
\begin{align}
\displaystyle{\int_{0}^{x_{0}} \dfrac{d}{\ell}(x) \mathrm{d}x }  = + \infty.
\label{eq:deltacahnge}
\end{align}
\label{hyp:probamort}
\end{hyp}
\begin{hyp}[\textbf{Individual energy avoids $+ \infty$}]
For all $x_{0}>0$, for all $R \in \mathfrak{R}_{\infty}$, for all $\xi_{0} >0$ such that $g(x,R) >0$ for $x \in [\xi_{0},+\infty[$,
\begin{align}
K_{x_{0},R}^{k}\mathbf{1}(\xi_{0}) \xrightarrow[k \rightarrow + \infty]{} 0 \hspace{0.3cm} \mathrm{and} \hspace{0.3 cm} \displaystyle{\int_{\xi_{0}}^{+ \infty} \dfrac{(b_{x_{0}}+d)(x)}{g(x,R)} \mathrm{d}x }  = + \infty.
\label{eq:beplusdechang}
\end{align}
\label{hyp:tempsinf}
\end{hyp}
\noindent \textbf{Remark:} If $R$ and $\xi_{0}$ are chosen like in Assumption~\ref{hyp:tempsinf}, we can use the change of variables $x=A_{\xi_{0},R}(u)$ and write for any bounded function $f : \mathbb{R}^{*}_{+} \rightarrow \mathbb{R}$ and such $\xi_{0},R$,
\[ K_{x_{0},R}f(\xi_{0}) = \displaystyle{\int_{\xi_{0}}^{+ \infty}\dfrac{b_{x_{0}}(x)}{g(x,R)} e^{-\displaystyle{ \int_{\xi_{0}}^{x}}\frac{(b_{x_{0}}+d)(\tau)}{g(\tau,R)}\mathrm{d}\tau}f(x-x_{0}) \mathrm{d}x }.\]
\noindent We precise in the following theorem how these assumptions lead to a process $\xi$ with almost surely biologically relevant paths. 
\begin{theorem}
Under the general setting of Section~\ref{subsec:gensett}, we have
\begin{align*}
Assumptions \hspace{0.1 cm} \textit{\ref{hyp:probamort}} \hspace{0.1 cm} and \hspace{0.1cm} \textit{\ref{hyp:tempsinf}} & \Leftrightarrow (\forall x_{0} >0, \forall \xi_{0}>0, \forall R \geq 0, \quad \mathbb{P}_{x_{0}, R, \xi_{0}}(T_{d} < + \infty) = 1).
\end{align*} 
\label{theo:cns}
\end{theorem}
\noindent \textbf{Remark:} Theorem~\ref{theo:cns} states that Assumptions~\ref{hyp:probamort} and~\ref{hyp:tempsinf} are necessary and sufficient to obtain a biologically relevant individual process for every $R \geq 0$ in our setting. More precisely, Assumption~\ref{hyp:probamort} almost surely prevents individual energy from reaching 0 (we only need a condition in the worst possible case, that is $R=0$). We give an equivalent formulation in Proposition~\ref{eq:delta} of Section~\ref{subsec:delta}. Assumption~\ref{hyp:tempsinf} almost surely prevents individual energy from reaching $+ \infty$ (in the case where it would be possible, that is $R \in \mathfrak{R}_{\infty}$). We give an equivalent formulation, under Assumption~\ref{hyp:probamort}, in Proposition~\ref{prop:tinfyinfini} of Section~\ref{subsec:teun}. To prove Theorem~\ref{theo:cns}, we will remark at the end of Section~\ref{subsec:teun} that it is the combination of Proposition~\ref{eq:delta}, Lemma~\ref{lemme:premierr} and Proposition~\ref{prop:tinfyinfini}.
\\\\
Finally, we provide the extended generator (see Definition 5.2. in \cite{davis1984piecewise}) of our process. We denote by $\mathcal{C}^{1}(\mathbb{R}^{*}_{+} \cup \{ \partial \},\mathbb{R})$ the space of functions in $\mathcal{C}^{1}(\mathbb{R}^{*}_{+},\mathbb{R})$ that take some real value at $\partial$. We define the operator $\mathcal{L}$ such that for every $\varphi \in \mathcal{C}^{1}(\mathbb{R}^{*}_{+}\cup\{\partial\},\mathbb{R})$,
\begin{align}
\mathcal{L}\varphi : x \mapsto \left\{
    \begin{array}{cl}
       g(x,R)\varphi'(x) + b_{x_{0}}(x)\left( \varphi(x-x_{0}) - \varphi(x) \right) +d(x)(\varphi(\partial) -\varphi(x)) & \mathrm{if} \hspace{0.1 cm} x \neq \partial, \\
        \varphi(\partial) & \mathrm{if} \hspace{0.1 cm} x= \partial.
    \end{array}
\right.
\label{eq:dybnkin}
\end{align}
We write $\mathcal{D}$ for the subset of $\mathcal{C}^{1}(\mathbb{R}^{*}_{+} \cup \{ \partial \},\mathbb{R})$, such that $\varphi \in \mathcal{D}$, if and only if $\varphi$ and $\mathcal{L}\varphi$ are bounded on $\mathbb{R}^{*}_{+} \cup \{ \partial \}$. Under Assumptions~\ref{hyp:probamort} and~\ref{hyp:tempsinf}, thanks to Theorem~\ref{theo:cns}, the process $(\xi_{t})_{t \geq 0}$ almost surely avoids $\flat$, so we can almost surely write $\varphi(\xi_{t})$ and $\mathcal{L}\varphi(\xi_{t})$ for every $t \geq 0$ and $\varphi \in \mathcal{C}^{1}(\mathbb{R}^{*}_{+}\cup\{\partial\},\mathbb{R})$. 
\begin{prop}
Under the general setting of Section~\ref{subsec:gensett}, and under Assumptions~\ref{hyp:probamort} and~\ref{hyp:tempsinf}, the individual process $(\xi_{t,x_{0},R,\xi_{0}})_{t \geq 0}$ (formally defined in Section~\ref{subsec:construc}) is a Piecewise Deterministic Markov Process (PDMP) with extended generator $\mathcal{L}$, where the domain $D(\mathcal{L})$ of this generator contains $\mathcal{D}$.
\\\\
In particular, it means that for every $\varphi \in \mathcal{D}$, the process $(\mathfrak{M}_{t})_{t \geq 0}$ defined for $t \geq 0$ by
$$  \mathfrak{M}_{t} := \varphi(\xi_{t}) - \varphi(\xi_{0}) - \displaystyle{\int_{0}^{t} \mathcal{L}\varphi(\xi_{s}) \mathrm{d}s}$$
is a martingale.
\label{prop:infinitesimal}
\end{prop}

\begin{proof}
We refer to the formalization of PDMPs by Davis in \cite{davis1984piecewise}. Under Assumptions~\ref{hyp:probamort} and~\ref{hyp:tempsinf}, our process almost surely does not explode or touch 0 in finite time. Furthermore, Equation~\eqref{eq:indivenergy} admits a unique positive local solution, starting from any positive energy $\xi_{0}$ thanks to the regularity of $g$. Jump rates $b_{x_{0}}$ and $d$ are $L^{1}_{\mathrm{loc}}(\mathbb{R}^{*}_{+})$, hence classical techniques using Gronwall inequality show that we verify Assumption 3.1. of \cite{davis1984piecewise}. Thus, our result is a consequence of Theorem 5.5. in \cite{davis1984piecewise}, which characterizes the full domain $D(\mathcal{L})$. The fact that $\mathcal{D} \subseteq D(\mathcal{L})$ can be checked using classical arguments. Indeed, Theorem 5.5. in \cite{davis1984piecewise} immediately entails that $\mathcal{C}^{\infty}$ functions with compact support are in $D(\mathcal{L})$, and we then use localisation and approximation by these test functions.
\end{proof}
\noindent \textbf{Remark:} In Theorem 5.5. of \cite{davis1984piecewise}, the domain $D(\mathcal{L})$ is entirely characterized. We only presented here a partial result on $\mathcal{D}$, for the sake of simplicity.

\subsubsection{The measure-valued population process}
\label{subsec:whole}
In the previous section, we described an individual process $\xi$ for which we can interpret random jumps as birth or death events. Our goal here is to define a measure-valued population process $(\mu_{t})_{t}$ based on these individual trajectories. In the following, we fix $x_{0}>0$, $R \geq 0$, and we work under the general setting of Section~\ref{subsec:gensett}. To enumerate individuals in the population, we use the Ulam-Harris-Neveu notation. We define
$$ \mathcal{U} := \bigcup_{n \in \mathbb{N}}(\mathbb{N}^{*})^{n+1}.$$
Over time, every individual will have an index of the form $ u := u_{0}...u_{n}$ with some $n \geq 0$, and some positive integers $u_{0},..., u_{n}$, \textit{i.e.} some $u \in \mathcal{U}$. The \textit{generation} of $u$ is $|u| := n$. We denote as $V_{t}$ the set of indices of individuals alive at time $t$ (\textit{i.e.} individuals born before $t$, and whose energy is not in $\{ \partial, \flat \}$ at time $t$). The size of the population at time $t$ is then naturally $\mathrm{Card}(V_{t})$.
\\\\
Let $\mathcal{M}_{P}$ be the set of finite point measures on $\mathbb{R}^{*}_{+}$. At time $t=0$, we pick a random variable $\mu_{0} \in \mathcal{M}_{P}$. Then, there exists a random variable $\mathcal{C}_{0} \in \mathbb{N}$ and a random vector $(\xi_{0}^{1},..., \xi_{0}^{\mathcal{C}_{0}}) \in (\mathbb{R}^{*}_{+})^{\mathcal{C}_{0}}$ such that
\begin{center}
        $\mu_{0} = \sum\limits_{u \in V_{0}}{\delta_{\xi_{0}^{u}}},$
\end{center}
with $V_{0} := \{1,...,\mathcal{C}_{0}\}$. Initial individuals have single labels $ u \in V_{0}$, their trajectories are independent conditionally to $\mu_{0}$, distributed as the process $\xi$ of Section~\ref{subsec:gensett} started from $\xi_{0}^{u}$ at time $\tau_{u} := 0$, and are denoted by $(\xi_{t}^{u,\tau_{u}})_{t \geq 0}$.
\\\\
Then, we define a family of i.i.d. processes $(\xi^{u})_{u \in \mathcal{U} \setminus \mathbb{N}^{*}}$, such that for every $u \in \mathcal{U} \setminus \mathbb{N}^{*}$, $\xi^{u}$ is distributed as the process $\xi$ of Section~\ref{subsec:gensett} started from $x_{0}$ (which corresponds to the fixed amount of energy transferred to offspring). Also, these processes are taken independent from the initial processes $(\xi^{u,0})_{u \in V_{0}}$. We also define, for every $\tau > 0$ and $u \in \mathcal{U} \setminus \mathbb{N}^{*} $, the shifted process $(\xi^{u,\tau}_{t})_{t \geq \tau}$ with $\xi^{u,\tau}_{t} := \xi^{u}_{t-\tau}$ for every $t \geq \tau$. When a birth event occurs at some random time $\tau$ among one of the individuals in $V_{\tau-}$, we create a new individual trajectory. We index the offspring in the following manner: if the parent has the label $u:=u_{0}u_{1}...u_{n}$ and has already given birth to $k$ children for $k \geq 0$, the index of the new offspring is taken as $uk:=u_{0}u_{1}...u_{n}u_{n+1}$ with $u_{n+1} = k$, and we add the index $uk$ to $V_{\tau}$. Also, we define $\tau_{uk} := \tau$ and the offspring energy follows the process $\xi^{uk,\tau_{uk}}$. When an individual dies or if its energy reaches $\flat$ at time $\tau,$ we remove this individual from $V_{\tau}$. In this manner, over time, we describe the population with a point measure written as
\begin{center}
        $\mu_{t} = \sum\limits_{u \in V_{t}}{\delta_{\xi_{t}^{u,\tau_{u}}}}$,
\end{center}
with $V_{t} \subseteq \mathcal{U}$. Since it relies on the construction of $\xi$ in Section~\ref{subsec:construc}, we do not further detail here the construction of $\mu$ and postpone it to Appendix~\ref{app:mesure}. In particular, $\mu_{t}$ is well-defined for $t \in [0,\overline{\Theta}[$, where $\overline{\Theta}< + \infty$ if there is an accumulation of jump times (and in that case, $\overline{\Theta}$ is the supremum of the jump times), and $\overline{\Theta}= + \infty$ otherwise. We prove in Appendix~\ref{app:mesure}, Proposition~\ref{prop:final}, that under Assumptions~\ref{hyp:probamort} and \ref{hyp:tempsinf}, then $\overline{\Theta} = + \infty$ almost surely for every random variable $\mu_{0} \in \mathcal{M}_{P}$. The study of such measure-valued processes is standard, the reader can refer to \cite{chi_08} for an age-structured measure-valued process; and for growth-fragmentation models, to \cite{coranicofab} in an evolutionary and ecological context, or to \cite{tomavsevic2022ergodic} for the spatial spreading of a filamentous fungus. 
\\\\
Finally, we denote as $\mathbb{Q}_{\mu_{0},x_{0}, R}$ the law associated to the population process $\mu$, with initial condition $\mu_{0}$, offspring energy $x_{0}$ and resource $R$. We define the event $\mathfrak{B} := \{ \exists t \geq 0, \exists u \in V_{t}, u \notin V_{0} \}$, \textit{i.e.} there is at least one birth event in the population. For any $x_{0}>0$, $R \geq 0$, we also define $\mathcal{M}_{P,x_{0},R} := \{ \mu_{0} \in \mathcal{M}_{P}, \mathbb{Q}_{\mu_{0},x_{0}, R}(\mathfrak{B}) >0 \}$. Remark that $\mathcal{M}_{P,x_{0},R} \neq \varnothing$, because $\delta_{2x_{0}} \in \mathcal{M}_{P,x_{0},R}$ for example. We will investigate in this article under which conditions on the model parameters the following assumption holds true.
\begin{hyp}[\textbf{Supercriticality of the population process}] For every $x_{0} >0$, there exists $R_{0} >0$, such that for all $R>R_{0}$, for every $\mu_{0} \in \mathcal{M}_{P,x_{0},R}$,
$$\left(\mathbb{Q}_{\mu_{0},x_{0},R}(\overline{\Theta} = + \infty)=1 \right) \hspace{0.1 cm} and \hspace{0.1 cm} \left(\mathbb{Q}_{\mu_{0},x_{0},R}(\forall t \geq 0, \hspace{0.1 cm} V_{t} \neq \varnothing ) >0 \right).$$
\label{ass:supercritical}
\end{hyp}
The first condition $\mathbb{Q}_{\mu_{0},x_{0},R}(\overline{\Theta} = + \infty)=1$ ensures that the population process is almost surely well-defined for all $t \geq 0$. On this event, Assumption~\ref{ass:supercritical} states that for any studied species, there exists an amount of resources from which the population can survive indefinitely with positive probability. Assumption~\ref{ass:supercritical} could be verified if at least one of the individuals in the population never dies. Thanks to Theorem~\ref{theo:cns}, we will not consider this case if we work under Assumptions \ref{hyp:probamort} and \ref{hyp:tempsinf}. Furthermore, Proposition~\ref{prop:final} in Appendix~\ref{app:mesure} implies that under Assumptions~\ref{hyp:probamort} and \ref{hyp:tempsinf}, $\mathbb{Q}_{\mu_{0},x_{0},R}(\overline{\Theta} = + \infty)=1$ for every $\mu_{0}, x_{0},R$. We want to work under Assumption~\ref{ass:supercritical}, to model the fact that living conditions can be enforced in laboratory studies, which allow populations of bacteria or algae to last as long as we want with high probability, if we maintain a sufficient level of resources \cite{malerba_2019}. We insist here on the fact that it is this particular choice of a fixed value for $R$ that fully supports our assumption. If there were competitive dynamics on $R$, the birth and death process could lead to almost sure extinction, and the best we can do in this case is to study existence and/or uniqueness of quasi-stationary distributions \cite{collet2011quasi}. 

\subsubsection{The generation process}
\label{subsec:generation}

Finally, we define a process $(\Upsilon_{n})_{n \geq 0}$ describing the size of generations of $\mu$, taking values in $\mathbb{N} \cup \{+ \infty\}$. We refer the reader to Corollary 2 in \cite{doney_72} for a general description of the embedding of a generation process into a continuous in time population process. For every $n \geq 0$, we define the random set
$$G_{n} := \{ u \in \mathcal{U}, |u|=n, \exists t \geq 0, u \in V_{t} \}, $$
which contains all the individuals of the $n$-th generation and let
$$ \Upsilon_{n} := \mathrm{Card}(G_{n}).$$ 
With this definition, $\Upsilon$ could \textit{a priori} be infinite at some point, but as soon as it should represent the generation sizes of a population, we want $\Upsilon_{n}$ to be almost surely finite for every $n \geq 0$, so we need the following assumption. For every $x_{0}, \xi_{0}>0$ and $R \geq 0$, we write $N_{x_{0},R,\xi_{0}}$ for the number of direct offspring of an individual following the process $(\xi_{t,x_{0},R,\xi_{0}})_{t \geq 0}$ during its life (\textit{i.e.} the number of birth events on its trajectory). We write $\nu_{x_{0},R,\xi_{0}}$ for the law of $N_{x_{0},R,\xi_{0}}$.  
\begin{hyp}
$\forall x_{0}>0, \forall R \geq 0, \forall \xi_{0}>0, \quad \mathbb{P}(N_{x_{0},R, \xi_{0}} < + \infty) = 1.$
\label{ass:pasinfini}
\end{hyp}
\noindent We will see in Corollary~\ref{corr:swagosss} of Section~\ref{subsec:teun} that Assumption~\ref{hyp:tempsinf} implies Assumption~\ref{ass:pasinfini}.
\begin{prop}
Under Assumption~\ref{ass:pasinfini}, $(\Upsilon_{n})_{n \geq 1}$ is a Galton-Watson process with offspring distribution $\nu_{x_{0},R,x_{0}}$.
\label{prop:jesuisungalton}
\end{prop}
The proof of this fact is developed in Appendix~\ref{appendix:embedding}. The Galton-Watson process $(\Upsilon_{n})_{n \geq 1}$ is either subcritical, critical, or supercritical. We use classical results on Galton-Watson processes to infer the following equivalent formulation of Assumption~\ref{ass:supercritical}. Let us write $m_{x_{0},R}(\xi_{0}) := \mathbb{E}(N_{x_{0},R,\xi_{0}})$, or simply $m_{x_{0}}(\xi_{0})$.

\begin{prop}
Under the general setting of Section~\ref{subsec:gensett},  under Assumptions~\ref{hyp:probamort} and \ref{hyp:tempsinf}, we have that
\begin{align*}
Assumption \hspace{0.1 cm} \textit{\ref{ass:supercritical}} \hspace{0.1 cm} \Leftrightarrow (\forall x_{0} >0, \exists R_{0} >0, \forall R>R_{0}, \quad m_{x_{0},R}(x_{0})>1).
\end{align*}
\label{prop:galtonwatson}
\end{prop}
The proof of Proposition~\ref{prop:galtonwatson} is given in Appendix~\ref{appendix:embedding}. Finally, we state an additional assumption.
\begin{hyp}
$\forall x_{0} > 0, \exists R_{0}>0, \forall R > R_{0},  \quad g(x_{0},R)>0 .$
\label{ass:gainenergy}
\end{hyp}
\noindent In Assumption~\ref{ass:gainenergy}, we assume that for every species characterized by $x_{0}$, there exists an amount of resource $R_{0}$ that allows individuals to grow. This assumption seems very natural if we think of our model as a general model able to allow the survival of any species within a broad range of characteristic energy $x_{0}$. We do not want that for some $x_{0}$, individuals can not grow, no matter the available resources. Thanks to the definition of the generation process $\Upsilon$, we can show that this last assumption is an immediate consequence of Assumptions~\ref{hyp:probamort}, \ref{hyp:tempsinf} and \ref{ass:supercritical}.

\begin{lemme}
Under the general setting of Section~\ref{subsec:gensett},
\begin{center}
$($Assumptions~\ref{hyp:probamort}, \ref{hyp:tempsinf} and \ref{ass:supercritical}\hspace{0.05 cm}$)$ $\Rightarrow$ Assumption~\ref{ass:gainenergy}.
\end{center}
\label{lemm:ass}
\end{lemme}
\begin{proof}
Assume that Assumptions~\ref{hyp:probamort} and \ref{hyp:tempsinf} hold true, and that Assumption~\ref{ass:gainenergy} is not verified. Considering Equation~\eqref{eq:indivenergy}, there exists $x_{0}$ such that no individual starting with this characteristic energy can reach higher energies, which leads to a birth rate equal to 0. Then, for any $R \geq 0$, $m_{x_{0},R}(x_{0})= 0$ and Assumption~\ref{ass:supercritical} is not verified by Proposition~\ref{prop:galtonwatson}. 
\end{proof}

\noindent Assumptions~\ref{hyp:probamort},~\ref{hyp:tempsinf} and~\ref{ass:supercritical} are the main assumptions of this paper. They imply Assumptions~\ref{ass:pasinfini} and~\ref{ass:gainenergy}, and they allow us to model biologically relevant processes, both at the individual and population levels. We will now define an allometric setting, and our motivation in the rest of this article will be to determine the allometric coefficients verifying Assumptions~\ref{hyp:probamort},~\ref{hyp:tempsinf} and~\ref{ass:supercritical}.

\subsection{The allometric setting}
\label{subsec:allomsett}
\noindent In this section, we highlight a specific choice for the functions $\tilde{b}$, $d$, $\ell$, $f$. The parameters $R$ and $x_{0}$ remain general constants, but we assume here allometric shapes on the four remaining functional parameters. More precisely, we impose:
\begin{enumerate}
\item $\tilde{b}(x):= C_{\beta}x^{\beta}$, so that $b_{x_{0}}(x):= \mathbb{1}_{x > x_{0}}C_{\beta}x^{\beta}$,
\item $d(x):= C_{\delta}x^{\delta}$,
\item $\ell(x) :=C_{\alpha}x^{\alpha}$,
\item $f(x,R) := \phi(R)C_{\gamma}x^{\gamma}$ (\textit{i.e.} $\psi(x) = C_{\gamma}x^{\gamma}$),
\end{enumerate}
with $(\beta, \delta) \in \mathbb{R}$ and $(C_{\beta}, C_{\delta}, C_{\alpha}, C_{\gamma}, \alpha, \gamma) \in \mathbb{R}^{*}_{+}$. Recall that $\ell$ and $\psi$ are increasing, which is why $\alpha$ and $\gamma$ are positive. In the following, we refer to all these specific assumptions for the model as the `allometric setting of Section~\ref{subsec:allomsett}'.

\section{Main results}
\label{sec:notations}

Now that our model is fully described, let us present the main results of this paper. First, in Section~\ref{sec:main}, we work under the allometric setting of Section~\ref{subsec:allomsett} with $0 < \alpha \leq 1$, and give in Theorem~\ref{theo:short} necessary conditions on the allometric coefficients $\beta$, $\delta$, $\gamma$ in order to verify Assumptions~\ref{hyp:probamort},~\ref{hyp:tempsinf} and~\ref{ass:supercritical}. Recall that these assumptions express conditions that need to be valid for every $x_{0} >0$. We give a partial result for the converse implication in Proposition~\ref{prop:theoreciprocal}. Finally, we present our conjectures about a necessary and sufficient condition on the allometric coefficients to verify Assumptions~\ref{hyp:probamort},~\ref{hyp:tempsinf} and~\ref{ass:supercritical}. We also present our results in the case $\alpha >1$ in Section~\ref{appendsix:swag}. The details of the proofs are given in Section~\ref{proofs}. Even if our aim is to have a biological interpretation for precise allometric relationships, most of our intermediate results in Section~\ref{proofs} are presented and valid under the general setting of Section~\ref{subsec:gensett}.

\subsection{Allometric constraints in the case $\alpha \leq 1$}
\label{sec:main}

\begin{figure}[h!]
\centering
\begin{tikzpicture}
\draw[-,color=red,thick] (-0.5,0) -- (1.5,0);

\draw (2.5,0.1) -- (2.5,-0.1) node[below]{0};
\draw (5.2,0.1) -- (5.2,-0.1) node[below]{1};
\draw[dotted] (3.3,0) -- (3.3,1);
\node[above] at (3.3,1) {$\alpha-1+\frac{C_{\delta}}{C_{\gamma}-C_{\alpha}}$};
\node[below=0.1 cm] at (1.5,0) {$\alpha-1$};
\draw (4.2,0.1) -- (4.2,-0.1) node[below=0.1 cm]{$\alpha$};
\node[above left,color=green] at (1.5,0) {$I_{1}$};
\node at (7,-0.25) {$\beta$};
\node[above left, color=green] at (6,0) {$I_{2}$};
\node at (-0.25,1) {$\delta= \alpha -1$};
\node at (-0.55,0.5) {$\gamma= \alpha$};

\draw[color=red,thick,-] (1.5,0) -- (3.3,0);
\draw[color=green,thick,->] (3.3,0) -- (6.9,0);

\draw[fill=green] (1.5,0) circle (0.1cm);
\end{tikzpicture}
\caption{Visual representation of Theorem~\ref{theo:short}. The green dot represents the singleton $I_{1}$, and the green line represents the set $I_{2}$. The non-admissible allometric coefficients $\beta$ for our model with $\alpha \leq 1$ are highlighted in red (and also, every $(\gamma,\delta,\beta)$ with $\delta \neq \alpha -1$ or $\gamma \neq \alpha$ is non-admissible).}
\label{fig:coeff}
\end{figure}
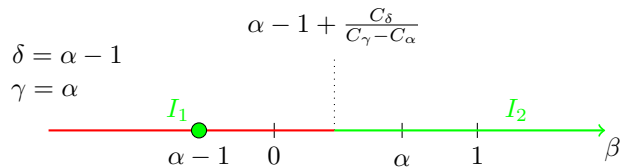

For some $0 < \alpha \leq 1$, we introduce two sets of allometric coefficients. The first one is a singleton, accounting for the usual allometric relationships highlighted in~\eqref{eq:equalallo}:
$$I_{1} := \{ (\gamma=\alpha, \delta= \alpha-1, \beta=\alpha-1) \}. $$
The second one contains a more complex condition on the allometric coefficient $\beta$ of the birth rate, and is well-defined if $C_{\gamma}> C_{\alpha}$:
$$I_{2} := \{ (\gamma=\alpha, \delta= \alpha-1, \beta), \hspace{0.1 cm} \beta \geq \alpha-1 + C_{\delta}/(C_{\gamma}-C_{\alpha}) \}. $$
\begin{theorem}
Let $0 <\alpha \leq 1$. Under the allometric setting of Section~\ref{subsec:allomsett}, under Assumptions~\ref{hyp:probamort} and~\ref{hyp:tempsinf} (biologically relevant model for every $x_{0}>0$) and Assumption~\ref{ass:supercritical} (supercriticality for every $x_{0}>0$), we have either
\begin{align}
(\gamma,\delta,\beta) \in I_{1}, \quad C_{\gamma} > C_{\alpha}, \quad C_{\beta} > C_{\delta}
\label{eq:iuncase}
\end{align}
or
\begin{align}
(\gamma,\delta,\beta) \in I_{2}, \quad C_{\gamma} > C_{\alpha}, \quad C_{\delta} \leq C_{\gamma}-C_{\alpha}.
\label{eq:itwocase}
\end{align}
\label{theo:short}
\end{theorem}
\noindent \textbf{Remark:} In Section~\ref{theo:sharp}, we give a detailed list of the constraints on the allometric coefficients implied by Assumptions~\ref{hyp:probamort},~\ref{hyp:tempsinf} and~\ref{ass:supercritical}, proved in the rest of Section~\ref{proofs}. These constraints together lead to Theorem~\ref{theo:short}. The reader can visualize the restrictions on allometric coefficients thanks to Figure~\ref{fig:coeff}. Note that the conditions on $C_{\gamma}$, $C_{\alpha}$, $C_{\beta}$ and $C_{\delta}$ are not represented on this graph. In the following, when the allometric coefficients are in $I_{1}$ (respectively $I_{2}$), we refer to it as `the $I_{1}$ (respectively $I_{2}$) case'.

An immediate consequence of Theorem~\ref{theo:short} is the following.
\begin{corr}
Under the allometric setting of Section~\ref{subsec:allomsett} and Assumptions~\ref{hyp:probamort},~\ref{hyp:tempsinf} and~\ref{ass:supercritical}, if $C_{\delta} > C_{\gamma} - C_{\alpha}$, then $$(\gamma,\delta,\beta) \in I_{1}, \quad C_{\gamma} > C_{\alpha}, \quad C_{\beta} > C_{\delta}.$$
\end{corr}
\noindent Hence, we recover in this specific case the usual allometric relationships highlighted by the Metabolic Theory of Ecology \cite{malerba_2019}, but it is supported by clear biological assumptions. We discuss more extensively the consequences of Theorem~\ref{theo:short} in Section~\ref{sec:discussion}. Finally, we give a partial result for the converse of Theorem~\ref{theo:short}.
\begin{prop}
Under the allometric setting of Section~\ref{subsec:allomsett}, if 
\begin{align}
(\gamma,\delta,\beta) \in I_{1}, \hspace{0.1 cm} C_{\gamma} > C_{\alpha}, \hspace{0.1 cm} C_{\beta} > (e-1)C_{\delta}, \hspace{0.1 cm} C_{\beta}+C_{\delta} < C_{\gamma} - C_{\alpha},
\label{eq:recipro}
\end{align} 
then Assumptions~\ref{hyp:probamort},~\ref{hyp:tempsinf} and~\ref{ass:supercritical} hold true.
\label{prop:theoreciprocal}
\end{prop}

\noindent We give the proof of this fact in Section~\ref{subsec:reciprocal}. Remark that~\eqref{eq:recipro} $\Rightarrow$~\eqref{eq:iuncase}, but we will show in Proposition~\ref{prop:conversenot} of Section~\ref{subsec:reciprocal} that~\eqref{eq:iuncase} is not sufficient to obtain Assumption~\ref{ass:supercritical}. Also, it is still an open question to know if~\eqref{eq:itwocase} is sufficient to verify Assumption~\ref{ass:supercritical}. We present detailed numerical simulations in Section~\ref{subsec:simulations}, that lead to the following conjectures about a necessary and sufficient condition to verify Assumptions~\ref{hyp:probamort},~\ref{hyp:tempsinf} and~\ref{ass:supercritical}, and the behavior of $m_{x_{0},R}(\xi_{0})$, in the $I_{2}$ case.
 
\begin{conj}
Let $0 <\alpha \leq 1$. Under the allometric setting of Section~\ref{subsec:allomsett}, we have
\begin{center}
$($Assumptions~\ref{hyp:probamort},~\ref{hyp:tempsinf} and~\ref{ass:supercritical}$)$ $\Leftrightarrow ($\eqref{eq:itwocase} $\mathrm{or}$~\eqref{eq:firstcondition}$)$,
\end{center}
where~\eqref{eq:itwocase} is already presented in Theorem~\ref{theo:short} (it is the $I_{2}$ case), and the second condition is
\begin{align}
(\gamma,\delta,\beta) \in I_{1}, \quad C_{\gamma} > C_{\alpha}, \quad C_{\delta} < C_{\gamma} - C_{\alpha}, \quad C_{\beta} > \Xi\left(\frac{C_{\delta}}{C_{\gamma} - C_{\alpha}}\right) C_{\delta},
\label{eq:firstcondition}
\end{align}
with $\Xi : ]0,1[ \rightarrow ]1, + \infty[$ a convex increasing function. 
\label{conj}
\end{conj}

\begin{conj}
Let $0 <\alpha \leq 1$. Under the allometric setting of Section~\ref{subsec:allomsett}, if \eqref{eq:itwocase} (in the $I_{2}$ case),
$$ \forall x_{0}>0, \forall R \geq 0, \forall \xi_{0}>0, \quad m_{x_{0},R}(\xi_{0}) = + \infty. $$
\label{conj:deux}
\end{conj}
\textbf{Remark:} Our last conjecture expresses that in the $I_{2}$ case, the average number of offspring starting from any $x_{0},R,\xi_{0}$ is infinite, which is indeed sufficient to verify Assumption~\ref{ass:supercritical} (supercriticality of the population process). This is very different from the $I_{1}$ case, where we can show that $m_{x_{0},R}(\xi_{0})$ is always finite (see for example the proof of Proposition~\ref{corr:pointquatre}).

\subsection{Allometric constraints in the case $\alpha > 1$}
\label{appendsix:swag}
We can adapt most of our reasoning from the case $\alpha  \leq 1$, and obtain the following theorem.

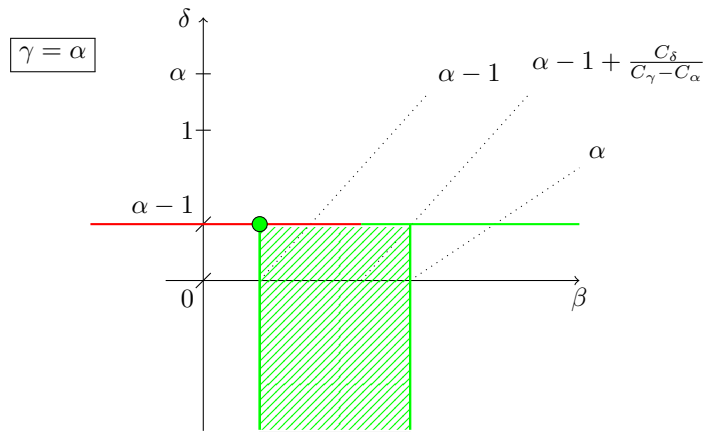
\begin{figure}[h!]
\centering
\begin{tikzpicture}
\draw[->] (-0.5,-2) -- (5,-2);
\draw[->] (0,-4) -- (0,1.5);

\draw (-0.1,-2.1) -- (0.1,-1.9) node[below left=0.1 cm]{0};
\draw (-0.1,0) -- (0.1,0) node[left=0.1 cm]{1};
\draw[dotted] (2.75,-2) -- (5,-0.5) node[above right] at (5,-0.5) {$\alpha$};
\draw[dotted] (2.1,-2) -- (4.35,0.5) node[above right] at (4.25,0.5){$\alpha-1+\frac{C_{\delta}}{C_{\gamma}-C_{\alpha}}$};
\draw[dotted] (0.75,-2) -- (3,0.5);
\node[above right] at (3,0.5) {$\alpha-1$};
\draw (-0.1,-1.35) -- (0.1,-1.15);
\node[above left] at (0,-1.25) {$\alpha-1$};
\draw (0.1,0.75) -- (-0.1,0.75) node[left]{$\alpha$};
\node at (5,-2.25) {$\beta$};
\node at (-0.25,1.5) {$\delta$};
\node at (-2,1) {$\boxed{\gamma=\alpha}$};

\draw[color=red,thick,-] (-1.5,-1.25) -- (2.1,-1.25);
\draw[color=green,thick,-] (2.1,-1.25) -- (5,-1.25);
\draw[color=green, thick,-] (2.75,-1.25) -- (2.75,-4);

\draw[color=green,thick,-] (0.75,-1.25) -- (0.75,-4);

\draw[pattern=north east lines, pattern color=green] (0.75,-1.3) -- (0.75,-4) -- (2.75,-4)--(2.75,-1.3);

\draw[color=green, thick,-] (2.75,-1.25) -- (2.75,-4);

\draw[color=green,thick,-] (0.75,-1.25) -- (0.75,-4);

\draw[color=white,thick,-] (0.7,-4) -- (2.8,-4);

\draw[fill=green] (0.75,-1.25) circle (0.1cm);
\end{tikzpicture}
\caption{Visual representation of Theorem~\ref{theo:sharpdeux}. All the coefficients $(\gamma, \delta, \beta)$ verifying one of the following conditions are non-admissible: 1) $\gamma \neq \alpha$; 2) $\delta \neq \alpha -1$, except for the ones in the hatched area and its border; 3) $\gamma=\alpha$, $\delta= \alpha-1$ and $\beta$ on the red line (except for the green dot).}
\label{fig:alphastrictun}
\end{figure}

\begin{theorem}
Let $\alpha > 1$. Under the allometric setting of Section~\ref{subsec:allomsett}, under Assumptions~\ref{hyp:probamort},~\ref{hyp:tempsinf} and~\ref{ass:supercritical}, we have:
\begin{itemize}
\item[\textbf{1.}] $\gamma = \alpha, \hspace{0.1 cm} C_{\gamma} > C_{\alpha}.$
\item[\textbf{2.}] $\delta \leq \alpha - 1 \leq \beta$.
\item[\textbf{3.}] If $\beta=\delta=\alpha-1$, then $C_{\beta} > C_{\delta}$.
\item[\textbf{4.}] If $(\delta=\alpha-1 \hspace{0.1cm} \mathrm{and} \hspace{0.1cm} \beta > \alpha-1)$, then $\beta \geq \alpha -1 + \dfrac{C_{\delta}}{C_{\gamma}-C_{\alpha}}$.
\item[\textbf{5.}] If $\beta > \alpha$, then $\delta \geq \alpha -1.$ Moreover, if $\beta > \alpha$ and $\delta =\alpha -1$, then $C_{\delta} \leq C_{\gamma}-C_{\alpha}$.
\end{itemize}
\label{theo:sharpdeux}
\end{theorem}
\noindent \textbf{Remark:} All our results until Section~\ref{sec:pta} are either valid under the general setting of Section~\ref{subsec:gensett}, or the computations can be adapted to the allometric setting of Section~\ref{subsec:allomsett} with $\alpha >1$. This is why in the following, we consider that the proofs for every point of Section~\ref{theo:sharp} give the similar points in Theorem~\ref{theo:sharpdeux}. However, we cannot use the conclusions of Section~\ref{sec:pta} (see the remark before Lemma~\ref{lemm:yt}). This explains the weaker conditions obtained in Theorem~\ref{theo:sharpdeux}, when we compare it to Theorem~\ref{theo:short}. Still, we think it is important to present our conclusions in that case. Indeed, DeLong and al. measured that metabolic rate scales with body mass superlinearly (\textit{i.e.} $\alpha >1$) for prokaryotes \cite{delong2010}. Malerba and Marshall also highlighted this superlinear scaling within a population of green algae \cite{malerba_2019}. We visualize the restrictions on allometric coefficients in the case $\alpha >1$ thanks to Figure~\ref{fig:alphastrictun}.

\section{Discussion}
\label{sec:discussion}
Before going into the mathematical details, we discuss the biological interpretation of our assumptions in Section~\ref{subsec:pasdenom}, and Section~\ref{discuss:final} gives our general conclusions and perspectives.

\subsection{Revisiting dimensional arguments of the Metabolic Theory of Ecology}
\label{subsec:pasdenom}

In Assumption~\ref{hyp:tempsinf}, we emphasize the ratio $(b_{x_{0}}+d)/g(.,R)$. The numerator accounts for the influence of birth and death rates, whereas the denominator represents the speed of acquisition or loss of energy by the individual in environment $R$. A classical idea is that $(b_{x_{0}}+d)/g(.,R)$ can be expressed as the inverse of an energy, because one can view $b_{x_{0}}$ and $d$ as dimensionless rates and $g$ as a gain of energy per unit of time \cite{savdee_08}, \cite{malerba_2019}. Usually, this is a justification for the choice $\beta=\delta=\alpha -1$. We underline here that this choice for the allometric coefficients is only a sufficient condition to obtain a biologically relevant individual process.
\\\\
Precisely, Assumption~\ref{hyp:probamort} and Assumption~\ref{hyp:tempsinf} indicate that the ratio $(b_{x_{0}}+d)/g(.,R)$ has divergent integrals, respectively near 0 and $+ \infty$. One possibility is to have the same allometric coefficient for $b_{x_{0}}$ and $d$, and in that case we necessarily have $\beta=\delta=\alpha-1$, so that the previous ratio is of order $1/x$. But there are plenty of other possibilities to obtain this divergent integral: for example, it happens when $\beta > \alpha -1$ and $\delta = \alpha -1$, and we see in Figure~\ref{fig:coeff} that the green area allows this kind of other reasonable choice. Hence, the choice $\beta=\delta=\alpha -1$ is sufficient but not necessary to verify Assumptions~\ref{hyp:probamort} and~\ref{hyp:tempsinf}. This is why we need an assumption about the population, which is Assumption~\ref{ass:supercritical}, to be more precise about the interspecific allometric constraints. This last assumption goes beyond the considerations about the ratio $(b_{x_{0}}+d)/g(.,R)$, so beyond the classical dimensional argument.

\subsection{General conclusion and perspectives}
\label{discuss:final}
To put it in a nutshell, here are the assumptions we can make for the conception of our model:
\begin{itemize}
\item The functions ruling energetic dynamics and the underlying birth and death process have allometric shapes.
\item An individual almost surely never reaches energy 0 before dying.
\item An individual almost surely never reaches energy $+ \infty$ before dying. Also, an individual has almost surely a finite lifetime.
\item Any population, characterized by the energy at birth $x_{0}$, survives with positive probability, at least for a sufficient amount of resources.
\end{itemize}
The way we put it here highlights that this is really the first allometric assumption that is restrictive and leads to the conclusions of Theorem~\ref{theo:short}. The other assumptions are strikingly basic, even if we underline that what we call a death event only accounts for intrinsic causes of death (\textit{i.e.} individuals only die because of the energy loss due to metabolism, and not competition factors). Notice that the previous assumptions also imply the following biologically relevant and basic features:
\begin{itemize}
\item An individual has almost surely a finite number of direct offspring.
\item An individual can always feed himself and grow if there are enough resources in the environment.
\end{itemize}
Recall here that a crucial aspect of our reasoning is that we want our model to respect the previous constraints for any kind of living species, tiny or gigantic. Thus, it is certainly adapted for the study of allometries involving several species, in the spirit of the guiding work of \cite{peters_1983} where the linear regressions are performed on a broad range of living species (nearly 20 orders of magnitude of body mass). The constraints we obtain are interspecific, and not intraspecific. One can refer to \cite{koz_1997} or \cite{malerba_2019} to understand this gap between interspecific allometries and within-species dynamics. An interesting opening would be to study the evolution of the allometric coefficient $\alpha$, adding a mutation mechanism in our previous setting. In other terms, we could adopt an evolutive point of view in line with the philosophy of adaptive dynamics \cite{metz1995adaptive, OdoDiekmann2003}.
\\\\
Our feeling is that the main phenomenon at stake here is the universality of $\alpha$ along $x_{0}$. It is eventually not a surprise that we get close to the results of the Metabolic Theory of Ecology. The simple case considered in this paper is meant to give an insight about what kind of mechanisms, other than competition for resources or trait evolution, are important to balance the energetic dynamics. To aim for the complex reality of an entire food web, one would add interaction between individuals and/or species and a dynamic evolution of the resource (see \cite{loeuille2005} or \cite{fritsch:hal-02129313} for a general approach). We leave this for future work.

\section{Construction of the process and proofs}
\label{proofs}
In Section~\ref{subsec:construc}, we construct the individual process highlighted in Section~\ref{subsec:gensett}. In Section~\ref{theo:sharp}, we give the different steps of the proof of Theorem~\ref{theo:short}. Even if we chose to highlight the allometric consequences of our reasoning in Section~\ref{sec:notations}, the results presented in the following are valid under the general setting of Section~\ref{subsec:gensett}, so their interpretation goes beyond the allometric setting (except for Sections~\ref{subsec:laprevue} and ~\ref{sec:pta}).

\subsection{Construction of the individual process $(\xi_{t})_{t \geq 0}$}
\label{subsec:construc}

Let us fix $x_{0}, R, \xi_{0}$. Recall that $A : (\xi_{0},t) \mapsto A_{\xi_{0}}(t)$ is the homogeneous-in-time flow associated to Equation~\eqref{eq:indivenergy} ($R$ is implicit in the notations). For $\xi_{0}>0$, $A_{\xi_{0}}(.)$ is defined on $[0,t_{\mathrm{max}}(\xi_{0})[$, where $t_{\mathrm{max}}(\xi_{0})$ is the deterministic time when $A_{\xi_{0}}(t_{\mathrm{max}}(\xi_{0})) = 0$ or $+ \infty$ ($t_{\mathrm{max}}(\xi_{0})$ is equal to $+ \infty$ if this never happens). One way to construct $(\xi_{t})_{t \geq 0}$ is to use two ingredients: 1) we first construct an auxiliary process $(\xi^{\mathrm{aux}}_{t})_{t \geq 0}$ starting from $\xi^{\mathrm{aux}}_{0} = \xi_{0}$ with birth events at rate $(b_{x_{0}}+d)\mathbb{1}_{\{\xi_{t}^{\mathrm{aux}} >x_{0}\}}$ and death event at rate $d\mathbb{1}_{\{\xi_{t}^{\mathrm{aux}} \leq x_{0} \}}$, and then 2) kill this process at a random time to obtain the appropriate birth and death rates $b_{x_{0}}$ and $d$. This construction of the process $\xi$ is made in the spirit of a Gillespie algorithm \cite{GILLESPIE1976403}. Although this construction may not seem the most natural, it will be convenient for couplings in the proofs in the following. We insist on the fact that for this construction, we only work under the general setting of Section~\ref{subsec:gensett}, and in particular without Assumptions~\ref{hyp:probamort} and~\ref{hyp:tempsinf}. This is essential if we want Theorem~\ref{theo:cns} to be meaningful.
\begin{itemize}
\item Let $(E_{i})_{i \geq 1}$ be i.i.d. random variables with exponential laws of parameter 1. First, if
$$\int_{0}^{t_{\mathrm{max}}(\xi_{0})}(b_{x_{0}}+d)(A_{\xi_{0}}(s)) \mathrm{d}s \leq E_{1}, $$
we set $J_{1}:=+ \infty$. In that case, we set $\xi^{\mathrm{aux}}_{t} = A_{\xi_{0}}(t)$ for $t \in [0,t_{\mathrm{max}}(\xi_{0})[$ and $\xi^{\mathrm{aux}}_{t} = \flat$ for $t \geq t_{\mathrm{max}}(\xi_{0})$. Recall that reaching $\flat$ means that the individual energy exploded or reached 0 before time $t$. Otherwise, as $(b_{x_{0}}+d)$ is a positive function, we can define the first time of jump as 
$$ J_{1} := \inf \left\{ t \in [0,t_{\max}(\xi_{0})[, \int_{0}^{t}(b_{x_{0}}+d)(A_{\xi_{0}}(s)) \mathrm{d}s = E_{1} \right\}. $$
In that case, $A_{\xi_{0}}(J_{1})$ is well-defined, and we set $\xi^{\mathrm{aux}}_{t} = A_{\xi_{0}}(t)$ for $t \in [0,J_{1}[$, and $$\xi^{\mathrm{aux}}_{J_{1}} = \left( A_{\xi_{0}}(J_{1}) - x_{0} \right) \mathbb{1}_{\{A_{\xi_{0}}(J_{1}) > x_{0}\}} + \partial \mathbb{1}_{\{A_{\xi_{0}}(J_{1}) \leq x_{0}\}}.$$
Remark that we need to check if $A_{\xi_{0}}(J_{1}) \leq x_{0}$. If this occurs, the jump is necessarily a death jump, as soon as $b_{x_{0}} \equiv 0$ on $]0,x_{0}]$, and then we set $\xi^{\mathrm{aux}}_{t} = \partial$ for $t \geq J_{1}$. We can derive the distribution function of $J_{1}$. For every $u \in ]0, t_{\max}(\xi_{0})[$:
\begin{align}
\mathbb{P}_{x_{0}, R,\xi_{0}}(J_{1} < u) = \displaystyle{\int_{0}^{u}(b_{x_{0}}+d)(A_{\xi_{0}}(t)) e^{- \displaystyle{\int_{0}^{t} (b_{x_{0}}+d)(A_{\xi_{0}}(s)) \mathrm{d}s }} \mathrm{d}t},
\label{eq:jun}
\end{align} 
and $\mathbb{P}_{x_{0}, R,\xi_{0}}(J_{1} = + \infty) = 1 - \mathbb{P}_{x_{0}, R,\xi_{0}}(J_{1} < t_{\max}(\xi_{0}))$.
\\\\
Now, let us suppose that we defined $J_{n}$ and $\xi_{J_{n}}^{\mathrm{aux}}$ for some $n \geq 1$. If $J_{n}= + \infty$, or $J_{n}< + \infty$ and $\xi^{\mathrm{aux}}_{J_{n}} = \partial$, we simply set $J_{n+1}:= + \infty$. In that case, $\xi^{\mathrm{aux}}_{t}$ is already defined for all $t \geq 0$, because it already reached $\flat$ or $\partial$. Now, we suppose that $J_{n} < + \infty$ and $\xi^{\mathrm{aux}}_{J_{n}} \neq \partial$, and we distinguish between two remaining cases. First, if
$$\int_{J_{n}}^{J_{n} + t_{\mathrm{max}}(\xi^{\mathrm{aux}}_{J_{n}})}(b_{x_{0}}+d)(A_{\xi^{\mathrm{aux}}_{J_{n}}}(s-J_{n}))\mathrm{d}s \leq E_{n+1}, $$
we again set $J_{n+1}:=+ \infty$, $\xi^{\mathrm{aux}}_{t} = A_{\xi^{\mathrm{aux}}_{J_{n}}}(t-J_{n})$ for $t \in [J_{n},J_{n} +t_{\mathrm{max}}(\xi^{\mathrm{aux}}_{J_{n}})[$ and $\xi^{\mathrm{aux}}_{t} = \flat$ for $t \geq J_{n} +t_{\mathrm{max}}(\xi^{\mathrm{aux}}_{J_{n}})$. Otherwise, we can define the $(n+1)$-th time of jump as 
$$ J_{n+1} := \inf \left\{ t \in [J_{n},J_{n}+t_{\max}(\xi^{\mathrm{aux}}_{J_{n}})[, \int_{J_{n}}^{t}(b_{x_{0}}+d)(A_{\xi^{\mathrm{aux}}_{J_{n}}}(s-J_{n}))\mathrm{d}s = E_{n+1} \right\}. $$
In that last case, we set $\xi^{\mathrm{aux}}_{t} = A_{\xi^{\mathrm{aux}}_{J_{n}}}(t-J_{n})$ for $t \in [J_{n},J_{n+1}[$, and $$\xi^{\mathrm{aux}}_{J_{n+1}} = \left( A_{\xi^{\mathrm{aux}}_{J_{n}}}(J_{n+1}-J_{n}) - x_{0} \right) \mathbb{1}_{\{A_{\xi^{\mathrm{aux}}_{J_{n}}}(J_{n+1}-J_{n}) > x_{0}\}} + \partial \mathbb{1}_{\{A_{\xi^{\mathrm{aux}}_{J_{n}}}(J_{n+1}-J_{n}) \leq x_{0}\}}.$$ 
If $\xi^{\mathrm{aux}}_{J_{n+1}} = \partial$, we then set $\xi^{\mathrm{aux}}_{t} = \partial$ for $t \geq J_{n+1}$.
Eventually, it is possible, in the general setting of Section~\ref{subsec:gensett} and without further assumptions, that $\mathfrak{J} := \sup_{n \in \mathbb{N}^{*}} J_{n} <+ \infty$ (\textit{i.e.} we have an accumulation of jump times). In that final case, we set $\xi^{\mathrm{aux}}_{t} = \flat$ for $t \geq \mathfrak{J}$.
\item Then, we define a sequence $(U_{i})_{i \geq 1}$ of i.i.d. random variables with uniform laws on $[0,1]$, independent of $(E_{i})_{i \geq 1}$. We define the time of death $T_{d}$ as
\[ T_{d} := \inf\left\{J_{i}, \hspace{0.1 cm} \xi^{\mathrm{aux}}_{J_{i}} \neq \flat \hspace{0.1 cm} \mathrm{and} \hspace{0.1 cm}  U_{i} \leq \dfrac{d}{b_{x_{0}}+d}\left(\xi^{\mathrm{aux}}_{J_{i}-}\right) \right\},\]
with the convention $\inf(\varnothing) = + \infty$. We also define $T_{0} := \inf \{ t \geq 0, \xi^{\mathrm{aux}}_{t} = \flat, \xi^{\mathrm{aux}}_{s} \xrightarrow[s \rightarrow t-]{} 0 \}$ and $T_{\infty} := \inf \{ t \geq 0, \xi^{\mathrm{aux}}_{t} = \flat, \xi^{\mathrm{aux}}_{s} \xrightarrow[s \rightarrow t-]{} + \infty \}$. Remark that $\{ T_{d} <  T_{0} \wedge T_{\infty} \} \subseteq \{ \mathfrak{J} = + \infty \}$, so if our process is biologically relevant as defined in Section \ref{subsec:cns}, there is no accumulation of jump times. We then define our final process $\xi$ for $t \in \mathbb{R}^{+}$ as
$$ \xi_{t} := \xi^{\mathrm{aux}}_{t}\mathbb{1}_{\{ t < T_{d} \}} + \partial \mathbb{1}_{\{ t \geq T_{d} \}}. $$
\end{itemize}
The reader can check that if for some $n \in \mathbb{N}$, the condition $\{A_{\xi^{\mathrm{aux}}_{J_{n}}}(J_{n+1}-J_{n}) \leq x_{0}\}$ is verified, then we have $ \frac{d}{b_{x_{0}}+d}\left(\xi^{\mathrm{aux}}_{J_{n}-}\right) =1 \geq U_{n}$, so that our jumps to $\partial$ for $\xi^{\mathrm{aux}}$ are consistent with the definition of $T_{d}$. These jumps for $\xi^{\mathrm{aux}}$ to $\partial$ are necessary to keep a well-defined positive (when it is not $\partial$ or $\flat$) process $\xi^{\mathrm{aux}}$. They account for deaths of individuals with energy smaller than $x_{0}$, whereas the definition of $T_{d}$ also includes deaths with an energy higher than $x_{0}$. Finally, birth times for the process $\xi$ are exactly times $J_{i}$ such that $J_{i} < + \infty$ and $J_{i} \neq T_{d}$, so the number of direct offspring of an individual during its life, is defined as 
\begin{align}
N_{x_{0},R,\xi_{0}} := \sup \{ i \geq 0, J_{i} < T_{d} \}.
\label{eq:nombre}
\end{align}

\subsection{Sketch of the proof of Theorem~\ref{theo:cns} and Theorem~\ref{theo:short}}
\label{theo:sharp}

We prove Theorem~\ref{theo:cns} by splitting the different possible situations depending on the resource $R$. We consider the case $R \in \mathfrak{R}_{0} \setminus \mathfrak{R}_{\infty}$ in Proposition~\ref{eq:delta} of Section~\ref{subsec:delta}, the case $R \notin \mathfrak{R}_{0} \cup \mathfrak{R}_{\infty}$ in Lemma~\ref{lemme:premierr} of Section~\ref{subsec:teun}, and conclude for $R \in \mathfrak{R}_{\infty}$ in Proposition~\ref{prop:tinfyinfini} of Section~\ref{subsec:teun}. Also, we divide the proof of Theorem~\ref{theo:short} in eight distinct points. Let $0 < \alpha \leq 1$. Under the allometric setting of Section~\ref{subsec:allomsett}, we have:
\begin{itemize}
\item[\textbf{1.}]Assumption~\ref{ass:gainenergy} $ \Leftrightarrow (\gamma = \alpha \hspace{0.1 cm} \mathrm{and} \hspace{0.1 cm} C_{\gamma} > C_{\alpha})$.
\item[\textbf{2.}] Assumption~\ref{hyp:probamort} $\Leftrightarrow  \delta \leq \alpha-1$.
\item[\textbf{3.}] Assumptions~\ref{hyp:tempsinf} and~\ref{ass:gainenergy}  $\Rightarrow ($Assumption~\ref{ass:pasinfini} $\mathrm{and} \hspace{0.1 cm} \mathrm{max}(\beta, \delta) \geq \alpha-1)$.
\item[\textbf{4.}] $($Assumptions~\ref{hyp:probamort}, \ref{hyp:tempsinf}, \ref{ass:supercritical} $\mathrm{and} \hspace{0.1 cm} \mathrm{max}(\beta, \delta) \geq \alpha-1) \Rightarrow (\beta \geq \alpha-1$, and if $\beta=\delta=\alpha-1$, then $C_{\beta} > C_{\delta})$.
\item[\textbf{5.}] (Assumptions~\ref{hyp:tempsinf}, \ref{ass:supercritical} and $\delta=\alpha-1$) $\Rightarrow$ $\bigg($if $\beta > \alpha-1$, then $\beta \geq \alpha -1 + \dfrac{C_{\delta}}{C_{\gamma}-C_{\alpha}}\bigg)$.
\item[\textbf{6.}] Assumptions~\ref{hyp:probamort}, \ref{hyp:tempsinf}, \ref{ass:supercritical} $\Rightarrow$ $($if $\beta > \alpha$, then $\delta \geq \alpha -1)$.
\item[\textbf{7.}] (Assumptions~\ref{hyp:tempsinf}, \ref{ass:supercritical} and $\delta=\alpha-1$) $\Rightarrow$ $($if $ \beta > \alpha-1$, then $C_{\delta} \leq C_{\gamma}-C_{\alpha})$.
\item[\textbf{8.}] Assumptions~\ref{hyp:tempsinf} and~\ref{ass:gainenergy} $\Rightarrow$ $($if $ \beta \leq\alpha$, then $\delta \geq \alpha -1)$.
\end{itemize}
\noindent \textbf{Remark:} Theorem~\ref{theo:short} is a consequence of points \textbf{1} to \textbf{8}. Thanks to Lemma~\ref{lemm:ass} and point \textbf{3}, we know that in Theorem~\ref{theo:short}, we work in fact under Assumptions~\ref{hyp:probamort} to~\ref{ass:gainenergy}, and $\mathrm{max}(\beta, \delta) \geq \alpha-1$. Then, $\gamma = \alpha$ and $C_{\gamma} > C_{\alpha}$ because of \textbf{1}. Also, $\delta = \alpha-1$ because of \textbf{2, 6, 8}. Finally, in Theorem~\ref{theo:short}, the remaining conditions on $\beta, C_{\beta}, C_{\delta}, C_{\gamma}, C_{\alpha}$ are the combination of \textbf{4, 5, 7}. 
\\\\
For the previously enumerated results, we precise where the reader can find the associated proofs, with some insight on the mathematical tools involved. Except for points \textbf{6, 7, 8}, the following results are stated and valid under the general setting of Section~\ref{subsec:gensett}.
\begin{itemize}
\item[\textbf{1.}] See Proposition~\ref{eq:gammalpha} in Section~\ref{subsec:energain}, which is purely deterministic.
\item[\textbf{2.}] Under the allometric setting of Section~\ref{subsec:allomsett}, Equation~\eqref{eq:deltacahnge} is:
$\forall x_{0}>0, \displaystyle{\int_{0}^{x_{0}} \dfrac{C_{\delta}}{C_{\alpha}}x^{\delta-\alpha} \mathrm{d}x}= + \infty,$ which is equivalent to $\delta \leq \alpha -1$.
\item[\textbf{3.}] See Corollary~\ref{corr:swagosss} in Section~\ref{subsec:teun}. The proof is based on an operator point of view on our probabilistic questioning.
\item[\textbf{4.}] See Proposition~\ref{lemm:lemmaone} in Section~\ref{proof:super}, which uses basic probability calculus.
\item[\textbf{5.}] See Proposition~\ref{corr:pointquatre} in Section~\ref{subsec:deltaalphamoinsun}. The proof is based on Proposition~\ref{prop:beplusde} in Section~\ref{proof:beplusde}, which uses a useful coupling to compare different individual trajectories.
\item[\textbf{6.}] This will be proven in Sections~\ref{subsec:lowenergy} and \ref{subsec:highenergy}. We define several couplings of our individual process with simpler processes, and also track the maximal energy reached by an individual.
\item[\textbf{7.}] We apply the reasoning of \textbf{6.} in the case $\delta=\alpha-1$. This constitutes Corollary~\ref{corr:thirdone} in Section~\ref{sec:coeff}.
\item[\textbf{8.}] We use a representation of the individual process $(\xi_{t})_{t \geq 0}$ with Poisson point measures, and martingale techniques. We were inspired by the concept of asymptotic pseudotrajectory developed by Benaïm and Hirsch \cite{benaim_06}. The complete proof of this last point is developed in Section~\ref{sec:pta}.

\end{itemize}

\subsection{Proof of 1. in Section~\ref{theo:sharp}}
\label{subsec:energain}
We begin with a simple result that illustrates perfectly our way of reasoning. Also, this section is particular in the sense that our considerations here are purely deterministic. Assumption~\ref{ass:gainenergy} translates into a condition on functions $\psi$ and $\ell$, ruling the energy dynamics in Equation~\eqref{eq:indivenergy}.
\begin{prop}
\label{eq:gammalpha}
Under the general setting of Section~\ref{subsec:gensett}, Assumption~\ref{ass:gainenergy} is equivalent to
\begin{align}
\forall x>0, \quad \psi(x) > \ell(x).
\label{prop:gammalpha}
\end{align}

Under the allometric setting of Section~\ref{subsec:allomsett},~\eqref{prop:gammalpha} is equivalent to
\begin{align*}
\gamma= \alpha \hspace{0.1 cm}\mathrm{and} \hspace{0.1 cm} C_{\gamma} > C_{\alpha}.
\end{align*}
\end{prop}

\begin{proof}
Let $x>0$, we have $g(x,R)= \phi(R)\psi(x) - \ell(x)$. It is straightforward to verify that Assumption~\ref{ass:gainenergy} is equivalent to~\eqref{prop:gammalpha}, because $\phi$ is non-negative, non-decreasing and $\phi(R) \rightarrow 1$ when $R \rightarrow +\infty$.
\\
Under the allometric setting of Section~\ref{subsec:allomsett},  by letting $x \rightarrow + \infty$, we see that~\eqref{prop:gammalpha} cannot be verified unless $\gamma > \alpha \quad \mathrm{or} \quad (\gamma= \alpha \hspace{0.1 cm}\mathrm{and} \hspace{0.1 cm} C_{\gamma} > C_{\alpha}).$
\\
In the same vein, if we consider now $x \rightarrow 0$, we also obtain regarding~\eqref{prop:gammalpha} that $\gamma < \alpha \quad \mathrm{or} \quad (\gamma= \alpha \hspace{0.1 cm}\mathrm{and} \hspace{0.1 cm} C_{\gamma} > C_{\alpha})$, which ends the direct implication. It is easy to see that the converse also holds true.
\end{proof}
\noindent \textbf{Remark:} This result gives point \textbf{1.} of Section~\ref{theo:sharp}. We immediately observe one of the major feature of the allometric setting: power functions of the form $x \mapsto x^{\kappa}$ have an antagonistic behavior near 0 and near $+ \infty$, depending on the sign of $\kappa$. In order to enforce general principles valid for any $x_{0}$, one has to choose precise allometric coefficients.
\\\\
Under the allometric setting of Section~\ref{subsec:allomsett}, an immediate consequence of Proposition~\ref{eq:gammalpha} is that, under Assumption~\ref{ass:gainenergy}, Equation~\eqref{eq:indivenergy} becomes
\begin{align}
    \dfrac{\mathrm{d}\xi_{t}}{\mathrm{d}t} =  C_{R}\xi_{t}^{\alpha},
    \label{eq:indivenergymod}
\end{align}
with $C_{R} := \phi(R)C_{\gamma} - C_{\alpha} \leq C_{\gamma}-C_{\alpha}$. Remark that if $0 < \alpha \leq 1$, for every $(R,\xi_{0})$, the solution of Equation~\eqref{eq:indivenergy} starting from $\xi_{0}$ with avalaible resources $R$ cannot explode in finite time, but can possibly reach 0 depending on the value of $R$. In the case $\alpha >1$, a solution of~\eqref{eq:indivenergymod} with $C_{R}>0$ can reach $+ \infty$ in finite time, but not 0. These situations have been taken into account in the construction of the individual process $(\xi_{t})_{t \geq 0}$ in Section~\ref{subsec:construc}.

\subsection{Equivalent formulations of Assumption~\ref{hyp:probamort}}
\label{subsec:delta}
In the following, we give probabilistic interpretations of Assumption~\ref{hyp:probamort}.
\begin{lemme}
Under the general setting of Section~\ref{subsec:gensett}, Assumption~\ref{hyp:probamort} is equivalent to
\begin{align}
\forall x_{0} >0, \forall \xi_{0} >0, \quad \mathbb{P}_{x_{0}, R=0, \xi_{0}}(T_{d} < + \infty) = 1.
\label{eq:equivalencetre}
\end{align}
\label{lemme:latruite}
\end{lemme}
\begin{proof}
First, let us show that Assumption~\ref{hyp:probamort} is equivalent to the weaker condition
\begin{align}
\mathbb{P}_{x_{0}, R=0, \xi_{0}=x_{0}}(T_{d} < + \infty) = 1
\label{laproba}
\end{align}
for every $x_{0}>0$. We fix $x_{0}$, $R=0$ and $\xi_{0}=x_{0}$ and investigate under which condition~\eqref{laproba} holds true. In this situation, the energy only decreases over time, and the only random event possibly occurring is a death, so from the construction of Section~\ref{subsec:construc}, $\{ T_{d} < + \infty \} = \{ J_{1} < t_{\mathrm{max}}(x_{0},0) \}$ (with $t_{\mathrm{max}}(x_{0},0)$ possibly equal to $+ \infty$). We use~\eqref{eq:jun} with $(b_{x_{0}}+d) \equiv d$ on $]0,x_{0}]$ to obtain  
\begin{align*}
\mathbb{P}_{x_{0},0,x_{0}}(T_{d} <+ \infty)
& = \mathbb{P}_{x_{0},0,x_{0}}(J_{1} < t_{\mathrm{max}}(x_{0},0)) \\ & = \displaystyle{\int_{0}^{t_{\mathrm{max}}(x_{0},0)} d(A_{x_{0},0}(s)) e^{-\displaystyle{\int_{0}^{s} d(A_{x_{0},0}(\tau)) \mathrm{d}\tau}} \mathrm{d}s}\\ & = 1-e^{-\displaystyle{\int_{0}^{t_{\mathrm{max}}(x_{0},0)} d(A_{x_{0},0}(\tau)) \mathrm{d}\tau}}.
\end{align*}
Thus,~\eqref{laproba} is equivalent to ($e^{-\int_{0}^{t_{\mathrm{max}}(x_{0},0)} d(A_{x_{0},0}(\tau)) \mathrm{d}\tau} = 0$ for all $x_{0}>0$). The function $\ell$ is positive increasing and $A_{x_{0},0}(.)$ is the solution of~\eqref{eq:indivenergy} with $A_{x_{0},0}(0) = x_{0}$ and $R=0$ (so $g(x,R) = -\ell(x)$). Hence, $A_{x_{0},0}(t)$ goes to 0 when $t$ goes to $t_{\max}(x_{0},0)$, and we can perform the change of variables $u=A_{x_{0},0}(\tau)$, which gives exactly Assumption~\ref{hyp:probamort}. Now, we suppose that~\eqref{laproba} holds true and we show~\eqref{eq:equivalencetre}. The same reasoning as before still applies when $\xi_{0} \leq x_{0}$, because there are no birth events (we just replace $t_{\mathrm{max}}(x_{0},0)$ by $t_{\mathrm{max}}(\xi_{0},0)$). If $\xi_{0}>x_{0}$, we define the random time 
$$\tau_{x_{0},\xi_{0}} := \inf \{ t \in [0,T_{d}[, \hspace{0.1cm} \xi_{t,x_{0},0,\xi_{0}} \leq x_{0} \},$$ 
with the convention $\inf(\varnothing) = + \infty$. First if $\tau_{x_{0},\xi_{0}} = + \infty$, then $T_{d} < T_{0} = +\infty$ because before $T_{d}$, the energy is in $[x_{0},\xi_{0}]$ and the death rate is positive continuous, so lower bounded by a positive constant, on this segment. Else if $\tau_{x_{0},\xi_{0}} < + \infty$, by Markov property, we start at time $\tau_{x_{0},\xi_{0}}$ from a new initial condition $\xi_{1} \leq x_{0}$, for which the result is already proven. 
\end{proof}

\noindent In fact, Assumption~\ref{hyp:probamort} also ensures that individual trajectories almost surely avoid 0 for every $x_{0}>0$, $\xi_{0}>0$, and any amount of resources $R \geq 0$.

\begin{prop}
Under the general setting of Section~\ref{subsec:gensett} and under Assumption~\ref{hyp:probamort}, we have
\begin{align*}
\forall x_{0} >0, \forall \xi_{0} >0, \forall R \geq 0, \quad \mathbb{P}_{x_{0}, R, \xi_{0}}(T_{0} = + \infty) = 1.
\end{align*}
\label{prop:toinfini}
\end{prop}

\begin{proof}
Let $x_{0}, \xi_{0}>0$ and $R \geq 0$ and assume by contradiction that $\mathbb{P}_{x_{0},R,\xi_{0}}(T_{0}< + \infty) >0$. From the construction of Section \ref{subsec:construc}, on the event $\{ T_{0}< + \infty \}$, there exists a family $(E_{i})_{i \geq 1}$ of i.i.d. exponential random variables with parameter 1, a family $(U_{i})_{i \geq 1}$ of i.i.d. uniform random variables on $[0,1]$ and $n \in \mathbb{N}^{*}$ such that $J_{n-1}< + \infty$ and 
$$\int_{0}^{t_{\mathrm{max}}(\xi^{\mathrm{aux}}_{J_{n-1}},R)}(b_{x_{0}}+d)(A_{\xi^{\mathrm{aux}}_{J_{n-1}},R}(s))\mathrm{d}s \leq E_{n}, $$
with the convention $J_{0}=0$. Moreover, as $T_{0}< + \infty$, we necessarily have
$$ A_{\xi^{\mathrm{aux}}_{J_{n-1}},R}(s) \xrightarrow[s \rightarrow t_{\mathrm{max}}(\xi^{\mathrm{aux}}_{J_{n-1}},R)]{} 0,$$
meaning that the flow $s \mapsto A_{\xi^{\mathrm{aux}}_{J_{n-1}},R}(s) $ is decreasing and $g(x,R) <0$ for every $x \in ]0,\xi^{\mathrm{aux}}_{J_{n-1},R}]$, according to Equation \eqref{eq:indivenergy}. The change of variables $x=A_{\xi^{\mathrm{aux}}_{J_{n-1}},R}(s)$ is then licit and we have
$$\int_{0}^{\xi^{\mathrm{aux}}_{J_{n-1}}}\dfrac{(b_{x_{0}}+d)(x)}{-g(x,R)}\mathrm{d}x \leq E_{n} < + \infty. $$
As $g$ is increasing in $R$ and $b_{x_{0}} \equiv 0 $ on $]0,x_{0}]$, this leads to
$$\int_{0}^{\xi^{\mathrm{aux}}_{J_{n-1}}\wedge x_{0}}\dfrac{d(x)}{-g(x,0)}\mathrm{d}x < + \infty,$$
which contradicts Assumption \ref{hyp:probamort} (because $g(x,0) \equiv -\ell$) and concludes the proof.
\end{proof}

Finally, we can give a stronger equivalent formulation of Assumption \ref{hyp:probamort}.

\begin{prop}
Under the general setting of Section~\ref{subsec:gensett}, Assumption~\ref{hyp:probamort} is equivalent to
\begin{align*}
\forall x_{0} >0, \forall \xi_{0} >0, \forall R \in \mathfrak{R}_{0} \setminus \mathfrak{R}_{\infty}, \quad \mathbb{P}_{x_{0}, R, \xi_{0}}(T_{d} < + \infty) = 1.
\end{align*}
\label{eq:delta}
\end{prop}

\begin{proof}
First, remark that $0 \in \mathfrak{R}_{0} \setminus \mathfrak{R}_{\infty}$, so by Lemma~\ref{lemme:latruite}, we only have to show that Assumption~\ref{hyp:probamort} implies $\mathbb{P}_{x_{0}, R, \xi_{0}}(T_{d} < + \infty) = 1$ for $R >0$ in $\mathfrak{R}_{0} \setminus \mathfrak{R}_{\infty}$ and any $x_{0}>0$, $\xi_{0}>0$ that we fix for now. Because $R \notin \mathfrak{R}_{\infty}$, there exists $x_{2} >\xi_{0}$ such that $g(x_{2},R) \leq 0$. Then, $T_{\infty} = + \infty$, and for $t < T_{d} \wedge T_{0}$, $\xi_{t,x_{0},R,\xi_{0}} \leq x_{2}$. Also, by Proposition~\ref{prop:toinfini}, we also have $T_{0}=+\infty$ almost surely. However, $R \in \mathfrak{R}_{0}$, so there exists $x_{1} > 0$, such that $g(y,R)<0$ for $y \leq x_{1}$. We define the random time 
$$\tau := \inf \{ t \in [0,T_{d}[, \hspace{0.1cm} \xi_{t,x_{0},R,\xi_{0}} \leq x_{1} \wedge x_{0} \},$$ 
with the convention $\inf(\varnothing) = + \infty$. First, if $\tau = + \infty$, we have $\xi_{t,x_{0},R,\xi_{0}} \in [x_{1} \wedge x_{0},x_{2}]$ before $T_{d}$ and the death rate is positive continuous on this segment so $T_{d}< + \infty$. Else if $\tau < + \infty$, by Markov property, we start at time $\tau$ from a new initial condition $\xi_{1} \leq x_{1} \wedge x_{0}$. As $\frac{d(x)}{-g(x,R)} > \frac{d(x)}{\ell(x)} $ for $x \leq \xi_{1}$, the same reasoning as Lemma~\ref{lemme:latruite} applies, from the initial condition $\xi_{1}$ with resource $R$, so $T_{d}< + \infty$. 
\end{proof}

\subsection{Proof of 3. in Section~\ref{theo:sharp} and consequences of Assumption \ref{hyp:tempsinf}}
\label{subsec:teun}
First, we prove intermediate lemmas that lead to Corollary~\ref{corr:swagosss}, which entails point \textbf{3.} in Section~\ref{theo:sharp}. Then, we will give in Proposition~\ref{prop:tinfyinfini} an equivalent formulation of Assumption~\ref{hyp:tempsinf} under Assumption~\ref{hyp:probamort}, which completes the proof of Theorem~\ref{theo:cns} with Proposition~\ref{eq:delta} and Lemma~\ref{lemme:premierr}. We begin with a technical lemma on $\mathfrak{J}$, the supremum of jump times of $\xi^{\mathrm{aux}}$ defined in Section~\ref{subsec:construc}.
\begin{lemme}
Under the general setting of Section~\ref{subsec:gensett}, for every $x_{0},R,\xi_{0}$, we have almost surely
$$ \{ \mathfrak{J}<+\infty \} \subseteq \{\mathfrak{J}=T_{\infty} \}.$$
\label{lemme:junfini}
\end{lemme}
\begin{proof}
We fix $x_{0},R,\xi_{0}$ in the following. Remark that $ \{ \mathfrak{J}< + \infty \} \subseteq \left\{ \forall n \in \mathbb{N}^{*}, \xi^{\mathrm{aux}}_{J_{n}} \notin \{\partial, \flat\} \right\} \subseteq \{ \forall n \in \mathbb{N}^{*}, \; J_{n} < T_{0} \wedge T_{\infty} \wedge T_{d} \}$. Thus, on the event $\{ \mathfrak{J}< + \infty \}$, the only possible random jumps on the finite time interval $[0,\mathfrak{J}[$ are birth events, and there is an infinite amount of such birth events. Remark that the event $\{\mathfrak{J} < T_{\infty} \}$ is equal to
\begin{align}
\left\{ \exists M>0, \forall t \in [0, \mathfrak{J}[, \exists s \in ]t, \mathfrak{J}[, \xi^{\mathrm{aux}}_{s} \leq M \right\},
\label{eq:infinijn}
\end{align}
so we suppose by contradiction that $\{ \mathfrak{J} < + \infty\}$ and $\eqref{eq:infinijn}$ occur with positive probability and we work on these events in the following. We define $\tau := \sup \{t \in [0, \mathfrak{J}[, \xi^{\mathrm{aux}}_{t} > M+1\}$, and distinguish between two cases. First, if $\tau < \mathfrak{J}$, it means that there is an infinite number of birth events  and $\xi^{\mathrm{aux}}_{t} \leq M+1$ on the finite time interval $]\tau, \mathfrak{J}[$. This happens with probability 0, because the birth rate $b_{x_{0}}$ is bounded on $[0,M+1]$. Else if $\tau = \mathfrak{J}$, by \eqref{eq:infinijn}, we can define two random sequences $(s_{k})_{k \in \mathbb{N}}$ and $(t_{k})_{k \in \mathbb{N}}$ such that for all $k \in \mathbb{N}$, $0 \leq s_{k} < t_{k} < s_{k+1} < \mathfrak{J}$, with $\xi^{\mathrm{aux}}_{s_{k}} \leq M$ and $\xi^{\mathrm{aux}}_{t_{k}} > M+1$. By definition of the process, the only way $\xi^{\mathrm{aux}}$ can increase is deterministic, between random jumps. Hence, there exists a deterministic time $\sigma>0$ such that for all $k \geq \mathbb{N}$, $t_{k}-s_{k} \geq \sigma$. This is a contradiction and concludes the proof, because we would then have 
$$ + \infty = \sum\limits_{k \geq 0} \sigma \leq \sum\limits_{k \geq 0} (t_{k}-s_{k}) \leq   \sum\limits_{k \geq 0} (s_{k+1}-s_{k}) \leq \mathfrak{J} < + \infty.$$
\end{proof}
For $R\geq 0$, we define $\Omega_{R} := \{ \xi_{0}>0, \forall x \geq \xi_{0}, g(x,R)>0\}$, and notice that $R \in \mathfrak{R}_{\infty}$, if and only if $\Omega_{R} \neq \varnothing$.

\begin{lemme}
Under the general setting of Section~\ref{subsec:gensett}, we have
$$ \forall x_{0}>0, \forall R \notin \mathfrak{R}_{0} \cup \mathfrak{R}_{\infty}, \forall \xi_{0}>0, \quad \mathbb{P}_{x_{0},R,\xi_{0}}(T_{d} < + \infty) =1,$$
and 
$$ \forall x_{0}>0, \forall R \in \mathfrak{R}_{\infty} \setminus \mathfrak{R}_{0}, \forall \xi_{0} \notin \Omega_{R}, \quad \mathbb{P}_{x_{0},R,\xi_{0}}(T_{d} < + \infty) =1.$$
\label{lemme:premierr}
\end{lemme}
\begin{proof}
Let $\xi_{0}>0$, $x_{0}>0$ and $R \notin \mathfrak{R}_{0} \cup \mathfrak{R}_{\infty}$, then 
$$ \forall x >0, \exists y \in ]0,x], \quad g(y,R) \geq 0, $$
and
$$ \forall x>0, \exists y \geq x, \quad g(y,R) \leq 0. $$
In that case, according to~\eqref{eq:indivenergy}, we have $T_{0}=T_{\infty}= + \infty$. In particular, $\mathfrak{J} = + \infty$ almost surely by Lemma~\ref{lemme:junfini}, because if $\mathfrak{J}<+ \infty$, we would have $\mathfrak{J}=T_{\infty}<+ \infty$. We can find $\overline{y} \geq x_{0} \vee \xi_{0}$ such that $g(\overline{y},R) \leq 0$, so before $T_{d}$, we have $\xi_{t} \leq \overline{y}$. If $R \in \mathfrak{R}_{\infty} \setminus \mathfrak{R}_{0}$ and $\xi_{0} \notin \Omega_{R}$, we verify that we can work in the exact same setting. Now, suppose by contradiction that $\{T_{d}=+\infty \}$ occurs with positive probability, and work on this event. We define
$$ \underline{\xi} := \inf_{t \geq 0} \xi_{t}. $$
First, if $\underline{\xi}>0$, the process stays in $[\underline{\xi},\overline{y}]$ over time. Otherwise, $\underline{\xi}=0$, and the only way this could happen with $R \notin \mathfrak{R}_{0}$ is that there is an infinite number of birth events. Hence, we can construct two random sequences $(s_{k})_{k \in \mathbb{N}}$ and $(t_{k})_{k \in \mathbb{N}}$ such that for all $k \in \mathbb{N}$, $0 \leq s_{k} < t_{k} < s_{k+1}$, with $\xi^{\mathrm{aux}}_{s_{k}} \leq x_{0}/2$ and $\xi^{\mathrm{aux}}_{t_{k}} > x_{0}$. As in the proof of Lemma~\ref{lemme:junfini}, this implies that the process stays almost surely an infinite amount of time in the segment $[x_{0}/2,x_{0}]$. In both cases, the process stays an infinite amount of time in a segment, where the death rate $d$ is lower-bounded by a positive constant, so the time of death $T_{d}$ would be almost surely finite, which is a contradiction. 
\end{proof}
Recall that $K_{x_{0},R}$ is the operator defined by~\eqref{eq:kixzerooo}. In the two following lemmas, we give a probabilistic interpretation of $K^{k}_{x_{0},R}\mathbf{1}$ for $k \geq 1$, as the probability that the $k$ first jumps of the process are birth events. For $k \geq 1$, we write $M^{k}_{x_{0},R,\xi_{0}}$, or simply $M^{k}$, for the event $\{ N_{x_{0},R,\xi_{0}} \geq k \}$ (\textit{i.e.} an individual starting from $\xi_{0}$ with characteristic energy $x_{0}$ and resources $R$ has at least $k$ direct offspring during its life), with $N_{x_{0},R,\xi_{0}}$ defined by \eqref{eq:nombre}. We write $B_{x_{0},R,\xi_{0}}$, or simply $B$, for the event $\{ J_{1} < T_{d} \}$, (\textit{i.e.} the first jump event is a birth). Also, we write $D_{x_{0},R,\xi_{0}}$, or simply $D$ for the event $\{ J_{1} = T_{d} \} \cap \{ J_{1} < t_{\mathrm{max}}(\xi_{0},R)\}$, (\textit{i.e.} the first jump event is a death). 
\begin{lemme}
Under the general setting of Section~\ref{subsec:gensett}, for every $x_{0},\xi_{0} >0$, $R \geq 0$ and $0<u < t_{\mathrm{max}}(\xi_{0})$, we have
$$\mathbb{P}_{x_{0},R,\xi_{0}}(\{J_{1} < u\} \cap B) = \displaystyle{\int_{0}^{u}b_{x_{0}}(A_{\xi_{0},R}(t))e^{- \displaystyle{\int_{0}^{t} (b_{x_{0}}+d)(A_{\xi_{0},R}(s)) \mathrm{d}s }}\mathrm{d}t },$$
In particular,
$$ \mathbb{P}_{x_{0},R,\xi_{0}}(B) = K_{x_{0},R}\mathbf{1}(\xi_{0}).$$
\label{lemme:lemmatwo}
\end{lemme}
\begin{proof}
From the construction of Section~\ref{subsec:construc}, on the event $\left\{J_{1} < u < t_{\mathrm{max}}(\xi_{0})\right\}$, the event $B$ is exactly $$\left\{U_{1} > \dfrac{d}{b_{x_{0}}+d}\left(A_{\xi_{0},R}(J_{1}) \right) \right\},$$
with $U_{1}$ being a uniform random variable on $[0,1]$ independent of $J_{1}$. Hence, conditioning on the value of $J_{1}$ and using~\eqref{eq:jun}, we obtain
\begin{align*}
\mathbb{P}_{x_{0},R,\xi_{0}}(\{J_{1} < u\} \cap B) & = \displaystyle{\int_{0}^{u}\left(1-\dfrac{d}{b_{x_{0}}+d}(A_{\xi_{0},R}(t))\right)(b_{x_{0}}+d)(A_{\xi_{0},R}(t)) e^{- \displaystyle{\int_{0}^{t} (b_{x_{0}}+d)(A_{\xi_{0},R}(s)) \mathrm{d}s }} \mathrm{d}t},
\end{align*}
which leads to the result. The second equality of the lemma comes from the limit $u \rightarrow t_{\mathrm{max}}(\xi_{0},R)$ (increasing events on the left-hand side and increasing integral when $u$ increases on the right-hand side).
\end{proof}
\begin{lemme}
Under the general setting of Section~\ref{subsec:gensett}, we have 
\begin{align}
\forall k \geq 1, \quad \forall x_{0}>0, \forall \xi_{0} >0, \forall R \geq 0, \hspace{0.1 cm} \mathbb{P}(M^{k}_{x_{0},R,\xi_{0}}) = K_{x_{0},R}^{k}\mathbf{1}(\xi_{0}),
\label{eq:probakxo}
\end{align}
so that 
\begin{center}
Assumption~\ref{ass:pasinfini} $\Leftrightarrow \bigg(\forall x_{0}>0, \forall R \geq 0, \forall \xi_{0} >0, \hspace{0.1 cm} K_{x_{0},R}^{k}\mathbf{1}(\xi_{0}) \xrightarrow[k \rightarrow + \infty]{} 0\bigg).$
\end{center} 
\label{lemme:kixzero}
\end{lemme}
\begin{proof}
We prove~\eqref{eq:probakxo} by induction on $k \geq 1$. The base case is given by Lemma~\ref{lemme:lemmatwo}, because $M^{1} = B$. Now we suppose that the result holds true for some $k \geq 1$. We define $(\tilde{\xi}_{.,x_{0},R,x})_{x >0}$ a family of processes with same law as $\xi_{.,x_{0},R,x}$, but independent from our initial process $\xi_{.,x_{0},R,\xi_{0}}$, and in particular independent from $M^{1}_{x_{0},R,\xi_{0}} $. We denote by $\tilde{M}^{k}_{x_{0},R,x}$ the associated events, defined like $M^{k}_{x_{0},R,\xi_{0}}$ for $\xi_{.,x_{0},R,\xi_{0}}$. By Markov property, we have that
$$\mathbb{P}(M^{k+1}_{x_{0},R,\xi_{0}}) = \mathbb{P}(M^{1}_{x_{0},R,\xi_{0}} \cap \tilde{M}^{k}_{x_{0},R,A_{\xi_{0},R}(J_{1})-x_{0}}). $$
Conditioning on the value of $J_{1}$ (which is in $]0,t_{\mathrm{max}}(\xi_{0})[$ on the event $M^{1}$) and using~\eqref{eq:jun} like in the proof of Lemma~\ref{lemme:lemmatwo}, we obtain 
\begin{align*}
\mathbb{P}(M^{k+1}_{x_{0},R,\xi_{0}})& = \displaystyle{\int_{0}^{t_{\mathrm{max}}(\xi_{0},R)}b_{x_{0}}(A_{\xi_{0},R}(t))e^{- \displaystyle{\int_{0}^{t} (b_{x_{0}}+d)(A_{\xi_{0},R}(s)) \mathrm{d}s }} \mathbb{P}(M^{k}_{x_{0},R,A_{\xi_{0},R}(t)-x_{0}}) \mathrm{d}t },
\end{align*}
because $\tilde{\xi}_{.,x_{0},R,x}$ and $\xi_{.,x_{0},R,x}$ have the same law. This concludes by induction hypothesis. The equivalence of Lemma~\ref{lemme:kixzero} follows from 
\begin{align*}
\{ N_{x_{0},R,\xi_{0}} = + \infty \} = \bigcap_{k \geq 1}M^{k}_{x_{0},R,\xi_{0}}. 
\end{align*} 
\end{proof}
\noindent We highlight the following consequence of Lemmas~\ref{lemme:premierr} and \ref{lemme:kixzero}, which implies point \textbf{3.} in Section~\ref{theo:sharp}.
\begin{corr}
Under the general setting of Section~\ref{subsec:gensett}, 
\begin{center}
Assumption~\ref{hyp:tempsinf}  $\Rightarrow$ Assumption~\ref{ass:pasinfini}.
\end{center}
Under the allometric setting of Section~\ref{subsec:allomsett}, we also have
\begin{center}
Assumptions~\ref{hyp:tempsinf} and \ref{ass:gainenergy}  $\Rightarrow$ $\max(\beta,\delta) \geq \alpha-1$.
\end{center}
\label{corr:swagosss}
\end{corr}
\begin{proof}
For the first point of the corollary, thanks to Lemma~\ref{lemme:kixzero}, the only thing to verify is that for any $x_{0}>0$, $\mathbb{P}(N_{x_{0},R,\xi_{0}} < + \infty) =1$, for every $R \notin \mathfrak{R}_{\infty}$, or $R \in \mathfrak{R}_{\infty}$ and $\xi_{0}>0$ such that there exists $x_{2}>\xi_{0}$, with $g(x_{2},R) \leq 0$ (\textit{i.e.} $\xi_{0} \notin \Omega_{R}$). Remark that $R \notin \mathfrak{R}_{\infty}$ implies the existence of such $x_{2}$ for any $\xi_{0}>0$. We work with such parameters in the following, and remark that in both cases, $T_{\infty}=+\infty$, so $\mathfrak{J}=+ \infty$ by Lemma~\ref{lemme:junfini}. First, if $R \notin \mathfrak{R}_{0}$, $T_{d} < + \infty$ by Lemma \ref{lemme:premierr} so $N_{x_{0},R,\xi_{0}} < + \infty$ by definition, see \eqref{eq:nombre}. Then if $R \in \mathfrak{R}_{0}$, there exists $x_{1}>0$ such that $g(y,R)<0$ if $y \leq x_{1}$. We define
$$ \tau := \inf \{ t \geq 0, \xi_{t,x_{0},R,\xi_{0}} \leq x_{0} \wedge x_{1} \},$$
with the convention $\inf(\varnothing) = + \infty$. We have again two cases. First if $\tau = + \infty$, then $\xi_{t,x_{0},R,\xi_{0}} $ stays in $[x_{0}\wedge x_{1},x_{2}] $ before $T_{d}$, and the death rate is positive continuous on this segment, so $T_{d}< +\infty$, hence $N_{x_{0},R,\xi_{0}} < + \infty$. Else if $\tau < + \infty$, there are no possible births after time $\tau$ because individual energy remains under $x_{0}$, so $N_{x_{0},R,\xi_{0}} < + \infty$.
\\\\
Under the allometric setting of Section~\ref{subsec:allomsett} and Assumption~\ref{ass:gainenergy}, according to Proposition~\ref{eq:gammalpha}, $\gamma=\alpha$ and $\mathfrak{R}_{\infty} \neq \varnothing$. We take $R$ in this set and then for all $y >0$, $g(y,R)>0$. The integral condition of Assumption~\ref{hyp:tempsinf} then leads to the second point.
\end{proof}
Now, we give another consequence of Assumption~\ref{hyp:tempsinf}. 
\begin{prop}
Under the general setting of Section~\ref{subsec:gensett},
\begin{align*}
Assumption \hspace{0.1 cm} \textit{\ref{hyp:tempsinf}} & \Rightarrow (\forall x_{0} >0, \forall
R \in \mathfrak{R}_{\infty} \setminus \mathfrak{R}_{0}, \forall \xi_{0} > 0, \quad \mathbb{P}_{x_{0}, R, \xi_{0}}( T_{d} < + \infty ) = 1).
\end{align*}
\label{prop:cns}
\end{prop}
\begin{proof}
We fix any $x_{0}>0$, $R \in \mathfrak{R}_{\infty} \setminus \mathfrak{R}_{0}$, $\xi_{0} >0 $, and we compute $\mathbb{P}_{x_{0}, R, \xi_{0}}(T_{d} < + \infty)$. In the following, we write $p_{k} : x >0 \mapsto \mathbb{P}_{x_{0},R,x}(\{ T_{d} = J_{k} \} \cap \{T_{d} < + \infty \})$. On the event $\{T_{d} < + \infty \}$, $T_{d}$ is one of the $J_{k}$, so we immediately have
\begin{align}
\mathbb{P}_{x_{0}, R, \xi_{0}}(T_{d} < + \infty) = \sum\limits_{k \geq 1} p_{k}(\xi_{0}) .
\label{eq:sommeproba}
\end{align} 
First, we have from the construction of Section~\ref{subsec:construc} that
\begin{align*}
p_{1}(\xi_{0}) = \mathbb{P}_{x_{0},R,\xi_{0}}(T_{d} = J_{1}, T_{d} < + \infty) & = \mathbb{P}_{x_{0},R,\xi_{0}}(D) = \mathbb{P}_{x_{0},R,\xi_{0}}(J_{1}<t_{\mathrm{max}}(\xi_{0})) -  \mathbb{P}_{x_{0},R,\xi_{0}}(B),
\end{align*}
which gives thanks to~\eqref{eq:jun} and Lemma~\ref{lemme:lemmatwo},
\begin{align}
p_{1}(\xi_{0}) & =  1 - \sigma_{x_{0},R,\xi_{0}} - K_{x_{0},R}\mathbf{1}(\xi_{0}),
\label{eq:omegaun}
\end{align}
with $\sigma_{x_{0},R,\xi_{0}}:= \exp\left(-\displaystyle{ \int_{0}^{t_{\max}(\xi_{0},R)}(b_{x_{0}}+d)(A_{\xi_{0},R}(\tau))\mathrm{d}\tau}\right)$. Then, for $k \geq 1$, we use Markov property as in the proof of Lemma~\ref{lemme:kixzero}, but with the events $\{ T_{d} = J_{k} \} \cap \{T_{d} < + \infty \}$ instead of $M^{k}_{x_{0},R,\xi_{0}}$, to obtain 
\begin{align*}
p_{k+1}(\xi_{0})& = \displaystyle{\int_{0}^{t_{\max}(\xi_{0},R)}b_{x_{0}}(A_{\xi_{0},R}(u))e^{-\displaystyle{ \int_{0}^{u}(b_{x_{0}}+d)(A_{\xi_{0},R}(\tau))\mathrm{d}\tau}}} p_{k}(A_{\xi_{0},R}(u)-x_{0}) \mathrm{d}u = K_{x_{0},R}p_{k}(\xi_{0}),
\end{align*}
and this is valid for every $\xi_{0}>0$ and $k \geq 1$, so 
\begin{align}
p_{k+1} = K_{x_{0},R}p_{k} = K_{x_{0},R}^{k}p_{1}.
\label{eq:formulerec}
\end{align} 
Now, thanks to Equations~\eqref{eq:sommeproba},~\eqref{eq:omegaun},~\eqref{eq:formulerec}, we obtain
\begin{align*}
\mathbb{P}_{x_{0}, R, \xi_{0}}(T_{d} < + \infty) & = \sum\limits_{k \geq 1}K_{x_{0},R}^{k-1}(\mathbf{1} - \sigma_{x_{0},R,.} - K_{x_{0},R}\mathbf{1})(\xi_{0}) \\
& = \sum\limits_{k \geq 1}K_{x_{0},R}^{k-1}(\mathbf{1} - K_{x_{0},R}\mathbf{1})(\xi_{0}) - \sum\limits_{k \geq 1}K_{x_{0},R}^{k-1}\sigma_{x_{0},R,.}(\xi_{0})\\ & = 1 - \lim_{k \rightarrow + \infty} K_{x_{0},R}^{k}\mathbf{1}(\xi_{0}) -\sum\limits_{k \geq 1}K_{x_{0},R}^{k-1}\sigma_{x_{0},R,.}(\xi_{0}).
\end{align*} 
We can split the sums from the first line to the second, because one of them is telescopic (and $\lim_{k \rightarrow + \infty} K_{x_{0},R}^{k}\mathbf{1}(\xi_{0})$ exists as the limit of a decreasing non-negative sequence) and the first line is the wanted probability. Now, if $\xi_{0} \in \Omega_{R}$, then $g(x,R)>0$ for $x \geq \xi_{0}$, and we can make the change of variables $x=A_{\xi_{0},R}(\tau)$ in the expression of $\sigma_{x_{0},R,\xi_{0}}$, with $A_{\xi_{0},R}(\tau)$ going to $+ \infty$ when $\tau$ goes to $t_{\max}(\xi_{0},R)$, so \begin{align*}
\sigma_{x_{0},R,\xi_{0}}:= \exp\left(-\displaystyle{ \int_{\xi_{0}}^{+ \infty}\dfrac{(b_{x_{0}}+d)(x)}{g(x,R)}\mathrm{d}x}\right).
\end{align*}
Remark that for $\xi_{0} \in \Omega_{R}$, if $\sigma_{x_{0},R,\xi_{0}} = 0$, then $\sigma_{x_{0},R,x} = 0$ for every $x \in \Omega_{R}$, by continuity of $b_{x_{0}}$, $d$ and $g(.,R)$. Also, if $\xi_{0} \notin \Omega_{R}$, recall that $R \notin \mathfrak{R}_{0}$, so $t_{\max}(\xi_{0},R) = + \infty$ and $A_{\xi_{0},R}(\tau)$ remains in a compact for $\tau \geq 0$, hence $\sigma_{x_{0},R,\xi_{0}}=0$. This implies that for every $\xi_{0} \in \Omega_{R}$, ($\sigma_{x_{0},R,\xi_{0}} = 0 \Leftrightarrow \sigma_{x_{0},R,.} \equiv 0$).
Thus, for $\xi_{0} \in \Omega_{R}$, $(\mathbb{P}_{x_{0},R,\xi_{0}}(T_{d}=+ \infty)=1)$ is equivalent to ($\sigma_{x_{0},R,\xi_{0}} = 0$ and $\lim_{k \rightarrow + \infty} K_{x_{0},R}^{k}\mathbf{1}(\xi_{0})=0$). Assuming Assumption~\ref{hyp:tempsinf} implies that ($\sigma_{x_{0},R,\xi_{0}} = 0$ and $\lim_{k \rightarrow + \infty} K_{x_{0},R}^{k}\mathbf{1}(\xi_{0})=0$) for every $R \in~ \mathfrak{R}_{\infty} \setminus \mathfrak{R}_{0}$ and $\xi_{0} \in \Omega_{R}$. Also, if $R \in \mathfrak{R}_{\infty} \setminus \mathfrak{R}_{0}$ and $\xi_{0} \notin \Omega_{R}$, then $\mathbb{P}_{x_{0},R,\xi_{0}}(T_{d}=+ \infty)=1$ by Lemma~\ref{lemme:premierr}, which ends the proof.
\end{proof}

Finally, we complete the proof of Theorem~\ref{theo:cns}. 

\begin{prop}
Under the general setting of Section~\ref{subsec:gensett}, and under Assumption~\ref{hyp:probamort}, we have
\begin{align*}
Assumption \hspace{0.1 cm} \textit{\ref{hyp:tempsinf}} & \Leftrightarrow (\forall x_{0} >0, \forall
R \in \mathfrak{R}_{\infty}, \forall \xi_{0} \in \Omega_{R}, \quad \mathbb{P}_{x_{0}, R, \xi_{0}}( T_{d} < + \infty ) = 1).
\end{align*}
Also, under Assumption~\ref{hyp:probamort}, for every $x_{0}>0$, $R \in \mathfrak{R}_{\infty}$ and $\xi_{0} \notin \Omega_{R}$, we have $\mathbb{P}_{x_{0}, R, \xi_{0}}( T_{d} < + \infty ) = 1$.
\label{prop:tinfyinfini}
\end{prop}
\begin{proof}
We fix $R \in \mathfrak{R}_{\infty}$ in the following. By Assumption~\ref{hyp:probamort} and Proposition~\ref{prop:toinfini}, we have almost surely $T_{0}=+ \infty$. We can work on this event and apply the exact same technique as in the proof of Proposition~\ref{prop:cns} (in particular, we verify that we still have $\sigma_{x_{0},R,\xi_{0}} = 0$ for every $\xi_{0} \notin \Omega_{R}$) to prove that $(\mathbb{P}_{x_{0},\xi_{0},R}(T_{d}=+ \infty)=1)$ is equivalent to ($\sigma_{x_{0},R,\xi_{0}} = 0$ and $\lim_{k \rightarrow + \infty} K_{x_{0},R}^{k}\mathbf{1}(\xi_{0})=0$) for every $R \in \mathfrak{R}_{\infty}$ and $\xi_{0} \in \Omega_{R}$, which concludes for the equivalence. Finally, if $\xi_{0} \notin \Omega_{R}$, it is as if we work with $R \notin \mathfrak{R}_{\infty}$, and under Assumption~\ref{hyp:probamort}, we can use the same techniques as in the proofs of Proposition~\ref{eq:delta} or Lemma~\ref{lemme:premierr} to obtain the second part of Proposition~\ref{prop:tinfyinfini}.
\end{proof}
\textbf{Remark:} What is important for the equivalence in the previous proof is the use of Assumption~\ref{hyp:probamort} first to ensure that $T_{0}=+ \infty$ almost surely, even if $R \in \mathfrak{R}_{0}$. Also, we need to consider $\xi_{0} \in \Omega_{R}$ to perform the change of variables $x=A_{\xi_{0},R}(\tau)$ in the proof of Proposition~\ref{prop:cns}, and to ensure $t_{\mathrm{\max}}(\xi_{0},R)=+ \infty$. At this point, the reader can verify that the combination of Proposition~\ref{eq:delta}, Lemma~\ref{lemme:premierr} and Proposition~\ref{prop:tinfyinfini} implies Theorem~\ref{theo:cns}.
\subsection{Proof of 4. in Section~\ref{theo:sharp}}
\label{proof:super}

A first immediate consequence of Assumption~\ref{ass:supercritical} is the following, which entails \textbf{4.} in Section~\ref{theo:sharp}.
\begin{prop}
Under the general setting of Section~\ref{subsec:gensett}, Assumptions~\ref{hyp:probamort}, \ref{hyp:tempsinf} and \ref{ass:supercritical}, we have 
\begin{align}
\forall x >0, \exists y \geq x, \quad \tilde{b}(y) > d(y).
\label{eq:bsupd}
\end{align}  
Under the allometric setting of Section~\ref{subsec:allomsett},~\eqref{eq:bsupd} is equivalent to
\[ \beta > \delta \hspace{0.1 cm} \mathrm{or} \hspace{0.1 cm} (\beta= \delta, \hspace{0.1 cm} C_{\beta} > C_{\delta}). \]
\label{lemm:lemmaone}
\end{prop}
\begin{proof}
Let us suppose that
\[ \exists x >0, \forall y \geq x, \quad \tilde{b}(y) \leq d(y). \]
Then, recall that $b_{x_{0}} \equiv 0$ on $]0,x_{0}]$, so for every $x_{0} \geq x$,
\[ \forall y>0, \hspace{0.1 cm} \dfrac{b_{x_{0}}(y)}{b_{x_{0}}(y)+d(y)} \leq 1/2.\]  
But in that case, no matter the choice of $R$ and the instant at which the individual energy jumps, the probability of giving birth instead of dying will be lower than or equal to $1/2$. In other terms, we have $$\forall x_{0} \geq x, \forall R \geq 0 , \hspace{0.1cm} m_{x_{0},R}(x_{0}) \leq 1,$$ which contradicts Assumption~\ref{ass:supercritical}, by Proposition~\ref{prop:galtonwatson}. The allometric consequence is straightforward.  
\end{proof}
\noindent \textbf{Remark:} It is natural to compare the birth and death rates to conclude about the supercriticality of the process. Informally, Proposition~\ref{lemm:lemmaone} says that the birth rate should dominate the death rate, at least in the high energy regime. It seems biologically natural to assume that a well-fed individual is more likely to give birth to offspring than to die. 
\\\\
An important difficulty in the sequel is to conclude for small $x_{0}$. It could be possible to have a lower birth rate than the death rate when the energy is low, but higher if the individual energy increases (as it is the case in the allometric setting of Section~\ref{subsec:allomsett} if $\beta > \delta$). Then, we have to control the probability for individuals to reach the favorable high energies, to know better about the expected number of offspring starting from $x_{0}$, which is $m_{x_{0},R}(x_{0})$. This is the main goal of the next subsections.

\subsection{Proof of 5. in Section~\ref{theo:sharp}}
\label{proof:beplusde}
In this section, we first introduce a coupling between the processes $\xi_{.,x_{0},R,\xi_{0}}$ for different values of $x_{0}$ and/or $\xi_{0}$, in the same environment $R$. The main result of this section is Proposition~\ref{prop:beplusde}, valid under the general setting of Section~\ref{subsec:gensett}. Then, we apply the same techniques in the allometric setting of Section~\ref{subsec:allomsett} and the particular case $\delta= \alpha-1$, to obtain \textbf{5.} in Section~\ref{theo:sharp}.

\subsubsection{Preliminary results}

In the following, for $x>0$, we say that $R \in \mathcal{R}_{x}$ if for all $y \geq x$, $g(y,R) >0$ (note that $R \in \mathfrak{R}_{\infty}$, if and only if there exists $x >0$ such that $R \in \mathcal{R}_{x}$). We naturally extend this notation to $\mathcal{R}_{0} :=  \{ R \geq 0, \forall y>0, g(y,R) >0 \}$ (be careful that $\mathcal{R}_{0}$ is distinct from the previous notation $\mathfrak{R}_{0}$). These sets are possibly equal to $\varnothing$ under the general setting~\ref{subsec:gensett}, and in that case most of our following results do not apply. However, the following lemma states that under the allometric setting of Section~\ref{subsec:allomsett} and Assumption~\ref{ass:gainenergy}, $\mathcal{R}_{0} \neq \varnothing$.

\begin{lemme}
Under the general setting of Section~\ref{subsec:gensett}, for $x > 0$, if there exists $R >0$ such that $R \in \mathcal{R}_{x}$ then
$$[R, + \infty [ \subseteq  \mathcal{R}_{x}.$$
In addition, under the allometric setting of Section~\ref{subsec:allomsett} and Assumption~\ref{ass:gainenergy}, there exists $R_{0}>0$ such that  
$$]R_{0}, + \infty [ =  \mathcal{R}_{0}.$$
\label{cestlafindesharicots}
\end{lemme}
\begin{proof}
The first point comes from the fact that $g(x,.)$ is increasing on $\mathbb{R}_{+}$ for any $x>0$. The allometric consequence is straightforward by Proposition~\ref{eq:gammalpha} and considering Equation~\eqref{eq:indivenergymod}.
\end{proof}
\noindent In the following, we will use the notation $B_{x_{0},R,\xi_{0}}$ from Section~\ref{subsec:teun}. In this section, we present our results with the assumption that $d/\tilde{b}$ is non-increasing. In fact, under the allometric setting of Section~\ref{subsec:allomsett}, Assumptions~\ref{hyp:probamort}, \ref{hyp:tempsinf} and \ref{ass:supercritical}, this assumption is a natural consequence of Proposition~\ref{lemm:lemmaone}. The following lemma highlights a useful coupling between the processes $\xi_{.,x_{0},R,\xi_{0}}$ for different values of $x_{0}$ and/or $\xi_{0}$ and the same well-chosen $R$.

\begin{lemme}
We work under the general setting of Section \ref{subsec:gensett}. We suppose Assumptions~\ref{hyp:probamort} and~\ref{hyp:tempsinf}, and that $d/\tilde{b}$ is non-increasing. Let $(x_{0}, \hat{x}_{0}, \xi_{0}, \hat{\xi}_{0}) >0$ and $R \in \mathcal{R}_{\xi_{0} \wedge \hat{\xi}_{0}}$. Let $\xi_{.,x_{0},R,\xi_{0}}$ and $\xi_{.,\hat{x}_{0},R,\hat{\xi}_{0}}$ be individual trajectories with law described in Section~\ref{subsec:construc}, respectively characterized by $x_{0}$ and $\hat{x}_{0}$, and respectively starting from $\xi_{0}$ and $\hat{\xi}_{0}$, in the same environment $R$. Then there exists a coupling of these random processes such that
\begin{itemize}
\item If $\max(x_{0},\xi_{0}) \geq \max(\hat{x}_{0},\hat{\xi}_{0})$, then on the event $B_{x_{0},R,\xi_{0}} \cap B_{\hat{x}_{0},R,\hat{\xi}_{0}}$, $$\xi_{J_{1}^{\hat{x}_{0},R,\hat{\xi}_{0}},\hat{x}_{0},R,\hat{\xi}_{0}} + \hat{x}_{0} \leq \xi_{J_{1}^{x_{0},R,\xi_{0}},x_{0},R,\xi_{0}} + x_{0},$$
where $J_{1}^{x_{0},R,\xi_{0}}$ and $J_{1}^{\hat{x}_{0},R,\hat{\xi}_{0}}$ denotes the first time of jump, respectively associated to $\xi_{.,x_{0},R,\xi_{0}}$ and $\xi_{.,\hat{x}_{0},R,\hat{\xi}_{0}}$.
\item If $\xi_{0} \geq \hat{\xi}_{0}$ and $\hat{x}_{0} \geq x_{0}$, then $$ B_{\hat{x}_{0},R,\hat{\xi}_{0}} \subseteq B_{x_{0},R,\xi_{0}}.$$
\end{itemize}
\label{lemm:coupling}
\end{lemme}
\begin{proof}
The proof is divided in three steps.
\\\\
\textbf{Step 1: construction of a coupling}
\\\\
First, we couple $\xi_{.,x_{0},R,\xi_{0}}$ and $\xi_{.,\hat{x}_{0},R,\hat{\xi}_{0}}$, and more generally define simultaneously all the $\xi_{.,x,R,y}$ for $x > 0$ and $y > 0$. The reader can check that this new definition of the individual process and the one of Section~\ref{subsec:construc} will have same distribution of sample paths. Indeed, instead of constructing jump times at inhomogeneous exponential rate depending on $b_{x_{0}}+d$, we will separate it into inhomogeneous exponential rates depending on $b_{x_{0}}$ and $d$ separately, and the fact that the jump times have the same law comes from the usual property for the law of the minimum of independent exponential random variables. The following construction is very similar to the one given in Section~\ref{subsec:construc}. We use it because it is convenient for our coupling purposes.
\\\\
We pick i.i.d.~random variables $(F_{i})_{i \geq 0}$, following exponential laws with parameter 1, and define auxiliary processes $(\xi^{\mathrm{aux}}_{t,x,R,y})_{t \geq 0}$, with random jumps occuring at the following times. 
\begin{itemize}
\item[-] First, if
$$\int_{0}^{t_{\mathrm{max}}(y,R)}b_{x}(A_{y,R}(s)) \mathrm{d}s \leq F_{1}, $$
we set $J_{1}^{x,R,y}:=+ \infty$, $\xi^{\mathrm{aux}}_{t,x,R,y} = A_{y,R}(t)$ for $t \in [0,t_{\mathrm{max}}(\xi_{0},R)[$ and $\xi^{\mathrm{aux}}_{t,x,R,y} = \flat$ for $t \geq t_{\mathrm{max}}(\xi_{0},R)$. Otherwise, we define the first time of jump as 
$$ J_{1}^{x,R,y} := \inf \left\{ t \in [0,t_{\max}(y,R)[, \int_{0}^{t}b_{x}(A_{y,R}(s)) \mathrm{d}s = F_{1} \right\}, $$
we set $\xi^{\mathrm{aux}}_{t,x,R,y} = A_{y,R}(t)$ for $t \in [0,J_{1}^{x,R,y}[$, and $$\xi^{\mathrm{aux}}_{J_{1}^{x,R,y},x,R,y} = A_{y,R}(J_{1}^{x,R,y}) - x .$$

\item[-] Now, let us suppose that we defined $J_{n}^{x,R,y}$ for some $n \geq 1$. If $J_{n}^{x,R,y}= + \infty$, we simply set $J_{n+1}^{x,R,y}:= + \infty$ and $\xi^{\mathrm{aux}}_{t,x,R,y}$ is already defined for all $t \geq 0$, it already reached $\flat$. Else if
$$\int_{J_{n}^{x,R,y}}^{J_{n}^{x,R,y} + t_{\mathrm{max}}\left(\xi^{\mathrm{aux}}_{J_{n}^{x,R,y},x,R,y}\right)}b_{x}\left(A_{\xi^{\mathrm{aux}}_{J_{n}^{x,R,y},x,R,y}}(s-J_{n}^{x,R,y})\right)\mathrm{d}s \leq F_{n+1}, $$
we also set $J_{n+1}^{x,R,y}:=+ \infty$, $\xi^{\mathrm{aux}}_{t,x,R,y} = A_{\xi^{\mathrm{aux}}_{J_{n}^{x,R,y},x,R,y}}(t-J_{n}^{x,R,y})$ for $t \in [J_{n}^{x,R,y},J_{n}^{x,R,y} +t_{\mathrm{max}}\left(\xi^{\mathrm{aux}}_{J_{n}^{x,R,y},x,R,y}\right)[$ and $\xi^{\mathrm{aux}}_{t,x,R,y} = \flat$ for $t \geq J_{n}^{x,R,y} +t_{\mathrm{max}}\left(\xi^{\mathrm{aux}}_{J_{n}^{x,R,y},x,R,y}\right)$. Otherwise, we can define the $(n+1)$-th time of jump as 
$$ J_{n+1}^{x,R,y} := \inf \left\{ t \in \left[J_{n}^{x,R,y},J_{n}^{x,R,y}+t_{\mathrm{max}}\left(\xi^{\mathrm{aux}}_{J_{n}^{x,R,y},x,R,y}\right)\right[, \int_{J_{n}}^{t}b_{x}\left(A_{\xi^{\mathrm{aux}}_{J_{n}^{x,R,y},x,R,y}}(s-J_{n}^{x,R,y})\right)\mathrm{d}s = F_{n+1} \right\}, $$
we set $\xi^{\mathrm{aux}}_{t,x,R,y} = A_{\xi^{\mathrm{aux}}_{J_{n}^{x,R,y},x,R,y}}(t-J_{n}^{x,R,y})$ for $t \in [J_{n}^{x,R,y},J_{n+1}^{x,R,y}[$, and 
\begin{align*}
\xi^{\mathrm{aux}}_{J_{n+1}^{x,R,y},x,R,y} = & A_{\xi^{\mathrm{aux}}_{J_{n}^{x,R,y},x,R,y}}(J_{n+1}^{x,R,y}-J_{n}^{x,R,y}) - x. 
\end{align*}
Notice that the jump rate of $\xi^{\mathrm{aux}}$ is $b_{x}$, so $\xi^{\mathrm{aux}}_{J_{n+1}^{x,R,y},x,R,y} > 0$. 
\end{itemize}
Then, we define the time of death: $$T_{d,x,R,y} := \inf \left\{ t \geq 0, \displaystyle{\int_{0}^{t} d(\xi^{\mathrm{aux}}_{s,x,R,y})\mathrm{d}s = F_{0}}\right\},$$
and finally set $\xi_{t,x,R,y} := \xi^{\mathrm{aux}}_{t,x,R,y}\mathbb{1}_{\{t < T_{d,x,R,y}\}} + \partial \mathbb{1}_{\{t \geq T_{d,x,R,y}\}} $ for all $t \geq 0$. Under Assumptions~\ref{hyp:probamort} and~\ref{hyp:tempsinf}, this process is almost surely biologically relevant thanks to Theorem~\ref{theo:cns}. 
\\\\
\textbf{Step 2: proof of the first point of the lemma.}
\\\\
Our coupling allows us to obtain immediately
\[ \displaystyle{\int_{0}^{J_{1}^{x_{0},R,\xi_{0}}} b_{x_{0}}(A_{\xi_{0}}(s))\mathrm{d}s} = F_{1} = \displaystyle{\int_{0}^{J_{1}^{\hat{x}_{0},R,\hat{\xi}_{0}}} b_{\hat{x}_{0}}(A_{\hat{\xi}_{0}}(s))\mathrm{d}s}. \]
Thanks to the choice of $R$, $g(x,R)>0$ if $x \geq \xi_{0} \wedge \hat{\xi}_{0}$, so we can use the change of variables $u=A_{\xi_{0}}(s)$, or $u=A_{\hat{\xi}_{0}}(s)$ and obtain
\begin{align*}
& \displaystyle{\int_{\xi_{0}}^{A_{\xi_{0}}(J_{1}^{x_{0},R,\xi_{0}})} \dfrac{b_{x_{0}}(u)}{g(u,R)} \mathrm{d}u} = \displaystyle{\int_{\hat{\xi}_{0}}^{A_{\hat{\xi}_{0}}(J_{1}^{\hat{x}_{0},R,\hat{\xi}_{0}})} \dfrac{b_{\hat{x}_{0}}(u)}{g(u,R)} \mathrm{d}u}, 
\end{align*}
so
\begin{align*}
  & \displaystyle{\int_{\max(x_{0},\xi_{0})}^{A_{\xi_{0}}(J_{1}^{x_{0},R,\xi_{0}})} \dfrac{\tilde{b}(u)}{g(u,R)} \mathrm{d}u} = \displaystyle{\int_{\max(\hat{x}_{0},\hat{\xi}_{0})}^{A_{\hat{\xi}_{0}}(J_{1}^{\hat{x}_{0},R,\hat{\xi}_{0}})} \dfrac{\tilde{b}(u)}{g(u,R)} \mathrm{d}u}. 
\end{align*} 
As $\max(x_{0},\xi_{0}) \geq \max(\hat{x}_{0},\hat{\xi}_{0})$ and $\tilde{b}(.)/g(.,R)$ is positive on the integration intervals, this enforces $A_{\xi_{0}}(J_{1}^{x_{0},R,\xi_{0}}) \geq A_{\hat{\xi}_{0}}(J_{1}^{\hat{x}_{0},R,\hat{\xi}_{0}})$. Thus, on the event $B_{x_{0},R,\xi_{0}} \cap B_{\hat{x}_{0},R,\hat{\xi}_{0}}$, the first point follows from the above definition of $\xi_{J_{1}^{x_{0},R,\xi_{0}},x_{0},R,\xi_{0}}$ and $\xi_{J_{1}^{\hat{x}_{0},R,\hat{\xi}_{0}},\hat{x}_{0},R,\hat{\xi}_{0}}$.
\\\\
\textbf{Step 3: proof of the second point of the lemma.}
\\\\
Remark that if we fix $x > 0$ and consider the function $y \in  [\xi_{0} \wedge \hat{\xi}_{0}, + \infty[ \mapsto \displaystyle{\int_{y}^{A_{y}(J_{1}^{x,R,y})} \dfrac{b_{x}(u)}{g(u,R)} \mathrm{d}u} $, it is constant, equal to $F_{1}$. One can check that this function is almost surely differentiable thanks to Theorem~\ref{theo:cns}, the choice of $R$, and the regularity of the flow, and the previous fact means that its derivative vanishes, so that for every $x > 0$ and $y \geq  \xi_{0} \wedge \hat{\xi}_{0}$:
\begin{align}
\dfrac{\partial A_{y}(J_{1}^{x,R,y})}{\partial y} \times \dfrac{b_{x}(A_{y}(J_{1}^{x,R,y}))}{g(A_{y}(J_{1}^{x,R,y}), R)} - \dfrac{b_{x}(y)}{g(y,R)} = 0.
\label{eq:derivb}
\end{align}  
Now, for $x > 0$ and $y \geq  \xi_{0} \wedge \hat{\xi}_{0}$, the event $B_{x,R,y}$ translates into
$$ \displaystyle{\int_{0}^{J_{1}^{x,R,y}} d(A_{y}(s))\mathrm{d}s} < F_{0}, $$
because $A_{y}(.)$ and $\xi_{.,x,R,y}$ coincide before $J_{1}^{x,R,y}$. With the change of variables $u=A_{y}(s)$, we obtain
$$ \displaystyle{\int_{y}^{A_{y}(J_{1}^{x,R,y})} \dfrac{d(u)}{g(u,R)} \mathrm{d}u} < F_{0}.$$
Hence, to prove the second point, it suffices to establish that $$\displaystyle{\int_{\xi_{0}}^{A_{\xi_{0}}(J_{1}^{x_{0},R,\xi_{0}})} \dfrac{d(u)}{g(u,R)} \mathrm{d}u} \leq \displaystyle{\int_{\hat{\xi}_{0}}^{A_{\hat{\xi}_{0}}(J_{1}^{\hat{x}_{0},R,\hat{\xi}_{0}})} \dfrac{d(u)}{g(u,R)} \mathrm{d}u}.  $$ 
As $\hat{x}_{0} \geq x_{0}$, we have $b_{x_{0}}(A_{\xi_{0}}(s)) \geq b_{\hat{x}_{0}}(A_{\xi_{0}}(s))$ for $s \geq 0$, so $J_{1}^{x_{0},R,\xi_{0}} \leq J_{1}^{\hat{x}_{0},R,\xi_{0}}.$ The flow is increasing because $g(.,R)$ is positive so $$ A_{\xi_{0}}(J_{1}^{x_{0},R,\xi_{0}}) \leq A_{\xi_{0}}(J_{1}^{\hat{x}_{0},R,\xi_{0}}), $$
and $d(.)/g(.,R)$ is positive so it suffices to prove that
$$\displaystyle{\int_{\xi_{0}}^{A_{\xi_{0}}(J_{1}^{\hat{x}_{0},R,\xi_{0}})} \dfrac{d(u)}{g(u,R)} \mathrm{d}u} \leq \displaystyle{\int_{\hat{\xi}_{0}}^{A_{\hat{\xi}_{0}}(J_{1}^{\hat{x}_{0},R,\hat{\xi}_{0}})} \dfrac{d(u)}{g(u,R)} \mathrm{d}u}.$$ 
Finally, $\xi_{0} \geq \hat{\xi}_{0}$, so it suffices to show that $y \in  [\xi_{0} \wedge \hat{\xi}_{0},+ \infty[ \mapsto \displaystyle{\int_{y}^{A_{y}(J_{1}^{\hat{x}_{0},R,y})} \dfrac{d(u)}{g(u,R)} \mathrm{d}u}$ is non-increasing. Again, this function is differentiable, so we prove that 
$$ \forall y \geq  \xi_{0} \wedge \hat{\xi}_{0}, \quad \dfrac{\partial A_{y}(J_{1}^{\hat{x}_{0},R,y})}{\partial y} \times \dfrac{d(A_{y}(J_{1}^{\hat{x}_{0},R,y}))}{g(A_{y}(J_{1}^{\hat{x}_{0},R,y}),R)} - \dfrac{d(y)}{g(y,R)} \leq 0. $$
Equation \eqref{eq:derivb} leads to
\begin{align*}
\dfrac{\partial A_{y}(J_{1}^{\hat{x}_{0},R,y})}{\partial y} \times \dfrac{d(A_{y}(J_{1}^{\hat{x}_{0},R,y}))}{g(A_{y}(J_{1}^{\hat{x}_{0},R,y}),R)} - \dfrac{d(y)}{g(y,R)} & = \dfrac{b_{\hat{x}_{0}}(y)}{g(y,R)} \times \dfrac{d(A_{y}(J_{1}^{\hat{x}_{0},R,y}))}{b_{\hat{x}_{0}}(A_{y}(J_{1}^{\hat{x}_{0},R,y}))} - \dfrac{d(y)}{g(y,R)}
\end{align*}   
This is equal to $\dfrac{\tilde{b}(y)}{g(y,R)} \bigg( \dfrac{d(A_{y}(J_{1}^{\hat{x}_{0},R,y}))}{\tilde{b}(A_{y}(J_{1}^{\hat{x}_{0},R,y}))} - \dfrac{d(y)}{\tilde{b}(y)} \bigg)$ if $y > \hat{x}_{0}$, and to $- \dfrac{d(y)}{g(y,R)}$ otherwise.
\\\\
In both cases, this quantity is non-positive because $d/\tilde{b}$ is non-increasing, $\tilde{b}(.)/g(.,R)$ and $d(.)/g(.,R)$ are non-negative, which ends the proof. 
\end{proof}
\textbf{Remark:} The reader can check that the definition of the individual process $\xi$ in the proof of Lemma~\ref{lemm:coupling} and the one of Section~\ref{subsec:construc} have same distribution of sample paths, thanks to Lemma~\ref{lemme:split} of Section~\ref{subsec:laprevue}.
\\\\
We can still define the number of offspring of an individual during its life as in~\eqref{eq:nombre}: $$N_{x_{0},R,\xi_{0}} := \sup \{ i \geq 0, J_{i}^{x_{0},R,\xi_{0}} < T_{d,x_{0},R,\xi_{0}} \}, $$
and $m_{x_{0},R}(\xi_{0}) := \mathbb{E}(N_{x_{0},R,\xi_{0}})$. In the following, consistent with these notations, we write $m_{0,R}(\xi_{0})$ for the expected number of birth jumps of an individual starting from $\xi_{0}$, but losing no energy when a birth occurs. Indeed, the previous coupling makes perfect sense for $x_{0}=0$ and we implicitly extend it to this remaining value in the following. We insist here on the fact that it makes sense to work with an individual process with $x_{0}=0$, but we cannot relate it to a population process anymore, because we do not want to see offspring with energy 0 appear in the population. We continue to use the expression ``birth jump'', but in the case $x_{0}=0$, the reader should be aware that this does not relate to any underlying branching process. Lemma~\ref{lemm:coupling} still apply for $x_{0}=0$, and we use it to stochastically dominate $m_{\hat{x}_{0},R}(.)$ by $m_{0,R}(.)$ for $\hat{x}_{0}>0$ in the upcoming lemma.

\begin{lemme}
We work under the general setting of Section \ref{subsec:gensett}. We suppose Assumptions~\ref{hyp:probamort} and~\ref{hyp:tempsinf}, and that $d/\tilde{b}$ is non-increasing. We fix $R \in \mathcal{R}_{0}$, and $(\hat{x}_{0}, \xi_{0}, \hat{\xi}_{0})$, such that $\hat{x}_{0} > 0$ and $\xi_{0} \geq \max(\hat{x}_{0},\hat{\xi}_{0})$. Then, under the same coupling as in Lemma~\ref{lemm:coupling}, we have
$$ N_{0,R,\xi_{0}} \geq N_{\hat{x}_{0},R,\hat{\xi}_{0}}.$$
Taking the expectation immediately leads to
$$ m_{0,R}(\xi_{0}) \geq m_{\hat{x}_{0},R}(\hat{\xi}_{0}). $$
\label{prop:super}
\end{lemme}

\begin{proof}
The second point of Lemma~\ref{lemm:coupling} allows us to assess that if an individual starting from $(\hat{x}_{0},\hat{\xi}_{0})$ has a child, then an individual starting from $(0,\xi_{0})$ has one too. Moreover, the first point of Lemma~\ref{lemm:coupling} shows that the new starting energies after the birth verify
$$\xi_{J_{1}^{\hat{x}_{0},R,\hat{\xi}_{0}},\hat{x}_{0},R,\hat{\xi}_{0}} + \hat{x}_{0} \leq \xi_{J_{1}^{0,R,\xi_{0}},0,R,\xi_{0}}.$$
This leads to $\xi_{J_{1}^{0,R,\xi_{0}},0,R,\xi_{0}} \geq \max(\hat{x}_{0}, \xi_{J_{1}^{\hat{x}_{0},R,\hat{\xi}_{0}},\hat{x}_{0},R,\hat{\xi}_{0}})$, so the vector $(\hat{x}_{0},\xi_{J_{1}^{0,R,\xi_{0}},0,R,\xi_{0}},\xi_{J_{1}^{\hat{x}_{0},R,\hat{\xi}_{0}},\hat{x}_{0},R,\hat{\xi}_{0}})$ satisfies the assumptions of Lemma~\ref{prop:super}. Thanks to Markov property, we can apply again the previous reasoning and obtain the same result for the second time of birth. This leads by induction to the fact that for every $i\geq 1$, if an individual starting from $(\hat{x}_{0},\hat{\xi}_{0})$ has $i$ children, then an individual starting from $(0,\xi_{0})$ has at least $i$ children. This gives the wanted conclusion.
\end{proof}
\noindent Finally, we use this lemma to prove the following proposition.
\begin{prop}
\label{prop:beplusde}
We work under the general setting of Section~\ref{subsec:gensett}, with Assumptions~\ref{hyp:probamort},~\ref{hyp:tempsinf} and~\ref{ass:supercritical}. Also, we assume that $d/\tilde{b}$ is non-increasing, and that $R \in \mathcal{R}_{0}$ satisfies that $d/g(.,R)$ has a divergent integral near $+ \infty$. Then, we necessarily have
\begin{align}
\forall x>0, \quad \displaystyle{\int_{x}^{+\infty}} \bigg(\displaystyle{\int_{x}^{u} \dfrac{\tilde{b}(s)}{g(s,R)} \mathrm{d}s}\bigg)  \dfrac{d(u)}{g(u,R)}\mathrm{d}u  = + \infty.
\label{eq:couplingeq}
\end{align} 
Under the allometric setting of Section~\ref{subsec:allomsett}, with Assumptions~\ref{hyp:probamort},~\ref{hyp:tempsinf} and~\ref{ass:supercritical}, Equation~\eqref{eq:couplingeq} implies that if $\delta \geq \alpha -1$, then
\[ \beta + \delta \geq 2(\alpha -1).\]
\end{prop}
\begin{proof}
We take $x >0$ and $R \in \mathcal{R}_{0}$. Lemma~\ref{prop:super} gives in particular that
\begin{align} 
m_{x,R}(x) \leq m_{0,R}(x).
\label{eq:comparem}
\end{align}
Without any loss of energy when a birth occurs, we can compute $m_{0,R}(x)$. The resource $R$ has been chosen so that $$ \displaystyle{\int_{x}^{+\infty}} \dfrac{d(y)}{g(y,R)} \mathrm{d}y = + \infty,$$ so from the coupling of Lemma~\ref{lemm:coupling} with $x_{0}=0$, we have $T_{d}<+ \infty$ almost surely. We know the law of $T_{d}$ when $x_{0}=0$, it is an inhomogeneous exponential law with rate depending on $d$ and the flow $A_{x}(.)$. Conditionally to the time of death $T_{d}$, the number of jumps before $T_{d}$ is then Poisson with parameter $\int_{0}^{T_{d}}\tilde{b}(A_{x}(w))\mathrm{d}w$. Hence:
\begin{align}
m_{0,R}(x) & = \displaystyle{\int_{0}^{+\infty}\bigg(\int_{0}^{t}\tilde{b}(A_{x}(w))\mathrm{d}w\bigg)d(A_{x}(t))e^{-\int_{0}^{t}d(A_{x}(\tau))\mathrm{d}\tau}\mathrm{d}u}. 
\label{eq:enfantssansnmort}
\end{align}
From the change of variables $s=A_{x}(w)$, then $u=A_{x}(t)$, we obtain
\begin{align}
m_{0,R}(x) & = \displaystyle{\int_{x}^{+\infty} \bigg(\displaystyle{\int_{x}^{u} \dfrac{\tilde{b}(s)}{g(s,R)} \mathrm{d}s}\bigg)  \dfrac{d(u)}{g(u,R)}e^{-\displaystyle{\int_{x}^{u}\dfrac{d(\tau)}{g(\tau,R)}\mathrm{d}\tau}}}\mathrm{d}u  \leq \displaystyle{\int_{x}^{+\infty} \bigg(\displaystyle{\int_{x}^{u} \dfrac{\tilde{b}(s)}{g(s,R)} \mathrm{d}s}\bigg)  \dfrac{d(u)}{g(u,R)}\mathrm{d}u}.
\label{eq:majorantm}
\end{align}
Let us suppose by contradiction that the right-most integral above is finite for some $x >0$, then \begin{align}
\displaystyle{\int_{y}^{+\infty} \bigg(\displaystyle{\int_{x}^{u} \dfrac{\tilde{b}(s)}{g(s,R)} \mathrm{d}s}\bigg)  \dfrac{d(u)}{g(u,R)}\mathrm{d}u} \xrightarrow[y\to +\infty]{} 0,
\label{eq:uno}
\end{align}
as the remainder of a convergent integral. For $y>x$, as $\tilde{b}(.)/g(.,R)$ is positive for our choice of $R$, we also have \begin{align}
\displaystyle{\int_{y}^{+\infty} \bigg(\displaystyle{\int_{y}^{u} \dfrac{\tilde{b}(s)}{g(s,R)} \mathrm{d}s}\bigg)  \dfrac{d(u)}{g(u,R)}\mathrm{d}u} \leq \displaystyle{\int_{y}^{+\infty} \bigg(\displaystyle{\int_{x}^{u} \dfrac{\tilde{b}(s)}{g(s,R)} \mathrm{d}s}\bigg)  \dfrac{d(u)}{g(u,R)}\mathrm{d}u}.
\label{eq:deuzio}
\end{align}
From~\eqref{eq:comparem},~\eqref{eq:majorantm},~\eqref{eq:deuzio} and finally~\eqref{eq:uno}, we would get $$ m_{y,R}(y) \leq 1$$ for $y$ high enough. One can check that we still get this result for this specific $y$, if we take a higher $R$ because $g(y,.)$ is non-decreasing. This contradicts Assumption~\ref{ass:supercritical} by Proposition~\ref{prop:galtonwatson} and ends the proof.
\\\\
Under the allometric setting of Section \ref{subsec:allomsett} and Assumptions~\ref{hyp:probamort},~\ref{hyp:tempsinf} and~\ref{ass:supercritical}, the assumption $\delta \geq \alpha-1$ ensures that the integral of $d/g(.,R)$ diverges near $ +\infty$, because with point \textbf{1.} of Section~\ref{theo:sharp}, we have $\gamma=\alpha$. Then, it follows from computations that if $\beta + \delta < 2(\alpha-1)$, then \eqref{eq:couplingeq} does not hold, considering that with Proposition~\ref{lemm:lemmaone}, we have $\beta \geq \delta$.
\end{proof} 
\noindent  \textbf{Remark:} This result gives a first insight on the possible regimes that could allow $m_{x_{0},R}(x_{0}) >1$. First, we can be in the case $\beta=\delta=\alpha-1$, and in that case $C_{\beta} > C_{\delta}$ according to Proposition~\ref{lemm:lemmaone}, so that the birth rate is higher than the death rate for every energy. Otherwise, still according to Proposition~\ref{lemm:lemmaone}, we have to be in the case $\beta > \delta$, so small individuals are more likely to die. However, there is a small probability that they reach high energies, and in that case, they could give birth to a lot of offspring, which could balance the high death rate of these newborns. The allometric coefficient on the birth rate should then be sufficiently high compared to the death rate, which is expressed in Proposition~\ref{prop:beplusde}.

\subsubsection{Proof of 5. in Section~\ref{theo:sharp}}
\label{subsec:deltaalphamoinsun}
Under the allometric setting of Section~\ref{subsec:allomsett} and if $\delta = \alpha-1$, we can be more precise in our computations and extend the result of Proposition~\ref{prop:beplusde}. We begin with an intermediate result in the case $\beta=\delta=\alpha-1$.

\begin{lemme}
Under the allometric setting of Section~\ref{subsec:allomsett} and Assumption~\ref{ass:gainenergy}, if $\beta=\delta=\alpha-1$, we have that $\mathcal{R}_{0} \neq \varnothing$ and
$$ \forall R \in \mathcal{R}_{0}, \forall x >0, \quad m_{0,R}(x) = C_{\beta}/C_{\delta} .$$
\label{prop:egalite}
\end{lemme}
\begin{proof}
Thanks to Lemma~\ref{cestlafindesharicots}, we introduce $R_{0}$ such that $]R_{0}, + \infty [ =  \mathcal{R}_{0}$ and we work with $R > R_{0}$ in the following. Thanks to Proposition~\ref{eq:gammalpha}, we have $\beta=\delta=\alpha-1=\gamma-1$ and $C_{\gamma}>C_{\alpha}$. The result then follows from the equality in~\eqref{eq:majorantm} and a straightforward computation.
\end{proof}

We are now ready to prove Proposition~\ref{corr:pointquatre}, which entails point \textbf{5.} in Section~\ref{theo:sharp}. Notice that under the allometric setting of Section~\ref{subsec:allomsett}, if $\delta=\alpha-1$, then Assumption~\ref{hyp:probamort} is verified by point \textbf{2.} in Section~\ref{theo:sharp}. Also, with Assumptions~\ref{hyp:probamort}, \ref{hyp:tempsinf} and \ref{ass:supercritical}, $d/\tilde{b}$ is non-increasing by Proposition~\ref{lemm:lemmaone}, and $d/g(.,R)$ has a divergent integral near $+ \infty$ by Proposition~\ref{eq:gammalpha}.

\begin{prop}
Under the allometric setting of Section~\ref{subsec:allomsett} with $\delta = \alpha-1$, under Assumptions~\ref{hyp:tempsinf} and \ref{ass:supercritical}, we have 
$$(\beta = \alpha-1 \hspace{0.1 cm} \mathrm{and} \hspace{0.1cm} C_{\beta}>C_{\delta}) \hspace{0.2cm} \mathrm{or} \hspace{0.2cm} \bigg( \beta \geq \alpha -1 + \dfrac{C_{\delta}}{C_{\gamma}-C_{\alpha}} \bigg).$$
\label{corr:pointquatre}
\end{prop}

\begin{proof}
Thanks to Lemma~\ref{lemm:ass}, we work under Assumption~\ref{ass:gainenergy}. Thanks to Lemma~\ref{cestlafindesharicots}, we introduce $R_{0}$ such that $]R_{0}, + \infty [ =  \mathcal{R}_{0}$ and we work with $R > R_{0}$ and $x>0$ in the following. First, thanks to Proposition~\ref{lemm:lemmaone}, we necessarily have $\beta \geq \alpha-1$. If $\beta = \alpha-1$, we obtain with~\eqref{eq:comparem} and Lemma~\ref{prop:egalite} that $m_{x,R}(x) \leq m_{0,R}(x) = C_{\beta}/C_{\delta}$, so $C_{\beta} > C_{\delta}$ is a necessary condition to verify Assumption~\ref{ass:supercritical}.
\\
Now, let us suppose by contradiction that $\alpha -1 < \beta < \alpha -1 + \dfrac{C_{\delta}}{C_{\gamma}-C_{\alpha}}$. We recall that $C_{R} = \phi(R)C_{\gamma} - C_{\alpha} \leq C_{\gamma}-C_{\alpha}$ and we obtain from~\eqref{eq:majorantm}:
\begin{align*}
m_{0,R}(x) & = \dfrac{C_{\beta}C_{\delta}}{C_{R}^{2}} \displaystyle{\int_{x}^{+\infty} \bigg(\dfrac{u^{\beta-\alpha+1}-x^{\beta-\alpha+1}}{\beta-\alpha+1}\bigg)\dfrac{1}{u} \bigg(\dfrac{x}{u}\bigg)^{C_{\delta}/C_{R}} }\mathrm{d}u,
\end{align*}
so according to~\eqref{eq:comparem}, we have
\begin{align*}
m_{x,R}(x) & \leq \dfrac{C_{\beta}C_{\delta}}{C_{R}^{2}(\beta-\alpha+1)}x^{C_{\delta}/C_{R}} \displaystyle{\int_{x}^{+\infty} u^{\beta-\alpha-(C_{\delta}/C_{R})}\mathrm{d}u}.
\end{align*}
From the definition of $C_{R}$, $\beta-\alpha-\dfrac{C_{\delta}}{C_{R}} \leq \beta - \alpha - \dfrac{C_{\delta}}{C_{\gamma}-C_{\alpha}} < -1$, so that the integral above converges and we finally get
\begin{align*}
m_{x,R}(x) & \leq \dfrac{C_{\beta}C_{\delta}}{C_{R}^{2}(\beta-\alpha+1)(\alpha+\dfrac{C_{\delta}}{C_{R}}-\beta-1)}x^{\beta-\alpha+1}.
\end{align*}
We can pick $R_{1}$ high enough so that $C_{R} \geq \dfrac{C_{\gamma}-C_{\alpha}}{2}$ for $R>R_{1}$, because $\phi(R) \underset{R \rightarrow + \infty}{\longrightarrow} 1$, and we obtain for $R>R_{1}$:
\begin{align*}
m_{x,R}(x) & \leq \dfrac{4C_{\beta}C_{\delta}}{(C_{\gamma}-C_{\alpha})^{2}(\beta-\alpha+1)(\alpha+\dfrac{C_{\delta}}{C_{\gamma}-C_{\alpha}}-\beta-1)}x^{\beta-\alpha+1}.
\end{align*}
Finally, we realize that if $x$ is close enough to 0, $m_{x,R}(x) \leq 1$ for any $R>R_{1}$. This contradicts Assumption~\ref{ass:supercritical} by Proposition~\ref{prop:galtonwatson} and concludes.
\end{proof}

\subsection{Proof of 6. and 7. of Section~\ref{theo:sharp}}
\label{subsec:laprevue}
In this section, we work under the allometric setting of Section~\ref{subsec:allomsett}, with Assumptions~\ref{hyp:probamort}, \ref{hyp:tempsinf} and \ref{ass:supercritical}. The result of point \textbf{7.} of Section~\ref{theo:sharp} in the case $\alpha-1 < \beta \leq \alpha$ follows from Section~\ref{subsec:deltaalphamoinsun} (see the proof of Corollary~\ref{corr:thirdone}) and the fact that Assumption~\ref{hyp:probamort} is verified if $\delta=\alpha-1$ (point \textbf{2.} in Section~\ref{theo:sharp}). Then, in this section, except in Corollary~\ref{corr:thirdone}, we work with $\beta > \alpha$ and prove points \textbf{6.} and \textbf{7.} of Section~\ref{theo:sharp}. Recall that for $k \geq 1$, $M^{k}$ is the event $\{ N_{x_{0},\xi_{0}} \geq k \}$ (\textit{i.e.} an individual starting from $\xi_{0}$ with characteristic energy $x_{0}$ and resources $R$ has at least $k$ direct offspring during its life). Under Assumption~\ref{ass:pasinfini}, we can write
\begin{align*}
m_{x_{0},R}(\xi_{0}) = \mathbb{E}(N_{x_{0},R,\xi_{0}}) = \sum\limits_{k \geq 1} \mathbb{P}_{x_{0},R,\xi_{0}}(M^{k}).
\end{align*}
\noindent We will study the general term of this sum to obtain points \textbf{6.} and \textbf{7.} of Section~\ref{theo:sharp}. In the following proof, we assume by contradiction that $\beta > \alpha$ and $\delta < \alpha-1$. We will contradict Assumption~\ref{ass:supercritical}, \textit{i.e.} obtain a subcritical or critical process. We highlight two energy regimes: the low-energy regime and the high-energy regime. The contribution to the expected number of offspring in the low-energy regime goes to 0 when $x_{0}$ goes to 0, thanks to the assumption $\delta < \alpha-1$. For the high-energy regime, it is the assumption $\beta > \alpha$ that is important to obtain the same convergence. A biological interpretation is the following: in the low-energy regime, what prevents individuals to reproduce is that the death rate dominates the energy dynamics (which is expressed by $\delta < \alpha-1$) so that individuals are way too likely to die; in the high-energy regime, what is problematic for reproduction is that the birth rate dominates the energy dynamics ($\beta > \alpha$), so that individuals lose energy too frequently and get back to the low-energy regime.
\\\\
By Lemma~\ref{lemm:ass}, we work under Assumption~\ref{ass:gainenergy}, so Lemma~\ref{cestlafindesharicots} holds true. We fix $R>R_{0}$ once and for all and also work with some $x_{0}>0$, which will eventually converge to 0. Thanks to Assumptions~\ref{hyp:probamort} and \ref{hyp:tempsinf}, by Theorem~\ref{theo:cns}, our process is almost surely biologically relevant. Recall that this means that $\xi$ almost surely never reaches $\flat$ and dies in finite time. 
\\\\
In Section~\ref{subsec:lowenergy}, we will use results for the construction of Lemma~\ref{lemm:coupling}, with $\xi_{0}>0$. In Section~\ref{subsec:highenergy}, we will use results for the construction of Section~\ref{subsec:construc}. Thus, we prove now a technical result that allows us to use both constructions, started from the same $(x_{0},R,\xi_{0})$. In the following, for $k \geq 1$, we write $J_{k}$ for the $k$-th jump time associated to $\xi^{\mathrm{aux}}$ defined in the coupling of Lemma~\ref{lemm:coupling}, and $\hat{J}_{k}$ for the $k$-th jump time associated to $\xi^{\mathrm{aux}}$ defined in Section~\ref{subsec:construc}, $T_{d}$ and $\hat{T}_{d}$ for the associated death time, $M^{k}$ or $\hat{M}^{k}$ for the associated events ``the $k$-th first jumps are birth jumps'', and $N$ or $\hat{N}$ for the total number of birth jumps during the trajectory. Considering the auxiliary process $\xi^{\mathrm{aux}}$ in the coupling of Lemma \ref{lemm:coupling}, we write $S^{k}$, for the maximal value reached by $\xi^{\mathrm{aux}}$ before $J_{k}$. In our setting, $\xi^{\mathrm{aux}}$ almost surely avoids $\flat$, and the energy is increasing outside birth jumps. $S^{k}$ is then clearly the maximal energy reached at a jump time before $J_{k}$ (\textit{i.e.} $S^{k} := \underset{t < J_{k}}{\sup} \hspace{0.1 cm} \xi^{\mathrm{aux}}_{t} = \underset{1 \leq i \leq k}{\max} \hspace{0.1 cm} \xi^{\mathrm{aux}}_{J_{i}-}$). In the same manner, considering the auxiliary process $\xi^{\mathrm{aux}}$ from the construction of Section~\ref{subsec:construc}, we write $\hat{S}^{k}$, for the maximal value reached by $\xi^{\mathrm{aux}}$ before $\hat{J}_{k} $. We have to be more careful here, because $\xi^{\mathrm{aux}}$ can jump to $\partial$ in the construction of Section~\ref{subsec:construc} and some $\hat{J}_{i}$ can be equal to $+ \infty$. Thus, we set $\xi^{\mathrm{aux}}_{+ \infty} = \partial$ and $\partial < x$ for every $x >0$. With these conventions, $\hat{S}^{k}$ is again the maximal energy reached at a jump time before $\hat{J}_{k}$ (\textit{i.e.} $\hat{S}^{k} := \underset{t < \hat{J}_{k}}{\sup} \hspace{0.1 cm} \xi^{\mathrm{aux}}_{t} = \underset{1 \leq i \leq k}{\max} \hspace{0.1 cm} \xi^{\mathrm{aux}}_{\hat{J}_{i}-}$).
\begin{lemme}
Under the general setting of Section~\ref{subsec:gensett}, under Assumptions~\ref{hyp:probamort} and \ref{hyp:tempsinf}, for every $k \geq 1$, $J^{k}$ conditionally to $M^{k}$ and $\hat{J}^{k}$ conditionally to $\hat{M}^{k}$ have the same law. It follows immediately from the definition of $S^{k}$ and $\hat{S}^{k}$ that $S^{k}$ conditionally to $M^{k}$ and $\hat{S}^{k}$ conditionally to $\hat{M}^{k}$ have the same law.
\label{lemme:split}
\end{lemme}
\begin{proof}
We prove the result by induction on $k$. Thanks to Theorem~\ref{theo:cns}, the process is almost surely biologically relevant, so the process avoid 0 and $+ \infty$ and dies in finite time almost surely. In this setting, on the event $\hat{M}^{1}$, $\hat{J}_{1} \neq + \infty$, and the law of $\hat{J}_{1}$ is an inhomogeneous exponential law given in \eqref{eq:jun}, with rate depending on $b_{x_{0}}+d$. On the event $M^{1}$, $J_{1} = J_{1} \wedge T_{d}$, and from the construction of Lemma~\ref{lemm:coupling}, the law of $J_{1}$, respectively $T_{d}$, is an inhomogeneous exponential law with the same form as in \eqref{eq:jun}, but replacing $b_{x_{0}}+d$ with $b_{x_{0}}$, respectively $d$. Furthermore, $J_{1}$ and $T_{d}$ are independent. Hence, the law of $J_{1} \wedge T_{d}$ is an inhomogeneous exponential law whose rate is the sum of those of $J_{1}$ and $T_{d}$. This entails that $J_{1}$ given $M^{1}$ and $\hat{J}_{1}$ given $\hat{M}^{1}$ have the same law. 
\\
Now, suppose that the result is valid for some $k \geq 1$. By Markov property, on the event $M^{k+1}$, respectively $\hat{M}^{k+1}$, we can apply again this reasoning starting from time $J_{k}$, respectively $\hat{J}_{k}$ (which is finite on the event $\hat{M}^{k}$), hence the result. 
\end{proof}
For $j \geq 0$, we write $\heartsuit_{j} := x_{0}^{\frac{3}{\alpha-\beta}}+j$. We choose $x_{0}$ smaller than 1 and small enough to have $\heartsuit_{0} \geq 2x_{0}$ (this is possible because $3/(\alpha-\beta) <0$). Thanks to Corollary~\ref{corr:swagosss}, we work under Assumption~\ref{ass:pasinfini} and we have
\begin{multline}
\hspace{0.6cm} m_{x_{0},R}(x_{0})  = \sum\limits_{k \geq 1} \mathbb{P}_{x_{0},R,x_{0}}(M^{k}) \\
\hspace{1.4cm} = \sum\limits_{k \geq 1} \mathbb{P}_{x_{0},R,x_{0}}(M^{k} \cap \{S^{k} \leq \heartsuit_{0}\} ) + \sum\limits_{k \geq 1} \sum\limits_{j \geq 0} \mathbb{P}_{x_{0},R,x_{0}}(M^{k} \cap \{\heartsuit_{j} < S^{k} \leq \heartsuit_{j+1}\} ) \\
 = \sum\limits_{k \geq 1} \mathbb{P}_{x_{0},R,x_{0}}(M^{k} \cap \{S^{k} \leq \heartsuit_{0}\} ) + \sum\limits_{k \geq 1} \sum\limits_{j \geq 0} \mathbb{P}_{x_{0},R,x_{0}}(\hat{M}^{k} \cap \{\heartsuit_{j} < \hat{S}^{k} \leq \heartsuit_{j+1}\} ).
\label{eq:lowhigh}
\end{multline}
We have partitioned our events $M^{k}$ depending on the maximal energy $S^{k}$ and used Lemma~\ref{lemme:split} for the third line. In the following, we seek for a contradiction with Assumption~\ref{ass:supercritical}. Our proof is divided in two parts, corresponding to the study of the two sums in the right-hand side above. The left-most simple sum accounts for the low-energy regime, where individuals are more likely to die fast than to give birth to many offspring. The double-sum represents the high-energy regime, and we will show that the probability to reach such a level of energy is decreasing to 0 when $x_{0}$ goes to 0, if $\beta > \alpha$. Thus, we want to show that the previously mentioned sums both converge to 0 when $x_{0}$ goes to 0. If this happens, thanks to the previous decomposition, we would have $m_{x_{0},R}(x_{0}) <1$ for $x_{0}$ small enough, which contradicts Assumption~\ref{ass:supercritical} by Proposition~\ref{prop:galtonwatson}, and proves point \textbf{6.} of Section~\ref{theo:sharp}. In the following, we refer to all the previously described setting as `the setting of Section~\ref{subsec:laprevue}'.
\subsubsection{Low-energy regime}
\label{subsec:lowenergy}
We study the sum $$\sum\limits_{k \geq 1} \mathbb{P}_{x_{0},x_{0}}(M^{k} \cap \{S^{k} \leq \heartsuit_{0}\} ).$$
We use the construction of Lemma~\ref{lemm:coupling} and define another coupling between our initial individual process and a continuum of similar processes $(\zeta^{x}_{t})_{t \geq 0}$ for $x >0$.
 Precisely, for $x>0$, we define $b^{x}_{x_{0}} : u \mapsto b_{x_{0}}(u)\mathbb{1}_{u \leq x} + b_{x_{0}}(x)\mathbb{1}_{u > x}$, and similarly $d^{x} : u \mapsto d(u)\mathbb{1}_{u \leq x} + d(x)\mathbb{1}_{u > x}$. These are simply the functions $b_{x_{0}}$ and $d$ frozen after $x$. We define the process $(\zeta^{x}_{t})_{t \geq 0}$ in the exact same manner as $(\xi_{t})_{t \geq 0}$, using the same random variables $(F_{i})_{i \geq 0}$ and starting from $\xi_{0}$ at time 0, but simply replacing $\xi$, $b_{x_{0}}$ and $d$ by respectively $\zeta^{x}$, $b^{x}_{x_{0}}$ and $d^{x}$ in the procedure of Lemma~\ref{lemm:coupling}. This is indeed a coupling between $\xi$ and $\zeta^{x}$, because we use the same random variables to define the jump times. However, we insist on the fact that we possibly define different jump times $(J^{x}_{i})_{i \geq 1}$ for $\zeta^{x}$, even if we still use the same exponential variables as in the definition of the $(J_{i})_{i \geq 1}$. 
\\\\
We naturally write $N^{x}_{x_{0},R,\xi_{0}}$ for the number of births occuring for the process $\zeta^{x}$. In the same manner as in Section~\ref{proof:beplusde}, we also remark that it is possible to formally define $\zeta^{x}$ for $x_{0}=0$ and to compare it to the other processes. The previous coupling allows us to state the following lemma.
\begin{lemme}
Under the setting of Section \ref{subsec:laprevue}, for $R>R_{0}$, $x>0$, $x_{0}>0$, $\xi_{0}>0$, $k \geq 1$, $$\{ N_{x_{0},R,\xi_{0}}\mathbb{1}_{S^{k} \leq x} \geq k \} \subseteq \{ N^{x}_{x_{0},R,\xi_{0}} \geq k \} \subseteq \{ N^{x}_{0,R,\xi_{0}} \geq k \}.$$
\label{lemme:couplage}
\end{lemme}
\begin{proof}
On the event $\{ S^{k} \leq x \} $, since $b_{x_{0}}$ and $d$ coincide respectively with $b_{x_{0}}^{x}$ and $d^{x}$ before $x$, $\xi$ and $\zeta^{x}$ also coincide until time $J_{k}$. Indeed, they both start from $\xi_{0}$ and as soon as the energy does not exceed $x$, the birth and death rates and the energy dynamics are the same. Hence, on the event $\{ S^{k} \leq x \} $, if the $k$ first jumps are births for $\xi$, they also are births for $\zeta^{x}$. This justifies the left-most inclusion of the lemma.
\\\\
Now, remark that in the same manner as in Section~\ref{proof:beplusde}, we can compare $\zeta_{x_{0}}^{x}$ and $\zeta_{0}^{x}$ and obtain the same conclusion for the $N^{x}$ as for $N$ in Lemma~\ref{prop:super}. The only thing that differs between $\xi$ and $\zeta^{x}$ is the shape of the birth and death rates, but they still verify the necessary condition highlighted in Lemma~\ref{prop:super}, that is $d^{x}/b_{x_{0}}^{x}$ is non-increasing, since $\delta < \alpha < \beta$. One can check that the exact same reasoning applies and we obtain $$ N^{x}_{0,R,\xi_{0}} \geq N^{x}_{x_{0},R,\xi_{0}},$$
which concludes for the right-most inclusion.
\end{proof}
\noindent Applying Lemma~\ref{lemme:couplage} leads to
\begin{align*}
\sum\limits_{k \geq 1} \mathbb{P}_{x_{0},R,x_{0}}(M^{k} \cap \{S^{k} \leq \heartsuit_{0}\} ) & = \sum\limits_{k \geq 1} \mathbb{P}(N_{x_{0},R,x_{0}}\mathbb{1}_{S^{k} \leq \heartsuit_{0}} \geq k ) \\
& \leq \sum\limits_{k \geq 1} \mathbb{P}( N^{\heartsuit_{0}}_{0,R,x_{0}} \geq k ) \\
& = \mathbb{E}(N_{0,R,x_{0}}^{\heartsuit_{0}}).
\end{align*}
Remark that we did not used any information about $\heartsuit_{0}$, so in fact we have
\begin{align}
\forall x_{0}>0, \forall y >0, \quad \sum\limits_{k \geq 1} \mathbb{P}_{x_{0},R,x_{0}}(M^{k} \cap \{S^{k} \leq y\} )   \leq \mathbb{E}(N_{0,R,x_{0}}^{y}).
\label{eq:majorationnnn}
\end{align}
\noindent In the same manner as in the proof of Proposition~\ref{prop:beplusde}, we can compute this expectation because there is no energy loss anymore so it is easier to understand the trajectory of $\zeta_{.,0,R,x_{0}}^{\heartsuit_{0}}$. Since $x_{0} < \heartsuit_{0}$, there exists $t_{\heartsuit} >0$ such that $A_{x_{0}}(t_{\heartsuit}) = \heartsuit_{0}$. We obtain from Equation~\eqref{eq:enfantssansnmort} applied to the rates $b_{x_{0}}^{\heartsuit_{0}}$ and $d^{\heartsuit_{0}}$:
\begin{align*}
\mathbb{E}(N_{0,R,x_{0}}^{\heartsuit_{0}}) & = \displaystyle{\int_{0}^{+\infty}\bigg(\int_{0}^{t}b_{x_{0}}^{\heartsuit_{0}}(A_{x_{0}}(w))\mathrm{d}w\bigg)d^{\heartsuit_{0}}(A_{x_{0}}(t))e^{-\int_{0}^{t}d^{\heartsuit_{0}}(A_{x_{0}}(\tau))\mathrm{d}\tau}\mathrm{d}t} \\
& = \displaystyle{\int_{0}^{t_{\heartsuit}}\bigg(\int_{0}^{t}b_{x_{0}}(A_{x_{0}}(w))\mathrm{d}w\bigg)d(A_{x_{0}}(t))e^{-\int_{0}^{t}d(A_{x_{0}}(\tau))\mathrm{d}\tau}\mathrm{d}t} \\
& + \displaystyle{\int_{t_{\heartsuit}}^{+\infty}\bigg( \bigg(\int_{0}^{t_{\heartsuit}}b_{x_{0}}(A_{x_{0}}(w))\mathrm{d}w\bigg)  + (t-t_{\heartsuit})b_{x_{0}}(\heartsuit_{0}) \bigg)d(\heartsuit_{0})e^{-\int_{0}^{t_{\heartsuit}}d(A_{x_{0}}(\tau))\mathrm{d}\tau}  e^{-(t-t_{\heartsuit})d(\heartsuit_{0})} \mathrm{d}t.} 
\end{align*}
Now, we use the fact that $\int_{t_{\heartsuit}}^{+ \infty} (t-t_{\heartsuit})e^{-(t-t_{\heartsuit})d(\heartsuit_{0})} \mathrm{d}t = d(\heartsuit_{0})^{-2}$. We also use the definition of $d^{\heartsuit_{0}}$ and $b_{x_{0}}^{\heartsuit_{0}}$, and the change of variables $u=A_{x_{0}}(t)$ to obtain
\begin{multline}
\hspace{1.5cm }\mathbb{E}(N_{0,R,x_{0}}^{\heartsuit_{0}}) = \displaystyle{\int_{x_{0}}^{\heartsuit_{0}}\bigg(\int_{x_{0}}^{u}\dfrac{\tilde{b}}{g}(w)\mathrm{d}w\bigg)\dfrac{d}{g}(u) e^{-\displaystyle{\int_{x_{0}}^{u}}\dfrac{d}{g}(\tau)\mathrm{d}\tau}\mathrm{d}u} \\
+\bigg( \bigg( \int_{x_{0}}^{\heartsuit_{0}}\dfrac{\tilde{b}}{g}(w)\mathrm{d}w \bigg) + \dfrac{b_{x_{0}}}{d}(\heartsuit_{0}) \bigg)\exp\bigg(-\displaystyle{\int_{x_{0}}^{\heartsuit_{0}}\frac{d}{g}(\tau)\mathrm{d}\tau}\bigg),
\label{eq:ncoeur}
\end{multline}
\noindent where we write $g$ instead of $g(.,R)$ for the sake of simplicity. Also, from Proposition~\ref{eq:gammalpha}, as soon as $\beta > \alpha$, $\tilde{b}/g$ is increasing. Hence, we compute the following upper bound:

\begin{align*}
\mathbb{E}(N_{0,R,x_{0}}^{\heartsuit_{0}}) & \leq x_{0} \dfrac{\tilde{b}}{g}(2x_{0}) \displaystyle{\int_{x_{0}}^{2x_{0}}\dfrac{d}{g}(u) e^{-\displaystyle{\int_{x_{0}}^{u}}\dfrac{d}{g}(\tau)\mathrm{d}\tau}\mathrm{d}u} + \heartsuit_{0} \dfrac{\tilde{b}}{g}(\heartsuit_{0}) \displaystyle{\int_{2x_{0}}^{\heartsuit_{0}}\dfrac{d}{g}(u) e^{-\displaystyle{\int_{x_{0}}^{u}}\dfrac{d}{g}(\tau)\mathrm{d}\tau}\mathrm{d}u} \\
& + \bigg( \heartsuit_{0} \dfrac{\tilde{b}}{g}(\heartsuit_{0})  + \dfrac{\tilde{b}}{d}(\heartsuit_{0}) \bigg)e^{-\displaystyle{\int_{x_{0}}^{\heartsuit_{0}}}\dfrac{d}{g}(\tau)\mathrm{d}\tau} \\
& \leq x_{0} \dfrac{\tilde{b}}{g}(2x_{0}) + \heartsuit_{0} \dfrac{\tilde{b}}{g}(\heartsuit_{0}) e^{-\displaystyle{\int_{x_{0}}^{2x_{0}}}\dfrac{d}{g}(\tau)\mathrm{d}\tau} + \bigg(\heartsuit_{0} \dfrac{\tilde{b}}{g}(\heartsuit_{0})  + \dfrac{\tilde{b}}{d}(\heartsuit_{0}) \bigg)e^{-\displaystyle{\int_{x_{0}}^{\heartsuit_{0}}}\dfrac{d}{g}(\tau)\mathrm{d}\tau}.
\end{align*}
At this point, one can check that all the terms in the upper bound converge to 0 when $x_{0}$ goes to 0, under the allometric setting of Section~\ref{subsec:allomsett} and with $\beta > \alpha > \delta +1$. Indeed, the exponential factors with general term $d/g$ dominate all the other quantities near 0 (they both are of order $\exp(-cx_{0}^{\delta-\alpha+1})$ with some constant $c$, the other terms are powers of $x_{0}$), and $\frac{\tilde{b}}{g}(2x_{0})$ converges to 0 when $x_{0}$ goes to 0.

\subsubsection{High-energy regime}
\label{subsec:highenergy}
Now, we deal with the double sum:
$$\sum\limits_{k \geq 1} \sum\limits_{j \geq 0} \mathbb{P}_{x_{0},R,x_{0}}(\hat{M}^{k} \cap \{\heartsuit_{j} < \hat{S}^{k} \leq \heartsuit_{j+1}\} ).$$
We begin with another coupling argument. Recall that in this section, we use the construction of Section~\ref{subsec:construc}. For $x>0$, we define 
$$ \tilde{N}_{x} := \sup \{i \geq 1, \forall 1 \leq j \leq i, \quad U_{j} \leq r(x) \}, $$
with the convention $\sup(\varnothing) := 0$ and $r := \frac{b_{x_{0}}}{b_{x_{0}}+d}$. Remark that $\tilde{N}_{x}+1$ follows simply a geometric law with parameter $1-r(x)$, but it is coupled to $\hat{N}_{x_{0},R,\xi_{0}}$ because it depends on the same $(U_{i})_{i \geq 1}$. This coupling via the only variables $(U_{i})_{i \geq 1}$ allows us to assess two things.
\begin{lemme}
Under the setting of Section \ref{subsec:laprevue}, we have:
\begin{enumerate}
\item For $x, x_{0},\xi_{0}>0$, $k \geq 1$, $\mathbb{P}_{x_{0},R,\xi_{0}}(\hat{N}_{x_{0},R,\xi_{0}}\mathbb{1}_{\hat{S}^{k} \leq x} \geq k ) \leq \mathbb{P}_{x_{0},R,\xi_{0}}(\tilde{N}_{x} \geq k)$.
\item $\tilde{N}_{x}$ is independent from $\hat{S}^{k}$.
\end{enumerate}
\label{lemme:couplagedeux}
\end{lemme}
\begin{proof}
\begin{enumerate}
\item On the event $\{ \hat{S}^{k} \leq x \} $, since $d/b_{x_{0}}$ is non-increasing, $r$ is non-decreasing, so we can assess that for every $t \leq \hat{J}_{k} \wedge \hat{T}_{d}$, $1 -r(\xi^{\mathrm{aux}}_{t}) \geq 1-r(x)$. But on this event, if $\hat{N}_{x_{0},R,\xi_{0}} \geq k$ (\textit{i.e.} on the event $\hat{M}^{k}$), then $\hat{J}_{k} < \hat{T}_{d}$, and for every $i \leq k$, $U_{i} > 1 -r(\xi^{\mathrm{aux}}_{\hat{J}_{i}-}) \geq 1-r(x)$. Hence, 
\begin{align*}
\mathbb{P}_{x_{0},R,\xi_{0}}(\hat{N}_{x_{0},R,\xi_{0}}\mathbb{1}_{\hat{S}^{k} \leq x} \geq k) & \leq \mathbb{P}_{x_{0},R,\xi_{0}}(\forall 1 \leq i \leq k, U_{i} > 1 -r(x)) \\
& = \mathbb{P}_{x_{0},R,\xi_{0}}(\tilde{N}_{x} \geq k).
\end{align*}
\item $\tilde{N}_{x}$ is independent from $(E_{i})_{i \geq 1}$, so independent from the $(\hat{J}_{i})_{i \geq 1}$, hence independent from $\hat{S}^{k}$.
\end{enumerate}
\end{proof}
\noindent This last independence property is crucial for the third line of the following computation. We take $k \geq 1$ and $j \geq 0$, and by Lemma~\ref{lemme:couplagedeux}:
\begin{align*}
\mathbb{P}_{x_{0},R,x_{0}}(\hat{M}^{k} \cap \{\heartsuit_{j} < \hat{S}^{k} \leq \heartsuit_{j+1}\} ) & = \mathbb{P}_{x_{0},R,x_{0}}(\{\hat{N}_{x_{0},R,x_{0}}\mathbb{1}_{\hat{S}^{k} \leq \heartsuit_{j+1}} \geq k  \} \cap \{\heartsuit_{j} < \hat{S}^{k} \} ) \\
& \leq \mathbb{P}_{x_{0},R,x_{0}}(\{\tilde{N}_{\heartsuit_{j+1}} \geq k  \} \cap \{\heartsuit_{j} < \hat{S}^{k}\} ) \\
& = \mathbb{P}_{x_{0},R,x_{0}}(\tilde{N}_{\heartsuit_{j+1}} \geq k ) \mathbb{P}_{x_{0},R,x_{0}}(\heartsuit_{j} < \hat{S}^{k} ). \\
\end{align*}
Finally, interverting the two sums leads to
\begin{align}
\sum\limits_{k \geq 1} \sum\limits_{j \geq 0} \mathbb{P}_{x_{0},R,x_{0}}(\hat{M}^{k} \cap \{\heartsuit_{j} < \hat{S}^{k} \leq \heartsuit_{j+1}\} )& \leq \sum\limits_{j \geq 0} \sum\limits_{k \geq 1} \mathbb{P}_{x_{0},R,x_{0}}(\tilde{N}_{\heartsuit_{j+1}} \geq k ) \mathbb{P}_{x_{0},R,x_{0}}(\heartsuit_{j} < \hat{S}^{k} ).
\label{doubvlesuomme}
\end{align}
Let $j \geq 0$, $k \geq 1$, we have
\begin{align*}
\mathbb{P}_{x_{0},R,x_{0}}(\heartsuit_{j} < \hat{S}^{k} ) & = \mathbb{P}_{x_{0},R,x_{0}}( \underset{1 \leq i \leq k}{\max} \hspace{0.1 cm} \xi^{\mathrm{aux}}_{\hat{J}_{i}-}> \heartsuit_{j} ) \\
& = \sum\limits_{i=1}^{k} \mathbb{P}_{x_{0},R,x_{0}}(\{ \xi^{\mathrm{aux}}_{\hat{J}_{i}-}> \heartsuit_{j} \} \cap \{ \forall  n \in \llbracket 1,i-1 \rrbracket, \xi^{\mathrm{aux}}_{\hat{J}_{n}-} \leq \heartsuit_{j} \}).
\end{align*}
For $0 < x \leq y$, $i \geq 1$ we define
\[ q_{i}(x,y) := \mathbb{P}_{x_{0},R,x}(\{ \xi^{\mathrm{aux}}_{\hat{J}_{i}-}> y \} \cap \{ \forall n \in \llbracket 1,i-1 \rrbracket, \xi^{\mathrm{aux}}_{\hat{J}_{n}-} \leq y \}) \]
so that 
\[ \mathbb{P}_{x_{0},R,x_{0}}(\heartsuit_{j} < \hat{S}^{k} ) = \sum\limits_{i=1}^{k} q_{i}(x_{0},\heartsuit_{j}). \]
Remark that with our convention $\partial <x$ for all $x>0$, on the event $\{ \xi^{\mathrm{aux}}_{\hat{J}_{i}-}> y \}$, the $i$ first jumps are necessarily birth events, for any $i \geq 1$.

\begin{lemme}
Under the setting of Section~\ref{subsec:laprevue}, we have 
\[\forall i \geq 1, \forall y \geq x_{0},\forall x \in (0, y-x_{0}], \quad q_{i}(x,y) \leq \exp\left(-\displaystyle{\int_{y-x_{0}}^{y}}\dfrac{b_{x_{0}}+d}{g}(u) \mathrm{d}u \right). \]
\end{lemme}

\begin{proof}
We prove this lemma by induction on $i$. First, we know the law of $\hat{J}_{1}$, so for $y >x_{0}$ and $x \in (0, y-x_{0}]$: 
\begin{align*}
q_{1}(x,y) & = \mathbb{P}_{x_{0},R,x}(\{ \xi^{\mathrm{aux}}_{\hat{J}_{1}-}> y \}) \\
& = \displaystyle{\int_{y}^{+ \infty} \dfrac{b_{x_{0}}+d}{g}(u) \exp\left(-\displaystyle{\int_{x}^{u}}\dfrac{b_{x_{0}}+d}{g}(\tau) \mathrm{d}\tau \right) \mathrm{d}u } \\
& = \exp\left(-\displaystyle{\int_{x}^{y}}\dfrac{b_{x_{0}}+d}{g}(u) \mathrm{d}u \right) \\
& \leq \exp\left(-\displaystyle{\int_{y-x_{0}}^{y}}\dfrac{b_{x_{0}}+d}{g}(u) \mathrm{d}u \right),
\end{align*}
because $\int_{.}^{+ \infty }(b_{x_{0}}/d)/g = + \infty$ as soon as $\beta > \alpha$ under the setting of Section~\ref{subsec:laprevue}. Then, for some $i \geq 1$, we suppose that the property is true for $q_{i}$, and we use Markov property to assess
\begin{align*}
q_{i+1}(x,y) & = \mathbb{P}_{x_{0},R,x}(\{\xi^{\mathrm{aux}}_{\hat{J}_{1}-} \leq y \} \cap  \{ \xi^{\mathrm{aux}}_{\hat{J}_{i+1}-}> y \} \cap \{ \forall n \in \llbracket 2,i \rrbracket, \xi^{\mathrm{aux}}_{\hat{J}_{n}-} \leq y \})\\
& = \displaystyle{\int_{x}^{y}}\dfrac{b_{x_{0}}}{g}(u) \exp\left(-\displaystyle{\int_{x}^{u}}\dfrac{b_{x_{0}}+d}{g}(\tau) \mathrm{d}\tau \right) q_{i}(u-x_{0},y) \mathrm{d}u \\
& \leq \exp\left(-\displaystyle{\int_{y-x_{0}}^{y}}\dfrac{b_{x_{0}}+d}{g}(u) \mathrm{d}u \right) \displaystyle{\int_{x}^{y}}\dfrac{b_{x_{0}}+d}{g}(u) \exp\left(-\displaystyle{\int_{x}^{u}}\dfrac{b_{x_{0}}+d}{g}(\tau) \mathrm{d}\tau \right) \mathrm{d}u \\
& \leq \exp\left(-\displaystyle{\int_{y-x_{0}}^{y}}\dfrac{b_{x_{0}}+d}{g}(u) \mathrm{d}u \right),
\end{align*}
where the third line comes from the induction hypothesis.
\end{proof}
\noindent This finally gives, because $2x_{0} \leq \heartsuit_{j}-x_{0}$,
\begin{align*}
\mathbb{P}_{x_{0},R,x_{0}}(\heartsuit_{j} < \hat{S}^{k} ) & \leq k \exp\left(-\displaystyle{\int_{\heartsuit_{j}-x_{0}}^{\heartsuit_{j}}}\dfrac{b_{x_{0}}+d}{g}(u) \mathrm{d}u \right) .
\end{align*}
Hence, as $\tilde{N}_{\heartsuit_{j+1}}+1$ follows a geometric law with parameter $1-r(\heartsuit_{j+1})$, we obtain from~\eqref{doubvlesuomme} that
\begin{align*}
\sum\limits_{k \geq 1} \sum\limits_{j \geq 0} \mathbb{P}_{x_{0},R,x_{0}}(\hat{M}^{k} \cap \{\heartsuit_{j} < \hat{S}^{k} \leq \heartsuit_{j+1}\} ) & \leq  \sum\limits_{j \geq 0} \left( \sum\limits_{k \geq 1} k r(\heartsuit_{j+1})^{k} \right) \exp\left(-\displaystyle{\int_{\heartsuit_{j}-x_{0}}^{\heartsuit_{j}}}\dfrac{b_{x_{0}}+d}{g}(u) \mathrm{d}u \right) \\
& \leq \sum\limits_{j \geq 0} \dfrac{r(\heartsuit_{j+1})}{(1-r(\heartsuit_{j+1}))^{2}} \exp\left(-\displaystyle{\int_{\heartsuit_{j}-x_{0}}^{\heartsuit_{j}}}\dfrac{b_{x_{0}}}{g}(u) \mathrm{d}u \right) \\
& = \sum\limits_{j \geq 0} \dfrac{b_{x_{0}}^{2}+b_{x_{0}}d}{d^{2}}(\heartsuit_{j+1}) \exp\left(-\displaystyle{\int_{\heartsuit_{j}-x_{0}}^{\heartsuit_{j}}}\dfrac{b_{x_{0}}}{g}(u) \mathrm{d}u \right).
\end{align*}
\noindent By assumption, $\beta > \alpha$, so that $b_{x_{0}}/g$ is increasing, which leads to the first inequality in the following. Also, the reader can check that we can pick $x_{0}$ small enough so that for all $j \geq 0$, $(\heartsuit_{j}-x_{0})^{\beta-\alpha} \geq \heartsuit_{j}^{\beta-\alpha}/2 $, which entails the second line in the following.
\begin{align*}
\sum\limits_{k \geq 1} \sum\limits_{j \geq 0} \mathbb{P}_{x_{0},R,x_{0}}(\hat{M}^{k} \cap \{\heartsuit_{j} < \hat{S}^{k} \leq \heartsuit_{j+1}\} ) & \leq \sum\limits_{j \geq 0} \dfrac{b_{x_{0}}^{2}+b_{x_{0}}d}{d^{2}}(\heartsuit_{j+1}) \exp\left(-\dfrac{C_{\beta}}{C_{R}}x_{0} (\heartsuit_{j}-x_{0})^{\beta-\alpha} \right) \\
& \leq \sum\limits_{j \geq 0} \dfrac{b_{x_{0}}^{2}+b_{x_{0}}d}{d^{2}}(\heartsuit_{j+1}) \exp\left(-\dfrac{C_{\beta}}{2C_{R}}x_{0} \heartsuit_{j}^{\beta-\alpha} \right).
\end{align*}
\noindent The term for $j=0$, which we denote in the following as $\varsigma_{0}$, converges to 0 when $x_{0}$ goes to 0, because it is a power of $x_{0}$ multiplied by an exponential term of order $\exp(-x_{0}^{-2})$. Now for $j \geq 1$, recall that $\heartsuit_{j} := x_{0}^{\frac{3}{\alpha-\beta}}+j$. Hence, because $\beta > \alpha$,
\begin{equation}
\begin{split}
x_{0} \heartsuit_{j}^{\beta-\alpha} & = x_{0}(x_{0}^{\frac{3}{\alpha-\beta}}+j)^{\beta-\alpha} \\
& = x_{0}(x_{0}^{\frac{3}{\alpha-\beta}}+j)^{(\beta-\alpha)/3}(x_{0}^{\frac{3}{\alpha-\beta}}+j)^{(\beta-\alpha)/3}(x_{0}^{\frac{3}{\alpha-\beta}}+j)^{(\beta-\alpha)/3} \\
& \geq x_{0}\dfrac{1}{x_{0}}\dfrac{1}{x_{0}}j^{(\beta-\alpha)/3}.
\end{split}
\label{eq:splitting}
\end{equation}
Thus, we have
\begin{align*}
\sum\limits_{k \geq 1} \sum\limits_{j \geq 0} \mathbb{P}_{x_{0},R,x_{0}}(\hat{M}^{k} \cap \{\heartsuit_{j} < \hat{S}^{k} \leq \heartsuit_{j+1}\} ) & \leq \varsigma_{0} + \sum\limits_{j \geq 1} \dfrac{b_{x_{0}}^{2}+b_{x_{0}}d}{d^{2}}(\heartsuit_{j+1}) \exp\left(-\dfrac{C_{\beta}}{2C_{R}}\dfrac{1}{x_{0}} j^{\frac{\beta-\alpha}{3}} \right).
\end{align*}
The general term of the sum in the right-hand side simply converges to 0 when $x_{0}$ goes to 0. Also, for $x_{0} \leq 1$ and $j \geq 1$, we have
$$\exp\left(-\frac{C_{\beta}}{2C_{R}}\frac{1}{x_{0}} j^{\frac{\beta-\alpha}{3}} \right) = e^{-\frac{C_{\beta}}{4C_{R}}\frac{1}{x_{0}}j^{\frac{\beta-\alpha}{3}}}e^{-\frac{C_{\beta}}{4C_{R}}\frac{1}{x_{0}}j^{\frac{\beta-\alpha}{3}}} \leq e^{-\frac{C_{\beta}}{C_{R}}\frac{1}{4x_{0}}} e^{-\frac{C_{\beta}}{4C_{R}}j^{\frac{\beta-\alpha}{3}}}$$
We let the reader check that, uniformly on $x_{0} \leq 1$, the term $\dfrac{b_{x_{0}}^{2}+b_{x_{0}}d}{d^{2}}(\heartsuit_{j+1}) e^{-\frac{C_{\beta}}{C_{R}}\frac{1}{4x_{0}}}$ is dominated by $P(j)$ with some polynomial function $P$ that does not depend on $x_{0}$ (this comes essentially from the exponential term, that allows us to get rid of the dependence on $x_{0}$ by rough bounds).\\
This entails that our general term is dominated by $P(j)e^{-\frac{C_{\beta}}{4C_{R}}j^{\frac{\beta-\alpha}{3}}}$, which is the general term of a converging sum as $\beta > \alpha$. Hence, a dominated convergence argument allow us to conclude: the double sum in the left-hand side above goes to 0 when $x_{0}$ goes to 0.
\subsubsection{Proof of 7. in Section~\ref{theo:sharp}}
\label{sec:coeff}
In this section, we precise what happens for the previous decomposition if $\delta=\alpha-1$. Precisely, we prove the following result.

\begin{prop}
Under the allometric setting of Section~\ref{subsec:allomsett}, Assumptions~\ref{hyp:tempsinf} and~\ref{ass:supercritical}, we have
$$(\beta  > \alpha \hspace{0.1 cm} \mathrm{and} \hspace{0.1 cm} \delta = \alpha -1) \Rightarrow (C_{\delta} \leq C_{\gamma}-C_{\alpha}).$$
\label{prop:decoupagedelta}
\end{prop}

\begin{proof}
We work with $R>R_{0}$ and $x_{0}>0$. Let us suppose by contradiction that $\beta  > \alpha$, $\delta = \alpha -1$ and $C_{\delta} > C_{\gamma}-C_{\alpha}$.
We can still use the decomposition of~\eqref{eq:lowhigh}. We begin by the low-energy regime. The upper bound of Equation~\eqref{eq:majorationnnn} is still valid for $y=\heartsuit_{0}$, so our goal is to prove that $\mathbb{E}(N_{0,R,x_{0}}^{\heartsuit_{0}})$ goes to 0 when $x_{0}$ goes to 0. Equation~\eqref{eq:ncoeur} is still valid here, but we can be more subtle in the computations because $\delta=\alpha-1$. Also, we use Proposition~\ref{corr:pointquatre} (Assumption~\ref{hyp:probamort} is verified if $\delta=\alpha-1$) to consider only the case $\beta \geq \alpha-1 + C_{\delta}/(C_{\gamma}-C_{\alpha})$. Recall also that $C_{R}=\phi(R)(C_{\gamma}-C_{\alpha}) \leq C_{\gamma}-C_{\alpha}$. First if $\beta > \alpha-1 + C_{\delta}/(C_{\gamma}-C_{\alpha})$, we get
\begin{align*}
\mathbb{E}(N_{0,R,x_{0}}^{\heartsuit_{0}}) & = \displaystyle{\int_{x_{0}}^{\heartsuit_{0}}\bigg(\int_{x_{0}}^{u}\dfrac{C_{\beta}}{C_{R}}w^{\beta-\alpha}\mathrm{d}w\bigg)\dfrac{C_{\delta}}{C_{R}}\dfrac{1}{u} \left( \dfrac{x_{0}}{u} \right)^{C_{\delta}/C_{R}} \mathrm{d}u} \\
& + \bigg( \bigg( \int_{x_{0}}^{\heartsuit_{0}}\dfrac{C_{\beta}}{C_{R}}w^{\beta-\alpha}\mathrm{d}w \bigg) + \dfrac{C_{\beta}}{C_{\delta}}\heartsuit_{0}^{\beta-\alpha+1} \bigg)\left( \dfrac{x_{0}}{\heartsuit_{0}} \right)^{C_{\delta}/C_{R}} \\
& \leq \mathfrak{C} x_{0}^{\frac{C_{\delta}}{C_{R}}} \heartsuit_{0}^{\beta-\alpha+1-\frac{C_{\delta}}{C_{R}}},
\end{align*}
with some constant $\mathfrak{C}$ that does not depend on $x_{0}$. Finally, recall that $\heartsuit_{0} := x_{0}^{\frac{3}{\alpha-\beta}}$. The reader can check that all of our previous reasonings still hold true if $\heartsuit_{0}$ is of the form $x_{0}^{\frac{\omega}{\alpha-\beta}}$ with $\omega >1$, and we take such a $\heartsuit_{0}$ until the end of this proof. The previous upper bound then leads to
$$ \mathbb{E}(N_{0,R,x_{0}}^{\heartsuit_{0}}) \leq \mathfrak{C}x_{0}^{\frac{C_{\delta}}{C_{R}}\left( 1 + \frac{\omega}{\beta-\alpha} \right) - \omega\left( 1 + \frac{1}{\beta-\alpha} \right)}. $$
We supposed that $C_{\delta} > C_{\gamma}-C_{\alpha} \geq C_{R}$. Hence, we can choose $\omega$ close enough to 1 to obtain that for all $R>R_{0}$, $\mathbb{E}(N_{0,R,x_{0}}^{\heartsuit_{0}})$ converges to 0 when $x_{0}$ goes to 0. If now $\beta =\alpha-1 + C_{\delta}/(C_{\gamma}-C_{\alpha})$, and $C_{R} < C_{\gamma}-C_{\alpha}$, the same reasoning applies. Otherwise, if $C_{R}=C_{\gamma}-C_{\alpha}$, we obtain an upper bound of the form $\mathfrak{C} x_{0}^{C_{\delta}/C_{R}} \ln(x_{0})$, which still converges to 0 when $x_{0}$ goes to 0.
\\\\
The reasoning in the high-energy regime is exactly the same, since the crucial hypothesis in this case is $\beta > \alpha$ to use dominated convergence. However, we have to check that our new definition for $\heartsuit_{0}$ with $\omega >1$ is still consistent with our reasoning. Precisely, we have to replace the lower bound obtained in~\eqref{eq:splitting} by
\begin{equation*}
\begin{split}
x_{0} \heartsuit_{j}^{\beta-\alpha} & = x_{0}(x_{0}^{\frac{\omega}{\alpha-\beta}}+j)^{\beta-\alpha} \\
& = x_{0}(x_{0}^{\frac{\omega}{\alpha-\beta}}+j)^{(\beta-\alpha)/\omega}(x_{0}^{\frac{\omega}{\alpha-\beta}}+j)^{\frac{\beta-\alpha}{2}(1-\frac{1}{\omega})}(x_{0}^{\frac{\omega}{\alpha-\beta}}+j)^{\frac{\beta-\alpha}{2}(1-\frac{1}{\omega})} \\
& \geq x_{0}\dfrac{1}{x_{0}}x_{0}^{\frac{1-\omega}{2}}j^{\frac{\beta-\alpha}{2}(1-\frac{1}{\omega})},
\end{split}
\end{equation*}
and the reader can check that all our other arguments hold true with these new exponents (what is important is that $1- \omega <0$ and $1-\frac{1}{\omega} >0$).
\end{proof}

\noindent Finally, we highlight the following corollary, which completes the proof of point \textbf{7.} in Section~\ref{theo:sharp}.

\begin{corr}
Under the allometric setting of Section~\ref{subsec:allomsett}, Assumptions~\ref{hyp:probamort},~\ref{hyp:tempsinf} and~\ref{ass:supercritical}, we have that if $\delta = \alpha-1$,
$$(\beta > \alpha-1) \Rightarrow (C_{\delta} \leq C_{\gamma}-C_{\alpha}).$$
\label{corr:thirdone}
\end{corr}

\begin{proof}
First, if $\beta \leq \alpha$ and $C_{\delta} > C_{\gamma}-C_{\alpha}$, then Proposition~\ref{corr:pointquatre} implies $\beta = \alpha-1$, hence $\beta \in ]\alpha-1,\alpha]$ implies $C_{\delta} \leq C_{\gamma}-C_{\alpha}$. Then if $\beta > \alpha$, we conclude with Proposition~\ref{prop:decoupagedelta}. 
\end{proof}

\subsection{Proof of 8. in Section~\ref{theo:sharp}}
\label{sec:pta}
First, in Section~\ref{subsec:construpta}, we work under the general setting of Section~\ref{subsec:gensett}, under Assumptions~\ref{hyp:probamort} and \ref{hyp:tempsinf}, and we suppose in addition that $t \mapsto A_{\xi}(t)$ is well-defined on $\mathbb{R}^{+}$ for every $\xi>0$ to obtain martingale properties for the individual process. Then, in Sections~\ref{subsec:limit}, \ref{subsec:core} and \ref{subsec:coredeux}, we work under the allometric setting of Section~\ref{subsec:allomsett} with $\alpha \leq 1$, Assumptions~\ref{hyp:tempsinf} and \ref{ass:gainenergy} and $R \in \mathcal{R}_{0}$. For these last sections, Assumption~\ref{ass:gainenergy} ensures that Lemma~\ref{cestlafindesharicots} holds true, so $\mathcal{R}_{0} \neq \varnothing$. Remark that we do not need to suppose Assumption~\ref{hyp:probamort} anymore, because if $R \in \mathcal{R}_{0}$, then $T_{0}=+ \infty$. Getting back to the notations of Section~\ref{subsec:construc}, we obtain that, if it is not $\partial$ or $\flat$, $\xi^{\mathrm{aux}}_{t} \leq A_{\xi_{0}}(t)$ for $t \geq 0$ (this is due to jumps for reproduction). We will show that when $\beta \leq \alpha$, there is a positive probability for $\xi^{\mathrm{aux}}_{t}$ to be asymptotically lower bounded by $(1-\varepsilon)A_{\xi_{0}}(t)$, for some $0< \varepsilon <1$. If in addition $\delta < \alpha-1$, the death rate decreases too fast and there is a chance that the full individual process
$\xi_{t} \notin \{ \partial, \flat \}$ for every $t \geq 0$, so that an individual never dies. This would contradict Assumption~\ref{hyp:tempsinf} thanks to Proposition~\ref{prop:cns} and our choice of $R$.
\\\\
First, we will give a third way to construct our process $\xi$ thanks to a Poisson point process. Then, we will extract from this new definition some useful martingale properties. Finally, we were inspired by the concept of asymptotic pseudotrajectory developed by Benaïm and Hirsch \cite{benaim_06} to compare our trajectories to the integral curves associated to a deterministic flow.
\subsubsection{Another construction of the individual process $(\xi_{t})_{t \geq 0}$}
\label{subsec:construpta}

We give here another way to construct $(\xi_{t})_{t \geq 0}$, similar to the coupling of Lemma~\ref{lemm:coupling}, but using a Poisson point process for the construction of the jump events. For this construction, we fix $\xi_{0}>0$, $x_{0}>0$ and $R>0$. We also suppose Assumptions~\ref{hyp:probamort} and \ref{hyp:tempsinf}, so almost surely $T_{d}<+ \infty$ by Theorem~\ref{theo:cns}. In other terms, we work in this section with almost surely biologically relevant processes, so we ignore the value $\flat$ in the following construction, and the processes are well-defined for $t \geq 0$.
\\\\
First, we define an auxiliary process $(X_{t})_{t \geq 0}$ following all the wanted mechanisms except for deaths (so $(X_{t})_{t \geq 0}$ takes its values in $\mathbb{R}^{*}_{+}$). Then, we kill this process to obtain $\xi$.
\begin{enumerate}
\item Let $X_{0} =\xi_{0}$ be the initial energy.
\item Let $\mathcal{N}(\mathrm{d}s,\mathrm{d}h)$ be a Poisson point measure on $\mathbb{R}^{+} \times \mathbb{R}^{+}$, with intensity $\mathrm{d}s \times \mathrm{d}h$. We write $\mathcal{N}_{C}(\mathrm{d}s,\mathrm{d}h) := \mathcal{N}(\mathrm{d}s,\mathrm{d}h) - \mathrm{d}s \mathrm{d}h$ for the associated compensated Poisson point measure. We define the process $(X_{t})_{t \geq 0}$ by
\begin{align} 
X_{t} := &  A_{\xi_{0}}(t) + \displaystyle{\int_{0}^{t}\int_{\mathbb{R}^{+}}}\mathbb{1}_{\{ h \leq b_{x_{0}}(X_{s-})\}} \bigg(A_{X_{s-}-x_{0}}(t-s) - A_{X_{s-}}(t-s)\bigg) \mathcal{N}(\mathrm{d}s, \mathrm{d}h).
\label{eq:constructionofx}
\end{align}
\item We define $E$ a random variable following an exponential law with parameter 1, independent from $\mathcal{N}$. We then stop the process $X$ at the random time \begin{align}
T_{d} := \inf \left\{ t \geq 0, \hspace{0.1 cm} \int_{0}^{t}d(X_{s})\mathrm{d}s = E \right\},
\label{eq:definitiontd}
\end{align}  
and define $\xi_{t} := X_{t}\mathbb{1}_{t \leq T_{d}} + \partial \mathbb{1}_{t > T_{d}} $. 
\end{enumerate}
This construction using a Poisson point process is standard \cite{four_04} \cite{champagnat2005individualbased}. Let $(J_{n})_{n \geq 1}$ be the jump times for the process $X$, and $\mathcal{F}_{t}$ the canonical filtration associated to $\mathcal{N}$. Classical properties of Poisson point measures show that for $n \geq 1$ and given $\mathcal{F}_{J_{n}}$, the next jump time $J_{n+1}$ is distributed as in the construction of Lemma~\ref{lemm:coupling}. This shows that $X$ and $\xi^{\mathrm{aux}}$ of Lemma~\ref{lemm:coupling} have same distribution of sample paths. Hence, this new definition of $\xi$ is consistent with the constructions of Section~\ref{subsec:construc} and Lemma~\ref{lemm:coupling}.
\\\\
We prove martingale properties for our process. For this purpose, we suppose in addition to Assumptions~\ref{hyp:probamort} and \ref{hyp:tempsinf}, that $t \mapsto A_{\xi}(t)$ is well-defined on $\mathbb{R}^{+}$ for every $\xi>0$. In the following, if $F$ is a $\mathcal{C}^{1,1}$ function from $\mathbb{R}^{+} \times \mathbb{R}^{*}_{+}$ to $\mathbb{R}$, we write $\partial_{1} F$ and $\partial_{2} F$ the partial derivative according to the first and second variable, respectively. The quadratic variation of the process $(X_{t})_{t \geq 0}$ given in the following proposition is a predictable quadratic variation (see Theorem 4.2. in \cite{Jacod1987LimitTF}).
\begin{prop}
Let $F : (t,x) \mapsto F(t,x)$ be $\mathcal{C}^{1,1}$ and $\xi_{0}>0$. Under the general setting of Section~\ref{subsec:gensett}, under Assumptions~\ref{hyp:probamort} and \ref{hyp:tempsinf}, if $t \mapsto A_{\xi}(t)$ is well-defined on $\mathbb{R}^{+}$ for every $\xi>0$ and $R \geq 0$, for all $t \geq 0$, we have:
\begin{align*} 
F(t,X_{t}) = & \hspace{0.1 cm}  F(0,\xi_{0}) + \displaystyle{\int_{0}^{t}\bigg(\partial_{1} F(s,X_{s}) + \partial_{2} F(s,X_{s})g(X_{s},R) \bigg) \mathrm{d}s} \\
& + \displaystyle{\int_{0}^{t}} b_{x_{0}}(X_{s}) \bigg(F(s,X_{s}-x_{0})- F(s,X_{s})\bigg) \mathrm{d}s + M_{F,t},
\end{align*}
where 
\begin{align*}
M_{F,t} := \displaystyle{\int_{0}^{t}\int_{\mathbb{R}^{+}}}\mathbb{1}_{\{ h \leq b_{x_{0}}(X_{s-})\}} \bigg(F(s,X_{s-}-x_{0}) - F(s,X_{s-})\bigg) \mathcal{N}_{C}(\mathrm{d}s, \mathrm{d}h)
\end{align*}
is a martingale with predictable quadratic variation
\begin{align*}
\langle M_{F} \rangle_{t} := \displaystyle{\int_{0}^{t}} b_{x_{0}}(X_{s}) \bigg(F(s,X_{s}-x_{0})- F(s,X_{s})\bigg)^{2} \mathrm{d}s.
\end{align*}
\label{prop:martingaleprobx}
\end{prop}
The proof of Proposition~\ref{prop:martingaleprobx} can be found in Appendix~\ref{app:propsept}. The reader can deduce from this result that the process $(t,X_{t})_{t \geq 0}$ is a solution to the martingale problem for the operator $L$, where for every $\mathcal{C}^{1,1}$ function $F$, $t \geq 0$, $x>0$,
\[ LF := (t,x) \mapsto \partial_{1} F(t,x) + \partial_{2} F(t,x)g(x,R)  +b_{x_{0}}(x)\bigg(F(x-x_{0})-F(x)\bigg). \]
We will design a relevant renormalization and time shift $Z$ of $X$ that verifies a similar martingale problem. Then, we will control the martingale part of $Z$. First, we want to compare $X$ to the evolution of energy without any birth event, so we define
\[ Y_{t} := \dfrac{X_{t}}{A_{\xi_{0}}(t)}. \]
The process $(Y_{t})_{t \geq 0}$ takes values in $[0,1]$.
\\\\
\textbf{Remark:} The rescaling of $X_{t}$ by $A_{\xi_{0}}(t)$ is not appropriate under the allometric setting of Section~\ref{subsec:allomsett} in the case $\alpha>1$, precisely because $A_{\xi_{0}}(.)$ explodes in finite time, so $Y_{t}$ would not be well-defined for all $t \geq 0$. This is one reason why we cannot adapt immediately the reasoning of this section for Theorem~\ref{theo:sharpdeux}.
\begin{lemme}
Let $F$ be a $\mathcal{C}^{1}$ function $\mathbb{R}^{+} \rightarrow \mathbb{R}$ and $\xi_{0}>0$. Under the general setting of Section~\ref{subsec:gensett}, under Assumptions~\ref{hyp:probamort} and \ref{hyp:tempsinf}, if $t \mapsto A_{\xi}(t)$ is well-defined on $\mathbb{R}^{+}$ for every $\xi>0$ and $R\geq 0$, for all $t \geq 0$, we have:
\begin{align*} 
F(Y_{t}) = & \hspace{0.1cm} F(1) + \displaystyle{\int_{0}^{t}\dfrac{F'(Y_{s})}{A_{\xi_{0}}(s)}\bigg( g(A_{\xi_{0}}(s)Y_{s},R) - g(A_{\xi_{0}}(s),R)Y_{s} \bigg) \mathrm{d}s} \\
& + \displaystyle{\int_{0}^{t}\int_{\mathbb{R}^{+}}}\mathbb{1}_{\{ h \leq b_{x_{0}}(A_{\xi_{0}}(s)Y_{s-})\}} \left(F\bigg(Y_{s-}-\dfrac{x_{0}}{A_{\xi_{0}}(s)}\bigg) - F(Y_{s-})\right) \mathcal{N}(\mathrm{d}s, \mathrm{d}h).
\end{align*}
\label{lemm:yt}
\end{lemme}
\begin{proof}
This simply comes from Proposition~\ref{prop:martingaleprobx} applied to $(t,x) \mapsto F\bigg(\dfrac{x}{A_{\xi_{0}}(t)}\bigg)$.
\end{proof}

\subsubsection{Martingale properties under the allometric setting of Section~\ref{subsec:allomsett}}
\label{subsec:limit}

In the following sections, we work under the allometric setting of Section~\ref{subsec:allomsett} with $\alpha \leq 1$, Assumptions~\ref{hyp:tempsinf} and \ref{ass:gainenergy}, and $R \in \mathcal{R}_{0}$ ($\mathcal{R}_{0} \neq \varnothing$ thanks to Lemma~\ref{cestlafindesharicots}). Remark that this implies that $t \mapsto A_{\xi}(t)$ is well-defined on $\mathbb{R}^{+}$ for every $\xi>0$ by Proposition~\ref{eq:gammalpha} and \eqref{eq:indivenergymod}. Also, $R \in \mathcal{R}_{0}$ so $R \notin \mathfrak{R}_{0}$, and we do not need to add Assumption~\ref{hyp:probamort} for the process $(X_{t})_{t \geq 0}$ of Section~\ref{subsec:construpta} to be biologically relevant. Thus, for this special choice of $R$, we will freely use results from Section~\ref{subsec:construpta}, where Assumption~\ref{hyp:probamort} was necessary only to ensure that $T_{0}=+ \infty$. We consider the following change of time scale:
\[ Z_{t}:= Y_{\pi(t)}, \]
with $\pi'(t) = \dfrac{A_{\xi_{0}}(\pi(t))}{g(A_{\xi_{0}}(\pi(t)),R)}>0$  and $\pi(0)=0$ ($\pi'$ is well-defined thanks to the choice of $R$). In our setting, $g(x,R) = C_{R}x^{\alpha}$ as in Equation~\eqref{eq:indivenergymod}, thanks to Proposition~\ref{eq:gammalpha}. We also have the precise expression:
\begin{equation}
\begin{split}
\forall (\xi_{0},t) \in (\mathbb{R}^{+})^{2}, \quad 
A_{\xi_{0}}(t) = \left\{
    \begin{array}{ll}
        \left( (1-\alpha)C_{R}t + \xi_{0}^{1-\alpha} \right)^{1/1-\alpha} & \mbox{if} \hspace{0.1cm} \alpha <1, \\
        \xi_{0}e^{C_{R}t} & \mbox{if} \hspace{0.1 cm} \alpha=1.
    \end{array}
\right.
\end{split}
\label{eq:flowsol}
\end{equation} 
This leads to $\pi'(t) = A_{\xi_{0}}(\pi(t))^{1-\alpha}/C_{R} $, so
\begin{equation}
\begin{split}
\pi'(t)= \left\{
    \begin{array}{ll}
       (1-\alpha)\pi(t) + \xi_{0}^{1-\alpha}/C_{R} & \mbox{if} \hspace{0.1cm} \alpha <1, \\
        1/C_{R} & \mbox{if} \hspace{0.1 cm} \alpha=1.
    \end{array}
\right.
\end{split}
\label{eq:piprime}
\end{equation}   
Hence 
\begin{equation}
\begin{split}
\pi : t \mapsto \left\{
    \begin{array}{ll}
       \dfrac{\xi_{0}^{1-\alpha}}{C_{R}(1-\alpha)}(e^{(1-\alpha)t}-1) & \mbox{if} \hspace{0.1cm} \alpha <1, \\
        t/C_{R} & \mbox{if} \hspace{0.1 cm} \alpha=1.
    \end{array}
\right.
\end{split}
\label{eq:pipipi}
\end{equation}
In both cases, we obtain $A_{\xi_{0}}(\pi(u)) = \xi_{0}e^{u}$. Also, $\pi$ is a deterministic bijection from $\mathbb{R}^{+}$ to $\mathbb{R}^{+}$. To understand better the behavior of $Z$, we define another Poisson point measure adapted to our modified time scale. Let $(s_{i},h_{i})_{i \in \mathbb{N}}$ be an enumeration of the random points of $\mathcal{N}$ such that
$$ \mathcal{N} = \sum_{i \in \mathbb{N}}{\delta_{s_{i},h_{i}}}, $$
we define
$$ \widehat{\mathcal{N}} = \sum_{i \in \mathbb{N}}{\delta_{\pi^{-1}(s_{i}),h_{i}\pi'(\pi^{-1}(s_{i}))}}. $$
\begin{lemme}
$\widehat{\mathcal{N}}$ is a Poisson point measure, with intensity $\mathrm{d}s\mathrm{d}h$. Hence, $\mathcal{N}$ and $\widehat{\mathcal{N}}$ have the same law.
\label{lemm:equapoi}
\end{lemme}
\begin{proof}
From the definition of $\widehat{\mathcal{N}}$, for any $A$ measurable set of $\mathbb{R}^{+} \times \mathbb{R}^{+}$, we have $\widehat{\mathcal{N}}(A) = \mathcal{N}(\hat{A}),$ where $$\hat{A} := \{ (s,h), \hspace{0.1 cm} (\pi^{-1}(s),h\pi'(\pi^{-1}(s))) \in A \}.$$ 
Thanks to the fact that $\pi$ is a bijection, $\widehat{\mathcal{N}}$ is a Poisson point measure with some intensity measure $\lambda$ on $\mathbb{R}^{2}_{+}$. If we take now $A = [s_{1},s_{2}] \times [h_{1},h_{2}]$, we obtain $$\hat{A} = \bigg\{ \{s\} \times \left[\dfrac{h_{1}}{\pi'(\pi^{-1}(s))},\dfrac{h_{2}}{\pi'(\pi^{-1}(s))}\right], \hspace{0.1 cm} \pi(s_{1}) \leq s \leq \pi(s_{2}) \bigg\},$$
because $\pi$ is positive increasing. Hence,
\begin{align*} 
\lambda(A)=\mathrm{Leb}(\hat{A}) & = \displaystyle{\int_{\pi(s_{1})}^{\pi(s_{2})} \int_{h_{1}/\pi'(\pi^{-1}(s))}^{h_{2}/\pi'(\pi^{-1}(s))} \mathrm{d}h\mathrm{d}s } = \displaystyle{\int_{s_{1}}^{s_{2}} \int_{h_{1}}^{h_{2}} \mathrm{d}h\mathrm{d}u, }
\end{align*}
by the change of variables $u=\pi^{-1}(s)$, where $\mathrm{Leb}$ denotes the Lebesgue measure on $\mathbb{R}^{+}$. $\mathcal{N}$ and $\widehat{\mathcal{N}}$ thus have same intensity.
\end{proof}

\noindent In the following, we write $\widehat{\mathcal{N}}_{C}$ for the compensated Poisson point measure associated to $\widehat{\mathcal{N}}$.
\begin{lemme}
Let $F$ be a $\mathcal{C}^{1}$ function $\mathbb{R}^{+} \rightarrow \mathbb{R}$ and $\xi_{0}>0$. Under the general setting of Section~\ref{subsec:gensett}, under Assumptions~\ref{hyp:probamort} and \ref{hyp:tempsinf}, if $t \mapsto A_{\xi}(t)$ is well-defined on $\mathbb{R}^{+}$ for every $\xi>0$ and $R \in \mathcal{R}_{0}$, for all $t \geq 0$, we have:
\begin{align*} 
F(Z_{t})= &  \hspace{0.1cm} F(1) + \displaystyle{\int_{0}^{t}F'(Z_{u})\left( \dfrac{g(A_{\xi_{0}}(\pi(u))Z_{u},R)}{g(A_{\xi_{0}}(\pi(u)),R)}-Z_{u} \right) \mathrm{d}u} \\
& + \displaystyle{\int_{0}^{t} \pi'(u) b_{x_{0}}\left(A_{\xi_{0}}(\pi(u))Z_{u}\right)}  \left(F\bigg(Z_{u}-\dfrac{x_{0}}{A_{\xi_{0}}(\pi(u))}\bigg) - F(Z_{u})\right) \mathrm{d}u + \mathcal{M}_{F,t},
\end{align*}
where
\begin{align*}
\mathcal{M}_{F,t} := \displaystyle{\int_{0}^{t}\int_{\mathbb{R}^{+}}}\mathbb{1}_{\{ h \leq \pi'(u) b_{x_{0}}\left(A_{\xi_{0}}(\pi(u))Z_{u-}\right) \}} \left(F\bigg(Z_{u-}-\dfrac{x_{0}}{A_{\xi_{0}}(\pi(u))}\bigg) - F(Z_{u-})\right) \widehat{\mathcal{N}}_{C}(\mathrm{d}u, \mathrm{d}h)
\end{align*}
is a $L^{2}$-martingale with predictable quadratic variation
\begin{align*}
\langle \mathcal{M}_{F} \rangle_{t} := \displaystyle{\int_{0}^{t}} \pi'(u) b_{x_{0}}\left(A_{\xi_{0}}(\pi(u))Z_{u}\right) \left(F\bigg(Z_{u}-\dfrac{x_{0}}{A_{\xi_{0}}(\pi(u))}\bigg)- F(Z_{u})\right)^{2} \mathrm{d}u.
\end{align*}
\label{corr:martzt}
\end{lemme}
\noindent The proof of this lemma can be found in Appendix~\ref{app:lemmamr}. Under our setting, it translates into:
\begin{multline*} 
F(Z_{t})= \hspace{0.1cm} F(1) + \displaystyle{\int_{0}^{t}F'(Z_{u})\left( Z_{u}^{\alpha}-Z_{u} \right) \mathrm{d}u} \\
+ \displaystyle{\int_{0}^{t} \pi'(u) b_{x_{0}}\left(\xi_{0}e^{u}Z_{u}\right)}  \left(F\bigg(Z_{u}-\dfrac{x_{0}}{\xi_{0}}e^{-u}\bigg) - F(Z_{u})\right) \mathrm{d}u + \mathcal{M}_{F,t}.
\end{multline*}
\subsubsection{Proof of point 8. in Section~\ref{theo:sharp} for $\beta < \alpha$}
\label{subsec:core}

In this section, we study the case $\beta < \alpha$. We focus on the martingale part $\mathcal{M}_{\mathrm{Id}}$ of $Z$, where $\mathrm{Id}$ is the identity function on $\mathbb{R}^{+}$, and then deduce a useful control of this process. Lemma~\ref{corr:martzt} for $F= \mathrm{Id}$, along with Equations~\eqref{eq:piprime} and \eqref{eq:pipipi} gives

\begin{align}
Z_{t} & = \hspace{0.1cm} 1 + \displaystyle{\int_{0}^{t}\left( Z_{u}^{\alpha}-Z_{u} \right) \mathrm{d}u}- \dfrac{x_{0}}{\xi_{0}} \displaystyle{\int_{0}^{t} \dfrac{\xi_{0}^{1-\alpha}e^{(1-\alpha)u}}{C_{R}} b_{x_{0}}\left(\xi_{0}e^{u}Z_{u}\right)} e^{-u} \mathrm{d}u + \mathcal{M}_{\mathrm{Id},t},
\label{eq:ztid}
\end{align} 
with 
\begin{align*}
\mathcal{M}_{\mathrm{Id},t} = - \displaystyle{\int_{0}^{t}\int_{\mathbb{R}^{+}}}\mathbb{1}_{\left\{ h \leq \frac{\xi_{0}^{1-\alpha}e^{(1-\alpha)u}}{C_{R}} b_{x_{0}}\left(\xi_{0}e^{u}Z_{u-}\right) \right\}} \dfrac{x_{0}}{\xi_{0}e^{u}} \widehat{\mathcal{N}}_{C}(\mathrm{d}u, \mathrm{d}h)
\end{align*}
and
\begin{align}
\langle \mathcal{M}_{\mathrm{Id}} \rangle_{t} = \displaystyle{\int_{0}^{t}} \dfrac{\xi_{0}^{1-\alpha}e^{(1-\alpha)u}}{C_{R}} b_{x_{0}}\left(\xi_{0}e^{u}Z_{u}\right) \left(\dfrac{x_{0}}{\xi_{0}e^{u}}\right)^{2} \mathrm{d}u.
\label{eq:crochet}
\end{align}

\begin{lemme}
Let $\varepsilon>0$, $x_{0} >0$, $R \in \mathcal{R}_{0}$. For all $\xi_{0}$ high enough, we have $$ \mathbb{P}_{x_{0},R,\xi_{0}}(\forall t, \hspace{0.1cm} |\mathcal{M}_{\mathrm{Id},t}| < \varepsilon) > 0.$$
\label{lemm:martingale}
\end{lemme}

\begin{proof}
We take any $(t,x_{0},\xi_{0})$. By definition of the predictable quadratic variation, $\mathcal{M}_{\mathrm{Id}}^{2} - \langle \mathcal{M}_{\mathrm{Id}} \rangle$ is a local martingale. It is even a martingale, because of Lemma~\ref{corr:martzt}, and $\mathbb{E}(\mathcal{M}_{\mathrm{Id},t}^{2} - \langle \mathcal{M}_{\mathrm{Id}} \rangle_{t})= \mathbb{E}(\mathcal{M}_{\mathrm{Id},0}^{2} - \langle \mathcal{M}_{\mathrm{Id}} \rangle_{0}) = 0$. We use Doob's maximal inequality in its $L^{2}$ version (Corollary II.1.6 in \cite{revuz2004continuous}) to assess that
\begin{align}
\mathbb{P}_{x_{0},R,\xi_{0}}(\sup_{0 \leq s \leq t}|\mathcal{M}_{\mathrm{Id},s}| \geq \varepsilon) & \leq \dfrac{ \mathbb{E}( \mathcal{M}_{\mathrm{Id},t}^{2})}{\varepsilon^{2}} = \dfrac{ \mathbb{E}( \langle \mathcal{M}\rangle_{\mathrm{Id},t})}{\varepsilon^{2}}.
\label{eq:ineqmaxim}
\end{align}

Now, we consider \eqref{eq:crochet}. Remark that $0 \leq Z_{u} \leq 1$ (because this is the case for $Y$), and $b_{x_{0}}(x) = \mathbb{1}_{\{x > x_{0}\}} C_{\beta}x^{\beta}$. Hence, for $u \geq 0$, $b_{x_{0}}\left(\xi_{0}e^{u}Z_{u}\right) \leq C_{\beta}\xi_{0}^{\beta}e^{\beta u}$. Finally, we have a deterministic upper bound for $\langle \mathcal{M} \rangle_{\mathrm{Id},t}$, so from Equations~\eqref{eq:crochet} and \eqref{eq:ineqmaxim}, we get
\begin{align*}
\mathbb{P}_{x_{0},R,\xi_{0}}(\sup_{0 \leq s \leq t}|\mathcal{M}_{\mathrm{Id},s}| \geq \varepsilon)& \leq  \dfrac{C_{\beta}\xi_{0}^{\beta-\alpha-1}x_{0}^{2}}{C_{R} \varepsilon^{2}} \displaystyle{\int_{0}^{t}} e^{(\beta-\alpha-1)u} \mathrm{d}u.
\end{align*}
Since $\beta < \alpha$, we finally obtain $\mathbb{P}_{x_{0},R,\xi_{0}}(\sup_{0 \leq s \leq t}|\mathcal{M}_{\mathrm{Id},s}| \geq \varepsilon) \leq \dfrac{C_{\beta}\xi_{0}^{\beta-\alpha-1}x_{0}^{2}}{C_{R} \varepsilon^{2}(1+\alpha-\beta)}$, which does not depend on $t$, and converges to 0 as $\xi_{0}$ goes to $+ \infty$. Hence, for $\xi_{0}$ high enough, $$ \mathbb{P}_{x_{0},R,\xi_{0}}(\exists t, \hspace{0.1cm} |\mathcal{M}_{\mathrm{Id},t}| \geq \varepsilon) < 1.$$
\end{proof}
\noindent This uniform control of the martingale part of $Z$ allows us to prove the following proposition.
\begin{prop}
Let $\varepsilon>0$, $x_{0} >0$, $R \in \mathcal{R}_{0}$. For $\xi_{0}$ high enough, we have $$ \mathbb{P}_{x_{0},R,\xi_{0}}(\forall t, \hspace{0.1cm} Z_{t} \geq 1-\varepsilon) > 0.$$
\label{prop:zt}
\end{prop}

\begin{proof}
We choose $\xi_{0}$ high enough so that the result of Lemma~\ref{lemm:martingale} holds true for $\varepsilon/2$. Using \eqref{eq:ztid}, we obtain that on the event $\{ \forall t, \hspace{0.1cm} |\mathcal{M}_{\mathrm{Id},t}| < \varepsilon/2 \}$,
\begin{align*} 
Z_{t} & \geq 1 + \displaystyle{\int_{0}^{t}\left( Z_{u}^{\alpha}-Z_{u} \right) \mathrm{d}u} - \dfrac{C_{\beta}\xi_{0}^{\beta-\alpha}x_{0}}{C_{R}} \displaystyle{\int_{0}^{t} e^{(\beta-\alpha)u}  \mathrm{d}u} + \mathcal{M}_{\mathrm{Id},t} \\
& \geq 1- \varepsilon/2 - \dfrac{C_{\beta}\xi_{0}^{\beta-\alpha}x_{0}}{C_{R}(\alpha-\beta)} + \displaystyle{\int_{0}^{t}\left( Z_{u}^{\alpha}-Z_{u} \right) \mathrm{d}u} \\
& \geq 1- \varepsilon/2 - \dfrac{C_{\beta}\xi_{0}^{\beta-\alpha}x_{0}}{C_{R}(\alpha-\beta)},
\end{align*}
because $Z_{u} \leq 1$ for all $u$. We simply choose $\xi_{0}$ higher if necessary to enforce $\dfrac{C_{\beta}\xi_{0}^{\beta-\alpha}x_{0}}{C_{R}(\alpha-\beta)} \leq \varepsilon/2$.
\end{proof}
\noindent Now, we can try to evaluate what happens if we bring back death events. Let $\varepsilon >0$, $R \in \mathcal{R}_{0}$ and $x_{0} > 0$, thanks to Proposition~\ref{prop:zt}, we choose $\xi_{0}$ high enough so that the event $\{ \forall t, \hspace{0.1cm} Z_{t} \geq 1-\varepsilon \}$ occurs with positive probability, and we work on this event, which can be rewritten, as $\pi$ is a bijection, as
\begin{align}
\left\{
    \begin{array}{rl}
       \left\{ \forall t, \hspace{0.1cm} X_{t} \geq (1-\varepsilon)\left(\xi_{0}^{1-\alpha}+C_{R}(1-\alpha)t\right)^{1/(1-\alpha)} \right\} & \mathrm{if} \hspace{0.1 cm} \alpha <1, \\
        \left\{ \forall t, \hspace{0.1cm} X_{t} \geq (1-\varepsilon)\xi_{0}e^{C_{R}t} \right\} & \mathrm{if} \hspace{0.1 cm} \alpha =1.
    \end{array}
\right.
\label{eq:xtsup}
\end{align}
In order to prove point \textbf{8.} of Section~\ref{theo:sharp} for $\beta < \alpha$, we suppose by contradiction that $\delta < \alpha -1 \leq 0$. Then, $d$ is decreasing. In the following, $c$ is a constant depending on $\xi_{0}$ and $\varepsilon$, possibly varying from line to line. We first assume $\alpha<1$. Observe that on the event $\{ \forall t, \hspace{0.1cm} Z_{t} \geq 1-\varepsilon \}$, we get for $t \geq 0$,
\begin{align*}
\int_{0}^{t}d(X_{s})\mathrm{d}s & \leq \int_{0}^{t}d\left(c(1+s)^{1/(1-\alpha)}\right)\mathrm{d}s \leq \int_{0}^{t}c(1+s)^{\delta/(1-\alpha)}\mathrm{d}s.
\end{align*}
As we assumed that $\delta/(1-\alpha) < -1$, the above integral converges when $t \rightarrow + \infty$, and there exists a constant, still denoted as $c$, such that $\int_{0}^{t}d(X_{s})\mathrm{d}s < c$ with positive probability. One can check that the same result holds if $\alpha=1$. Finally, from the definition of $T_{d}$ in~\eqref{eq:definitiontd}, on the event $\{ \forall t, \hspace{0.1cm} Z_{t} \geq 1-\varepsilon \}$ (so on the event $\{ \int_{0}^{t}d(X_{s})\mathrm{d}s < c \}$  by the previous reasoning), we have $\{c <E \} \subseteq \{ T_{d}=+ \infty \}$, so we obtain 
\begin{align*}
\mathbb{P}_{x_{0},R,\xi_{0}}(T_{d} = + \infty) & \geq \mathbb{P}_{x_{0},R,\xi_{0}}(\{T_{d} = + \infty\} \cap \{ \forall t, \hspace{0.1cm} Z_{t} \geq 1-\varepsilon \}) \\ 
& \geq \mathbb{P}_{x_{0},R,\xi_{0}}\left(\{ c < E \} \cap \left\{\forall t, \hspace{0.1cm} Z_{t} \geq 1-\varepsilon  \right\} \right)\\
& = \mathbb{P}_{x_{0},R,\xi_{0}}(c< E) \mathbb{P}_{x_{0},R,\xi_{0}}\left( \left\{\forall t, \hspace{0.1cm} Z_{t} \geq 1-\varepsilon  \right\} \right),
\end{align*} 
because $E$ is independent from $X$, so from $Z$. This lower bound is positive thanks to the choice of $\xi_{0}$. By Theorem~\ref{theo:cns}, this contradicts Assumption~\ref{hyp:tempsinf} and ends the proof of point \textbf{8.} of Section~\ref{theo:sharp} for $\beta > \alpha$.

\subsubsection{Proof of point \textbf{8.} of Section~\ref{theo:sharp} for $\beta = \alpha$}
\label{subsec:coredeux}
Finally, in the case $\beta = \alpha$, we use another renormalization. Let $\kappa <1$, we define
\[ Y_{\kappa,t} := \dfrac{X_{t}}{A_{\xi_{0}}(t)^{\kappa}}. \]

\begin{lemme}
Let $F$ be a $\mathcal{C}^{1}$ function $\mathbb{R}^{+} \rightarrow \mathbb{R}$. Under the allometric setting of Section~\ref{subsec:allomsett} with $\alpha \leq 1$, Assumptions~\ref{hyp:tempsinf} and \ref{ass:gainenergy}, for all $R \in \mathcal{R}_{0}$, for all $t \geq 0$, we have:
\begin{align*} 
F(Y_{\kappa,t}) = & \hspace{0.1cm} F(\xi_{0}^{1-\kappa}) + \displaystyle{\int_{0}^{t}\dfrac{F'(Y_{\kappa,s})}{A_{\xi_{0}}(s)^{\kappa}}\bigg( g(A_{\xi_{0}}(s)^{\kappa}Y_{\kappa,s},R) - \kappa A_{\xi_{0}}(s)^{\kappa-1} g(A_{\xi_{0}}(s),R)Y_{\kappa,s} \bigg) \mathrm{d}s} \\
& + \displaystyle{\int_{0}^{t}\int_{\mathbb{R}^{+}}} b_{x_{0}}(A_{\xi_{0}}(s)^{\kappa}Y_{\kappa,s}) \left(F\bigg(Y_{\kappa,s}-\dfrac{x_{0}}{A_{\xi_{0}}(s)^{\kappa}}\bigg) - F(Y_{\kappa,s})\right) \mathrm{d}s \\
& + \displaystyle{\int_{0}^{t}\int_{\mathbb{R}^{+}}}\mathbb{1}_{\{ h \leq b_{x_{0}}(A_{\xi_{0}}(s)^{\kappa}Y_{\kappa,s-})\}} \left(F\bigg(Y_{\kappa,s-}-\dfrac{x_{0}}{A_{\xi_{0}}(s)^{\kappa}}\bigg) - F(Y_{\kappa,s-})\right) \mathcal{N}_{C}(\mathrm{d}s, \mathrm{d}h).
\end{align*}
\label{lemm:yyy}
\end{lemme}
\begin{proof}
This is Proposition~\ref{prop:martingaleprobx} applied to $(t,x) \mapsto F\bigg(\dfrac{x}{A_{\xi_{0}}^{\kappa}(t)}\bigg)$.
\end{proof}
\noindent Applying Lemma~\ref{lemm:yyy} to $F := \mathrm{Id}$ leads to
\begin{align*} 
Y_{\kappa,t} = & \hspace{0.1cm} \xi_{0}^{1-\kappa} + \displaystyle{\int_{0}^{t}\dfrac{1}{A_{\xi_{0}}(s)^{\kappa}}\bigg( g(A_{\xi_{0}}(s)^{\kappa}Y_{\kappa,s},R) - \kappa A_{\xi_{0}}(s)^{\kappa-1} g(A_{\xi_{0}}(s),R)Y_{\kappa,s} \bigg) \mathrm{d}s} \\
& - \displaystyle{\int_{0}^{t}\int_{\mathbb{R}^{+}}} b_{x_{0}}(A_{\xi_{0}}(s)^{\kappa}Y_{\kappa,s}) \dfrac{x_{0}}{A_{\xi_{0}}(s)^{\kappa}} \mathrm{d}s - \displaystyle{\int_{0}^{t}\int_{\mathbb{R}^{+}}}\mathbb{1}_{\{ h \leq b_{x_{0}}(A_{\xi_{0}}(s)^{\kappa}Y_{\kappa,s-})\}} \dfrac{x_{0}}{A_{\xi_{0}}(s)^{\kappa}} \mathcal{N}_{C}(\mathrm{d}s, \mathrm{d}h) \\
\geq & \hspace{0.1cm} \xi_{0}^{1-\kappa} + \displaystyle{\int_{0}^{t}\bigg( (C_{R}-x_{0}C_{\beta}) A_{\xi_{0}}(s)^{\kappa (\alpha-1)} Y_{\kappa,s}^{\alpha} - C_{R}\kappa A_{\xi_{0}}(s)^{\alpha-1} Y_{\kappa,s} \bigg) \mathrm{d}s} \\
& - \displaystyle{\int_{0}^{t}\int_{\mathbb{R}^{+}}}\mathbb{1}_{\{ h \leq b_{x_{0}}(A_{\xi_{0}}(s)^{\kappa}Y_{\kappa,s-})\}} \dfrac{x_{0}}{A_{\xi_{0}}(s)^{\kappa}} \mathcal{N}_{C}(\mathrm{d}s, \mathrm{d}h).
\end{align*}
We choose $x_{0}$ small enough such that $C_{R}-x_{0}C_{\beta} \geq  C_{R}\kappa,$ and we obtain, for all $t >0$:
\begin{align*}
Y_{\kappa,t} \geq \xi_{0}^{1-\kappa} + \displaystyle{\int_{0}^{t}C_{R}\kappa A_{\xi_{0}}(s)^{\alpha-1} \bigg(  A_{\xi_{0}}(s)^{(\kappa-1)(\alpha-1)}Y_{\kappa,s}^{\alpha} - Y_{\kappa,s} \bigg) \mathrm{d}s} -\tilde{\mathcal{M}}_{t},
\end{align*} 
with $\tilde{\mathcal{M}}_{t} := \displaystyle{\int_{0}^{t}\int_{\mathbb{R}^{+}}}\mathbb{1}_{\{ h \leq b_{x_{0}}(A_{\xi_{0}}(s)^{\kappa}Y_{\kappa,s-})\}} \dfrac{x_{0}}{A_{\xi_{0}}(s)^{\kappa}} \mathcal{N}_{C}(\mathrm{d}s, \mathrm{d}h)$. Also, we know that $X_{t} \leq A_{\xi_{0}}(t)$, so $Y_{\kappa,t} \leq A_{\xi_{0}}(t)^{1-\kappa}$ for $t \geq 0$. We verify that for every $x>0$, the function $y \mapsto x^{1-\alpha}y^{\alpha}-y$ is non-negative on $]0,x]$, so for every $s \geq 0$:
$$A_{\xi_{0}}(s)^{(\kappa-1)(\alpha-1)}Y_{\kappa,s}^{\alpha} - Y_{\kappa,s} \geq 0. $$
This entails that for $t \geq 0$ and $x_{0}$ small enough,
\begin{align}
Y_{\kappa,t} \geq \xi_{0}^{1-\kappa} -\tilde{\mathcal{M}}_{t}.
\label{eq:swag}
\end{align} 
\begin{lemme}
Let $\varepsilon>0$, $x_{0}>0$, $R \in \mathcal{R}_{0}$ and $\kappa > \dfrac{1}{2}$. Then we can take $\xi_{0}$ high enough to have $$ \mathbb{P}_{x_{0},R,\xi_{0}}(\forall t, \hspace{0.1cm} |\tilde{\mathcal{M}}_{t}| < \varepsilon) > 0.$$
\label{lemme:martkappa}
\end{lemme}
\begin{proof}
We take $t > 0$ and $\alpha <1$. As in the proof of Lemma~\ref{lemm:martingale}, we use Doob's maximal inequality in its $L^{2}$ version to obtain
\begin{align*}
\mathbb{P}_{x_{0},R,\xi_{0}}(\sup_{0 \leq s \leq t}|\tilde{\mathcal{M}}_{s}| \geq \varepsilon) & \leq \dfrac{\mathbb{E}(\langle \tilde{\mathcal{M}} \rangle_{t})}{\varepsilon^{2}} \\
& \leq \dfrac{x_{0}^{2}}{\varepsilon^{2}} \displaystyle{\int_{0}^{t}}  \dfrac{b_{x_{0}}(A_{\xi_{0}}(s)^{\kappa}Y_{\kappa,s})}{A_{\xi_{0}}(s)^{2 \kappa}} \mathrm{d}s \\
& \leq \dfrac{C_{\beta}x_{0}^{2}}{\varepsilon^{2}} \displaystyle{\int_{0}^{t}}  A_{\xi_{0}}(s)^{\alpha - 2 \kappa} \mathrm{d}s,
\end{align*}
as $b_{x_{0}}(x)\leq C_{\beta}x^{\beta}$, $Y_{\kappa,t} \leq A_{\xi_{0}}(t)^{1-\kappa}$ and $\beta = \alpha$. Also, we choose $\kappa > 1/2$, which implies, using Equation~\eqref{eq:flowsol},
\begin{align*}
\mathbb{P}_{x_{0},\xi_{0}}(\sup_{0 \leq s \leq t}|\tilde{\mathcal{M}}_{s}| \geq \varepsilon)& \leq \dfrac{C_{\beta}x_{0}^{2}\xi_{0}^{1-2\kappa}}{\varepsilon^{2}C_{R}(2\kappa - 1)},
\end{align*}
which does not depend on $t$, and converges to 0 as $\xi_{0}$ goes to $+ \infty$. Hence, for $\xi_{0}$ high enough, $$ \mathbb{P}_{x_{0},\xi_{0}}(\exists t, \hspace{0.1cm} |\tilde{\mathcal{M}}_{t}| \geq \varepsilon) < 1.$$
One can check that in the case $\alpha =1$, we obtain the same result.
\end{proof}

\begin{corr}
Let $\varepsilon>0$, $R \in \mathcal{R}_{0}$ and $\kappa > 1/2$. For $x_{0},\xi_{0}$ respectively small and high enough, we have $$ \mathbb{P}_{x_{0},\xi_{0}}(\forall t, \hspace{0.1cm} Y_{\kappa, t} \geq \xi_{0}^{1-\kappa}-\varepsilon) > 0.$$
\label{prop:ykappat}
\end{corr}

\begin{proof}
Let $x_{0}$ be small enough such that \eqref{eq:swag} holds true, and let $\xi_{0}$ be high enough so that the result of Lemma~\ref{lemme:martkappa} holds true for such $x_{0}$. Then, on the event $\{ \forall t, \hspace{0.1cm} |\tilde{\mathcal{M}}_{t}| < \varepsilon \}$, we have thanks to~\eqref{eq:swag}, for all $t \geq 0$,
\begin{align*} 
Y_{\kappa,t} \geq \xi_{0}^{1-\kappa}-\varepsilon.
\end{align*}
\end{proof}
\noindent In the same way as the previous section, we evaluate what happens if we bring back death events. Let $\varepsilon >0$  and $\kappa >1/2$. Thanks to Corollary~\ref{prop:ykappat}, we choose $x_{0}$ small enough and $\xi_{0}$ high enough (and even higher if necessary to have $\xi_{0}^{1-\kappa} \geq 1$), so that the event $\{ \forall t, \hspace{0.1cm} Y_{\kappa,t} \geq 1-\varepsilon \}$ occurs with positive probability, and we work on this event, which can be rewritten from Equation~\eqref{eq:flowsol} as 
\begin{align}
\left\{
    \begin{array}{rl}
       \left\{ \forall t, \hspace{0.1cm} X_{t} \geq (1-\varepsilon)\left(\xi_{0}^{1-\alpha}+C_{R}(1-\alpha)t\right)^{\kappa/(1-\alpha)} \right\} & \mathrm{if} \hspace{0.1 cm} \alpha <1, \\
        \left\{ \forall t, \hspace{0.1cm} X_{t} \geq (1-\varepsilon)\xi_{0}^{\kappa}e^{\kappa C_{R}t} \right\} & \mathrm{if} \hspace{0.1 cm} \alpha =1.
    \end{array}
\right.
\label{eq:xsup}
\end{align}
To prove point \textbf{8.} of Section~\ref{theo:sharp} when $\beta = \alpha$, we suppose by contradiction that $\delta < \alpha-1 \leq 0$, so $d$ is decreasing. Once this is done, we also choose $\kappa \in (1/2,1)$ such that $\delta < \frac{\alpha-1}{\kappa}$. In the following, $c$ is a constant depending on $\xi_{0}$ and $\varepsilon$, possibly varying from line to line. We first take $\alpha<1$. On the event~\eqref{eq:xsup}, we get for $t \geq 0$,
\begin{align*}
\int_{0}^{t}d(X_{s})\mathrm{d}s & \leq \int_{0}^{t}d(c(1+s)^{\kappa/(1-\alpha)})\mathrm{d}s \leq \int_{0}^{t}c(1+s)^{\delta \kappa/(1-\alpha)}\mathrm{d}s,
\end{align*}
with $\delta \kappa /(1-\alpha) < -1$. Hence, the above integral converges when $t \rightarrow + \infty$. One can check that the same result holds if $\alpha=1$. We conclude with the same arguments as in the end of Section~\ref{subsec:core}.
\subsection{Converse implication in Theorem~\ref{theo:short}}
\label{subsec:reciprocal}

In this section, we work under the allometric setting of Section~\ref{subsec:allomsett}. Also, we work under~\eqref{eq:iuncase} or~\eqref{eq:recipro}, and they both imply Assumption~\ref{ass:gainenergy} by Proposition~\ref{eq:gammalpha}. Hence, we pick $R_{0}$ verifying the result of Lemma~\ref{cestlafindesharicots}. We begin by two preliminary results on the operator $K_{x_{0},R}\mathbf{1}$ for $x_{0}>0$ and $R>R_{0}$.

\begin{lemme}
Let $x_{0}>0$, $k \geq 0$. Under the allometric setting of Section~\ref{subsec:allomsett} and~\eqref{eq:iuncase}, we have for $R>R_{0}$,
$$ \forall \xi \geq kx_{0}, \quad K_{x_{0},R}^{k}\mathbf{1}(\xi) = \bigg( \dfrac{C_{\beta}}{C_{\beta}+C_{\delta}} \bigg)^{k}.$$
\label{lemme:induvtion}
\end{lemme}

\begin{proof}
It is straightforward to prove this result by induction on $k$, using \eqref{eq:kixzerooo}, where we notice that under \eqref{eq:iuncase}, for $R >R_{0}$, we have $\displaystyle{\int_{0}^{t_{\max}(\xi_{0},R)} (b_{x_{0}}+d)(A_{\xi_{0},R}(u)) \mathrm{d}u }  = + \infty.$
\end{proof}

\begin{lemme}
Let $x_{0}>0$, $k \geq 1$. Under the allometric setting of Section~\ref{subsec:allomsett} and~\eqref{eq:iuncase}, we have for $R>R_{0}$,
$$  K_{x_{0},R}^{k}\mathbf{1}(x_{0}) \geq \bigg( \dfrac{C_{\beta}}{C_{\beta}+C_{\delta}} \bigg)^{k} k^{-\frac{C_{\beta}+C_{\delta}}{C_{R}}}.$$
\label{lemme:lelemme}
\end{lemme}

\begin{proof}
We compute for $k \geq 1$,
\begin{align*}
K_{x_{0},R}^{k}\mathbf{1}(x_{0}) & := \displaystyle{\int_{x_{0}}^{+ \infty}\dfrac{b_{x_{0}}(u)}{g(u,R)}e^{-\displaystyle{ \int_{x_{0}}^{u}}\dfrac{b_{x_{0}}(\tau)+d(\tau)}{g(\tau,R)}\mathrm{d}\tau}K_{x_{0},R}^{k-1}\mathbf{1}(u-x_{0}) \mathrm{d}u } \\
& \geq \displaystyle{\int_{kx_{0}}^{+ \infty}\dfrac{b_{x_{0}}(u)}{g(u,R)}e^{-\displaystyle{ \int_{x_{0}}^{u}}\dfrac{b_{x_{0}}(\tau)+d(\tau)}{g(\tau,R)}\mathrm{d}\tau}\bigg( \dfrac{C_{\beta}}{C_{\beta}+C_{\delta}} \bigg)^{k-1} \mathrm{d}u }
\end{align*}
by Lemma~\ref{lemme:induvtion}. Under the allometric setting of Section~\ref{subsec:allomsett} and~\eqref{eq:iuncase}, a straightforward computation leads to the result.
\end{proof}

\noindent We are now ready to prove the converse implication in Theorem~\ref{theo:short}.

\begin{prop}
Under the allometric setting of Section~\ref{subsec:allomsett}, under Assumptions~\ref{hyp:probamort} and~\ref{hyp:tempsinf}, we have
\begin{center}
~\eqref{eq:recipro} $\Rightarrow$ Assumption~\ref{ass:supercritical}.
\end{center} 
\label{prop:reciprocal}
\end{prop}

\begin{proof}
We pick $x_{0}>0$ and $R> R_{0}$. Under Assumptions~\ref{hyp:probamort} and~\ref{hyp:tempsinf}, by Lemma~\ref{lemme:kixzero} and Lemma~\ref{lemme:lelemme}, we can write 
\begin{align*}
m_{x_{0},R}(x_{0}) := \mathbb{E}(N_{x_{0},R,x_{0}}) = \sum\limits_{k \geq 1} \mathbb{P}_{x_{0},R,x_{0}}(M^{k}) & = \sum\limits_{k \geq 1} K_{x_{0},R}^{k}\mathbf{1}(x_{0}) \\
& \geq \sum\limits_{k \geq 1} \bigg( \dfrac{C_{\beta}}{C_{\beta}+C_{\delta}} \bigg)^{k} k^{-\frac{C_{\beta}+C_{\delta}}{C_{R}}}.
\end{align*}
Now, if~\eqref{eq:recipro}, we have $C_{\beta}+ C_{\delta} < C_{\gamma}-C_{\alpha}$ and we know that $C_{R} \underset{R \rightarrow + \infty}{\longrightarrow} C_{\gamma}-C_{\alpha}$, so we can choose $R>R_{0}$ high enough such that $\frac{C_{\beta}+C_{\delta}}{C_{R}} <1$, and the previous lower bound, with the fact that $C_{\beta} > (e-1)C_{\delta}$, gives
$$m_{x_{0},R}(x_{0}) > \sum\limits_{k \geq 1} \bigg( \dfrac{C_{\beta}}{C_{\beta}+C_{\delta}} \bigg)^{k} k^{-1} > \sum\limits_{k \geq 1} \dfrac{(1-e^{-1})^{k}}{k} =1, $$ 
which concludes since this is valid for every $x_{0}>0$.
\end{proof}

\noindent Finally, we prove that~\eqref{eq:iuncase} is not a sufficient condition to obtain Assumption~\ref{ass:supercritical}. We begin again by a preliminary lemma.

\begin{lemme}
Under the allometric setting of Section~\ref{subsec:allomsett}, if~\eqref{eq:iuncase}, then for $x_{0}>0$ and $R>R_{0}$, we have
\begin{align}
\forall k \geq 1, \forall \xi >0, \quad K_{x_{0},R}^{k}\mathbf{1}(\xi) \leq  \bigg(\dfrac{C_{\beta}}{C_{\beta}+C_{\delta}} \bigg)^{k},
\label{eq:touslesk}
\end{align}
and
\begin{align}
\exists \mathfrak{C}(C_{\beta},C_{\delta})<1, \quad K_{x_{0},R}^{2}\mathbf{1}(x_{0}) \leq \mathfrak{C}(C_{\beta},C_{\delta}) \bigg(\dfrac{C_{\beta}}{C_{\beta}+C_{\delta}} \bigg)^{2}.
\label{eq:kegaldeux}
\end{align} 
\label{lemm:trimar}
\end{lemme}

\begin{proof}
It is straightforward to prove~\eqref{eq:touslesk} by induction on $k$. Then, remark that
\begin{equation*}
\begin{split}
K_{x_{0},R}\mathbf{1}(\xi) \left\{
    \begin{array}{ll}
       = \frac{C_{\beta}}{C_{\beta}+C_{\delta}} & \mbox{if} \hspace{0.1cm} \xi \geq x_{0}, \\
        \leq \frac{C_{\beta}}{C_{\beta}+C_{\delta}} \bigg( \dfrac{\xi}{x_{0}} \bigg)^{\frac{C_{\delta}}{C_{\gamma}-C_{\alpha}}} & \mbox{else.}
    \end{array}
\right.
\end{split}
\end{equation*}   
We use this fact and~\eqref{eq:touslesk} to obtain
\begin{equation*}
\begin{split}
K_{x_{0},R}^{2}\mathbf{1}(x_{0}) & \leq \dfrac{C_{\beta}}{C_{\beta}+C_{\delta}} \bigg( \displaystyle{\int_{x_{0}}^{\frac{3}{2}x_{0}}\dfrac{b_{x_{0}}(u)}{g(u,R)}e^{-\displaystyle{ \int_{x_{0}}^{u}}\dfrac{b_{x_{0}}(\tau)+d(\tau)}{g(\tau,R)}\mathrm{d}\tau} \left(\dfrac{1}{2}\right)^{\frac{C_{\delta}}{C_{\gamma}-C_{\alpha}}} \mathrm{d}u }\\
& + \displaystyle{\int_{3/2x_{0}}^{+ \infty}\dfrac{b_{x_{0}}(u)}{g(u,R)}e^{-\displaystyle{ \int_{x_{0}}^{u}}\dfrac{b_{x_{0}}(\tau)+d(\tau)}{g(\tau,R)}\mathrm{d}\tau} \mathrm{d}u }\bigg).
\end{split}
\end{equation*}
Under the allometric setting of Section~\ref{subsec:allomsett} and~\eqref{eq:iuncase}, this gives
\begin{equation}
\begin{split}
K_{x_{0},R}^{2}\mathbf{1}(x_{0}) & \leq \bigg(\dfrac{C_{\beta}}{C_{\beta}+C_{\delta}} \bigg)^{2} \mathfrak{C}(C_{\beta},C_{\delta}),
\end{split}
\label{eq:kxodeux}
\end{equation}
with $\mathfrak{C}(C_{\beta},C_{\delta}) := \bigg( \left(\dfrac{1}{2}\right)^{\frac{C_{\delta}}{C_{\gamma}-C_{\alpha}}} \bigg( 1 - \bigg(\dfrac{2}{3} \bigg)^{\frac{C_{\beta }+C_{\delta}}{C_{\gamma}-C_{\alpha}}} \bigg) +  \bigg(\dfrac{2}{3} \bigg)^{\frac{C_{\beta} + C_{\delta}}{C_{\gamma}-C_{\alpha}}} \bigg) <1$, which concludes.
\end{proof}

\begin{prop}
Under the allometric setting of Section~\ref{subsec:allomsett}, under Assumptions~\ref{hyp:probamort},~\ref{hyp:tempsinf} and~\ref{ass:supercritical}, if $\beta=\delta=\alpha-1$, there exists $\Xi : \mathbb{R}^{+} \rightarrow ]1,+ \infty[$ such that
$$ \dfrac{C_{\beta}}{C_{\delta}} > \Xi(C_{\delta}) > 1.$$
If in addition, $C_{\delta} \geq C_{\gamma}-C_{\alpha}$, then
$$ \Xi(C_{\delta}) > 1.07.$$
In particular,~\eqref{eq:iuncase} is not a sufficient condition to verify Assumption~\ref{ass:supercritical} under Assumptions~\ref{hyp:probamort} and~\ref{hyp:tempsinf}.
\label{prop:conversenot}
\end{prop}

\begin{proof}
We pick $x_{0}>0$ and $R> R_{0}$. Under our assumptions, we necessarily have $C_{\beta}> C_{\delta}$ by point~\textbf{4.} of Section~\ref{theo:sharp}. Under Assumptions~\ref{hyp:probamort} and~\ref{hyp:tempsinf}, by Lemma~\ref{lemme:kixzero} and Lemma~\ref{lemm:trimar}, we can write 
\begin{align*}
m_{x_{0},R}(x_{0}) := \mathbb{E}(N_{x_{0},R,x_{0}}) = \sum\limits_{k \geq 1} \mathbb{P}_{x_{0},R,x_{0}}(M^{k}) & = \sum\limits_{k \geq 1} K_{x_{0},R}^{k}\mathbf{1}(x_{0}) \\
& \leq \sum\limits_{k \geq 1} \bigg( \dfrac{C_{\beta}}{C_{\beta}+C_{\delta}} \bigg)^{k} -(1-\mathfrak{C}(C_{\beta},C_{\delta}) )\bigg(\dfrac{C_{\beta}}{C_{\beta}+C_{\delta}} \bigg)^{2}\\
& = \dfrac{C_{\beta}}{C_{\delta}} -(1-\mathfrak{C}(C_{\beta},C_{\delta}) )\bigg(\dfrac{C_{\beta}}{C_{\beta}+C_{\delta}} \bigg)^{2}.
\end{align*}
For a given value of $C_{\delta}$, this upper bound goes to $1-\frac{L}{4}$, with some $L>0$, when $C_{\beta}$ converges to $C_{\delta}$ (see the explicit expression of $\mathfrak{C}(C_{\beta},C_{\delta}) $ in the proof of Lemma~\ref{lemm:trimar}). Hence, for a given value of $C_{\delta}$, if $C_{\beta}$ is too close from $C_{\delta}$, $m_{x_{0},R}(x_{0}) <1$ so Assumption~\ref{ass:supercritical} is not verified. This is valid for any $C_{\delta}$, hence the result.
\\\\
If in addition, we assume that $C_{\delta} \geq C_{\gamma}-C_{\alpha}$, we can go further in~\eqref{eq:kxodeux} and obtain
\begin{equation*}
\begin{split}
K_{x_{0},R}^{2}\mathbf{1}(x_{0}) & \leq \bigg(\dfrac{C_{\beta}}{C_{\beta}+C_{\delta}} \bigg)^{2}\bigg( \dfrac{1}{2} \bigg( 1 - \bigg(\dfrac{2}{3} \bigg)^{\frac{C_{\beta}}{C_{\delta}}+1} \bigg) +  \bigg(\dfrac{2}{3} \bigg)^{\frac{C_{\beta}}{C_{\delta}}+1} \bigg) \leq \dfrac{13}{18} \bigg(\dfrac{C_{\beta}}{C_{\beta}+C_{\delta}} \bigg)^{2}.
\end{split}
\end{equation*}
Hence, we can pick $\mathfrak{C}(C_{\beta},C_{\delta}) \equiv \frac{13}{18}$ in the previous reasoning, to obtain that 
\begin{align*}
m_{x_{0},R}(x_{0}) \leq \dfrac{C_{\beta}}{C_{\delta}} - \dfrac{5}{18}\bigg(\dfrac{C_{\beta}}{C_{\beta}+C_{\delta}} \bigg)^{2},
\end{align*}
which is less than or equal to 1 if $C_{\beta} \leq 1,07 C_{\delta}$.
\end{proof}

\section{Simulations}
\label{subsec:simulations}

In this section, we work under the allometric setting of Section~\ref{subsec:allomsett} and $\alpha \leq 1$. We present simulations of individual trajectories $(\xi_{t,x_{0},R,\xi_{0}})_{t \geq 0}$, and are mainly interested in the behavior of the average number of offspring $m_{x_{0},R}(x_{0})$, for different values of parameters. 
In particular, our main aim is to determine when we have $m_{x_{0},R}(x_{0})>1$, leading then to the supercriticality of the population process (see Proposition~\ref{prop:galtonwatson}).
Precisely, according to Theorem~\ref{theo:short}, we first fix the following parameters (unless specified otherwise):

\begin{center}
\begin{tabular}{|c|c|}
	\hline 
	Parameter & Value \\
	\hline
	$\alpha$ & 0.75 \\
	$\gamma$ & $\alpha$ \\
	$\delta$ & $\alpha -1$ \\
	$\phi(R)$ & $2/3$ \\
	$C_{\gamma}$ & 2 \\
	$C_{\alpha}$ & 1 \\
	\hline 
\end{tabular}
\end{center}
We picked the usual value $\alpha = 0.75$ highlighted by the Metabolic Theory of Ecology \cite{peters_1983, brown_04, savdee_08}. In this section, we will independently make vary the following parameters: $x_{0}$, $\beta$, $C_{\beta}$ and $C_{\delta}$. We estimate $m_{x_{0},R}(x_{0})$ by a Monte-Carlo method and use \verb+Python+ for the simulations. We present our numerical results, that lead to Conjecture~\ref{conj} on a necessary and sufficient condition to verify Assumption~\ref{ass:supercritical} under Assumptions~\ref{hyp:probamort} and~\ref{hyp:tempsinf}, and Conjecture~\ref{conj:deux} about the behavior of $m_{x_{0},R}(\xi_{0})$ in the $I_{2}$ case.
\\\\
Recall that the allometric coefficients $\beta < \alpha - 1$ are non-admissible because they contradict Assumption~\ref{hyp:tempsinf} (see Section~\ref{subsec:teun}), that is they allow individuals with $T_{d}=+ \infty$. These are obviously trajectories that we cannot observe with numerical simulations, so we focus in the following on the case $\beta \geq \alpha-1$.

\subsection{Typical look of individual trajectories}

\begin{figure}[h!]
\centering
\begin{subfigure}{0.4\textwidth}
    \includegraphics[width=\textwidth]{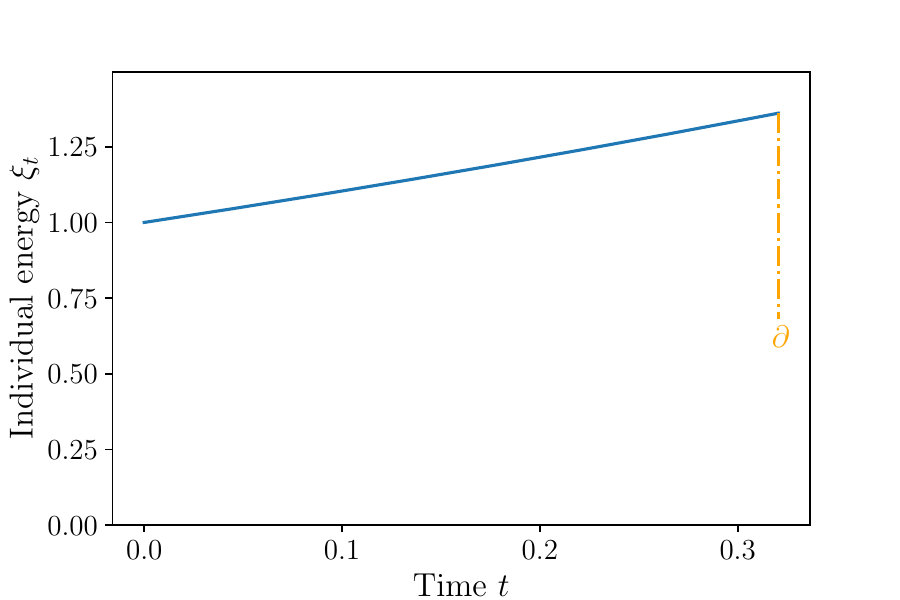}
    \caption{Trajectory 1: no birth events}
    \label{fig:first}
\end{subfigure}
 \hspace{0.5 cm}
\begin{subfigure}{0.4\textwidth}
    \includegraphics[width=\textwidth]{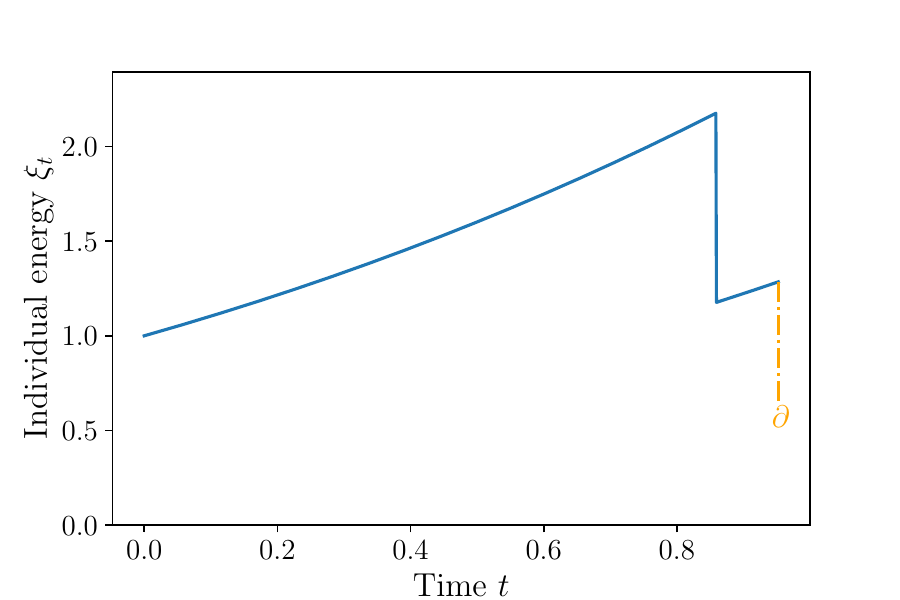} 
    \caption{Trajectory 2: 1 birth event}
    \label{fig:second}
\end{subfigure}
 \hspace{0.5 cm}
\begin{subfigure}{0.4\textwidth}
   \includegraphics[width=\textwidth]{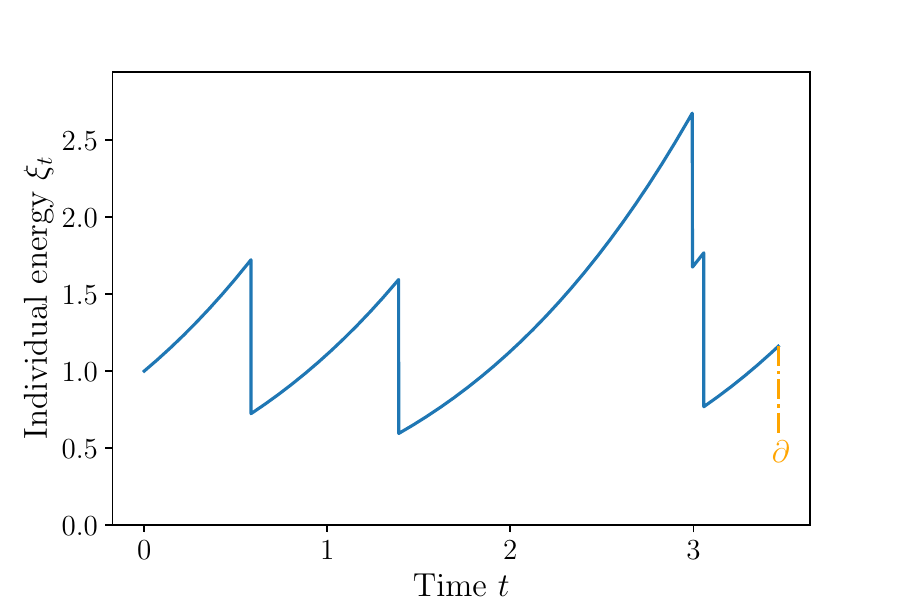}
    \caption{Trajectory 3: 4 birth events}
    \label{fig:third}
\end{subfigure}
\hspace{0.5 cm}
\begin{subfigure}{0.4\textwidth}
    \includegraphics[width=\textwidth]{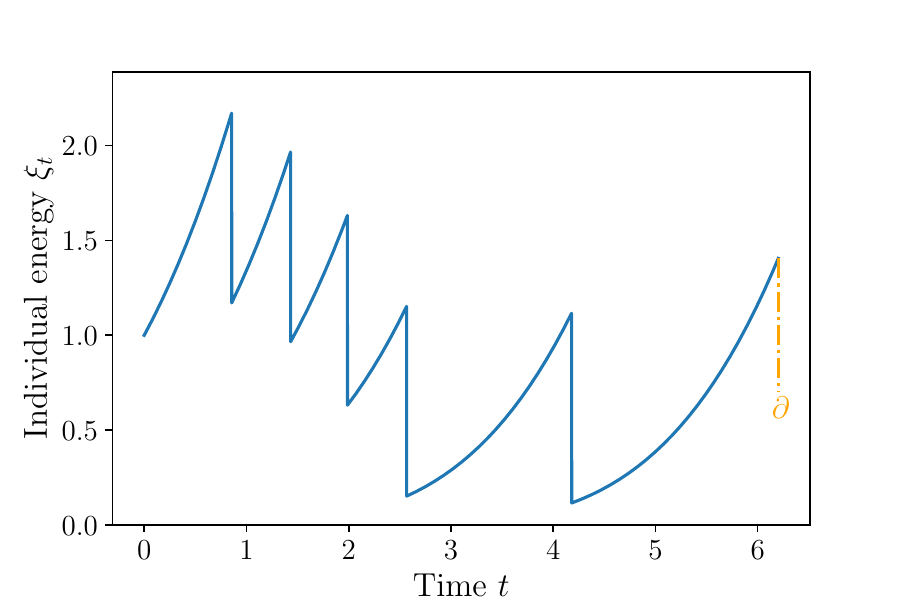}
    \caption{Trajectory 4: 6 birth events}
    \label{fig:fourth}
\end{subfigure}
 \hspace{0.5 cm}
\begin{subfigure}{0.4\textwidth}
    \includegraphics[width=\textwidth]{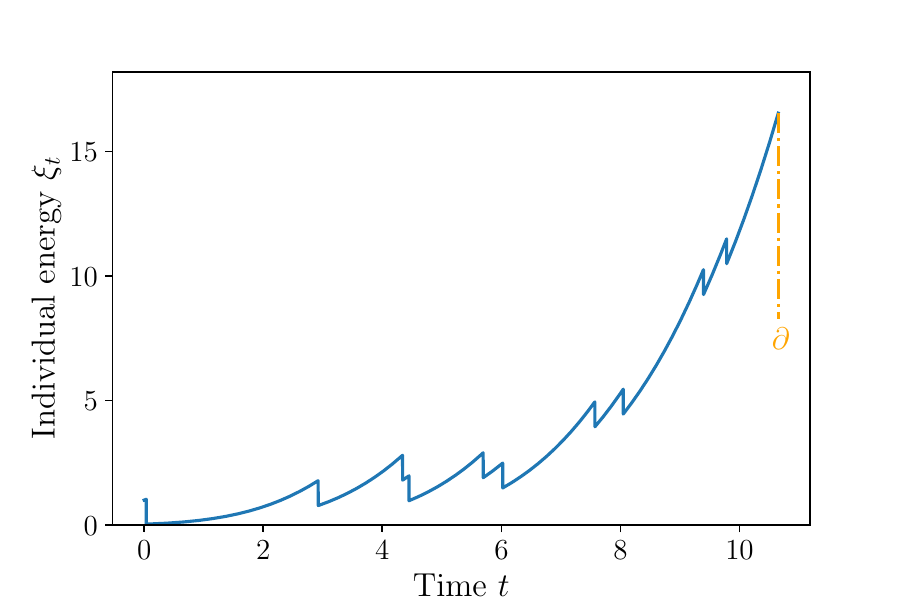} 
    \caption{Trajectory 5: 10 birth events}
    \label{fig:fifth}
\end{subfigure}
 \hspace{0.5 cm}
\begin{subfigure}{0.4\textwidth}
    \includegraphics[width=\textwidth]{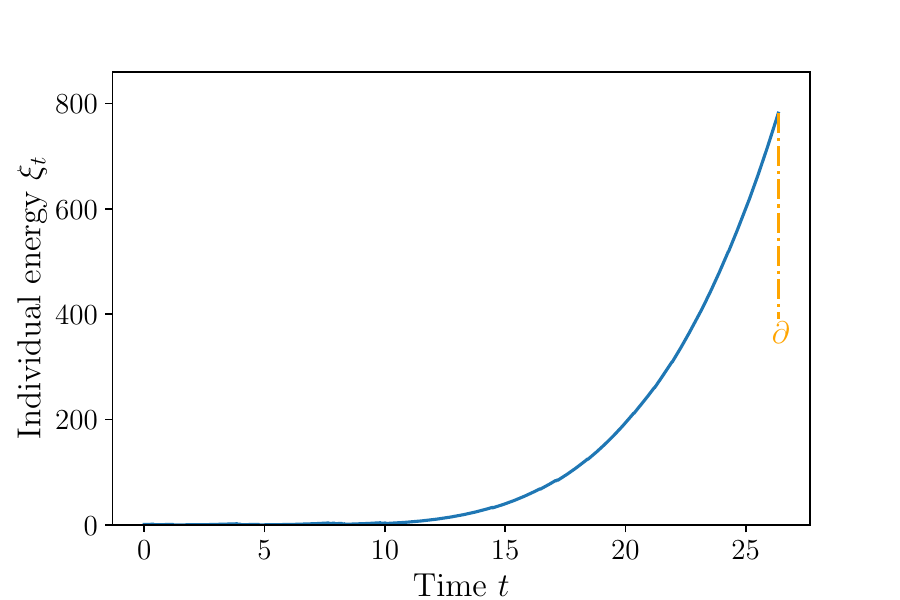}
    \caption{Trajectory 6: 25 birth events}
    \label{fig:sixth}
\end{subfigure}
\caption{Different shapes for individual trajectories, with parameters $x_{0}=\xi_{0}=1$, $\beta=-0.2$, $C_{\beta}=2$, $C_{\delta}=0.5$. Most of the time, we obtain the shapes of \ref{fig:first} and \ref{fig:second}, and the rarest trajectory is \ref{fig:sixth}. The dashed vertical orange line represents the time of death.}
\label{fig:indiv}
\end{figure}

We run independent simulations of the individual process with a set of parameters satisfying $(\gamma, \delta, \beta)\in I_2$ (in particular $\beta>\alpha-1$).
On Figure~\ref{fig:indiv}, we present six individual trajectories, illustrating the three distinct shapes that we observed in our simulations.
The most common trajectories are shown on Figures~\ref{fig:first} and \ref{fig:second}, where individuals die fast and have 0 or 1 offspring. Rarely, lifetime and energy increase and there is a huge amount of birth events like on Figures~\ref{fig:fifth} and \ref{fig:sixth} (on Figure~\ref{fig:sixth}, the energy increases so fast that we barely see the negative jumps due to birth events anymore). Finally, we have intermediate trajectories, as shown on Figures~\ref{fig:third} and \ref{fig:fourth}, where individual energy stays close to the initial energy $x_{0}$. They occur more often than \ref{fig:fifth} and \ref{fig:sixth}, but less often than \ref{fig:first} and \ref{fig:second}. Modifying the parameter $\beta>\alpha-1$ leads to similar observations. Now if $\beta=\alpha-1$, we still encounter the shapes described earlier, but the frequencies of observation of these three different shapes become relatively similar.
\newpage

\subsection{Subcriticality for $\beta > \alpha -1$, not in the $I_{2}$ case}
\begin{figure}[h!]
\centering
\begin{subfigure}{0.4\textwidth}
   \includegraphics[width=\textwidth]{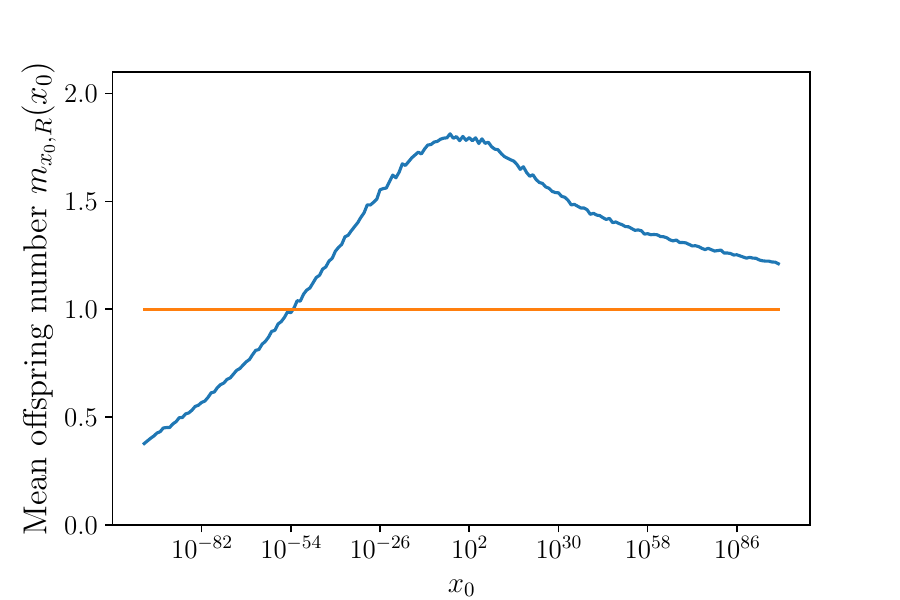} 
    \caption{$\beta=-0.24$}
    \label{fig:firstbeta}
\end{subfigure}
\hspace{0.01 cm}
\begin{subfigure}{0.4\textwidth}
    \includegraphics[width=\textwidth]{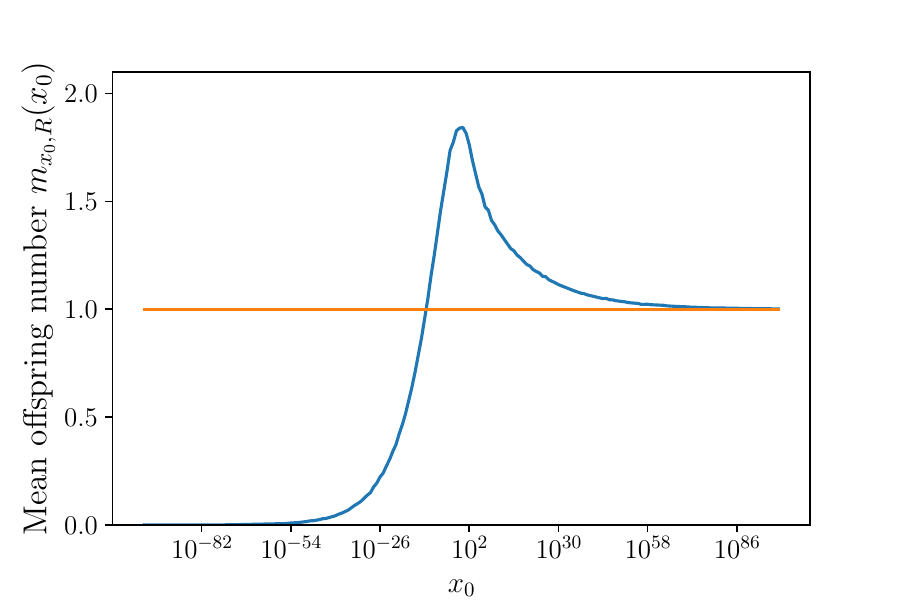}  
    \caption{$\beta=-0.2$}
    \label{fig:secondbeta}
\end{subfigure}
\hspace{0.01 cm}
\begin{subfigure}{0.4\textwidth}
   \includegraphics[width=\textwidth]{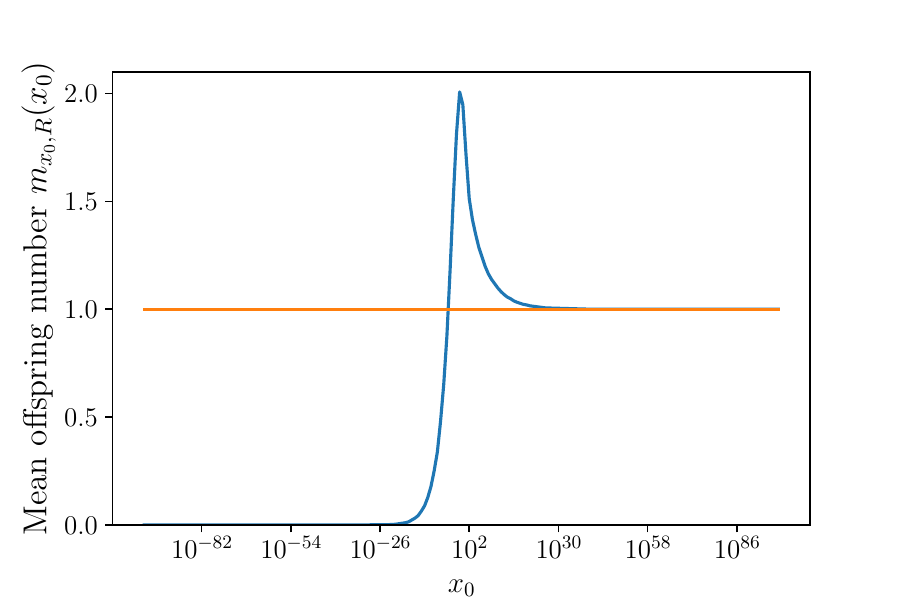}
    \caption{$\beta=-0.1$}
    \label{fig:thirdbeta}
\end{subfigure}
\caption{Monte Carlo estimation of the average number of offspring $m_{x_{0},R}(x_{0})$ for $x_{0} \in [10^{-100},10^{100}] $, plotted on a log-scale, with parameters $C_{\beta}=2$, $C_{\delta}=0.5$, and different values of $\alpha - 1 = -0.25 < \beta < -0.05 = \alpha -1 + \frac{C_{\delta}}{C_{\gamma}-C_{\alpha}}$. For every value of $x_{0}$, we simulated $n =50000$ Monte-Carlo samples to obtain an estimation of $m_{x_{0},R}(x_{0})$. If at some $x_{0}$, the blue line is above the orange line, the numerical simulation suggests that the population process is supercritical for the corresponding value of $x_{0}$, and if it is below, the population process is subcritical.}
\label{fig:beta}
\end{figure}

We illustrate that condition~\eqref{eq:itwocase} is necessary to obtain Assumption~\ref{ass:supercritical} if $\beta > \alpha -1$. To do so, we choose values of $\beta$ outside the $I_{2}$ case, and we observe that there exists some $x_{0}>0$ with $m_{x_{0},R}(x_{0}) \leq 1$. On Figure~\ref{fig:beta}, we represent $m_{x_{0},R}(x_{0})$ depending on $x_{0}$, for values of $\beta$ outside the $I_{2}$ case. We verify that $m_{x_{0},R}(x_{0}) \leq 1$ for small values of $x_{0}$. We also observe that we need to look at smaller and smaller values for $x_{0}$ when $\beta$ goes to $\alpha-1$ to obtain a subcritical regime. Moreover, it seems that over a given value for $x_{0}$, depending on $\beta$, the regime is supercritical, but the mean number of offspring goes to 1 when $x_{0}$ goes to $+ \infty$. There is a value for $x_{0}$ that maximizes the mean number of offspring and it does not seem to depend on $\beta$. 

\subsection{Explosion of $m_{x_{0},R}(x_{0})$ in the $I_{2}$ case}

\begin{figure}[h!]
\centering
\begin{subfigure}{0.4\textwidth}
    \includegraphics[width=\textwidth]{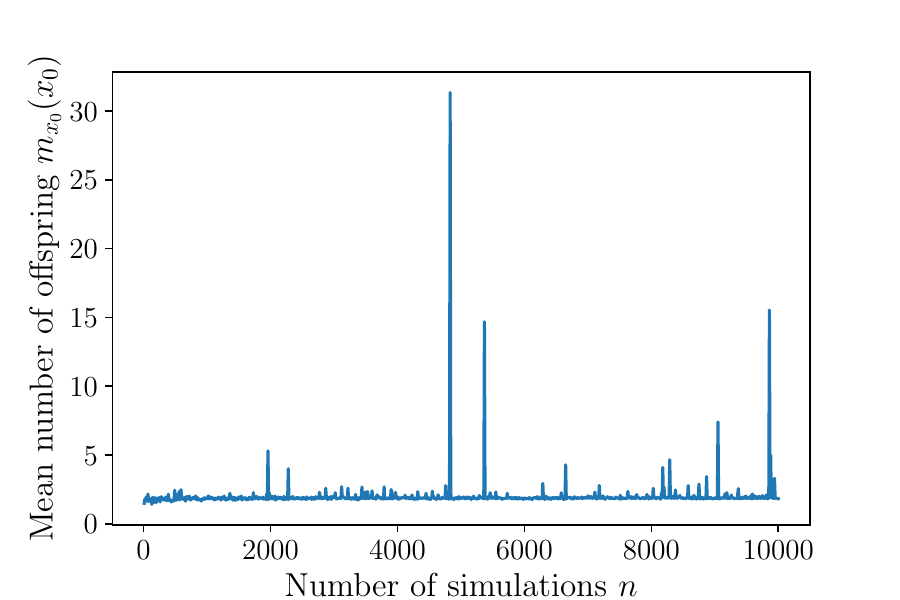}
    \caption{$\beta=0.26$}
    \label{fig:firstmeanos}
\end{subfigure}
\hspace{0.01 cm}
\begin{subfigure}{0.4\textwidth}
     \includegraphics[width=\textwidth]{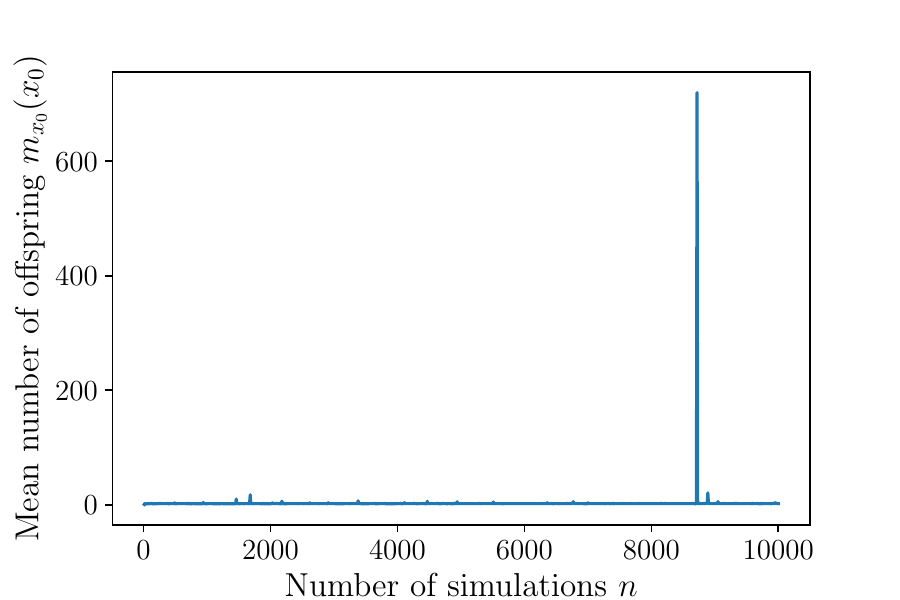}
    \caption{$\beta=1$}
    \label{fig:secondmeanos}
\end{subfigure}
\hspace{0.01 cm}
\begin{subfigure}{0.4\textwidth}
   \includegraphics[width=\textwidth]{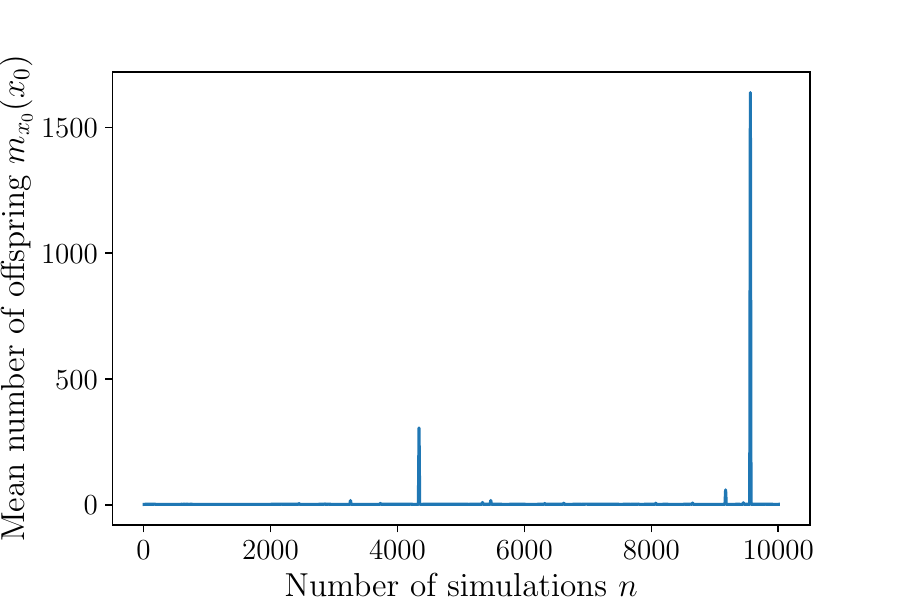}
    \caption{$\beta=3$}
    \label{fig:thirdmeanos}
\end{subfigure}
\caption{Monte-Carlo estimation of the mean number of offspring $m_{x_{0}}(x_{0})$ for $x_{0}=1$, $C_{\beta}=2$, $C_{\delta}=0.5$ and different values of $\beta$. The $x$-axis represents the number $n$ of independent individual trajectories simulated to estimate $m_{x_{0}}(x_{0})$.}
\label{fig:meanos}
\end{figure}

When we continue to increase the value of $\beta$ and look now at the $I_{2}$ case, we obtain an unstable behavior of the Monte Carlo estimate, which leads to Conjecture~\ref{conj:deux}.
\\\\
Let us write $n$ for the number of individual trajectories simulated to estimate $m_{x_{0},R}(x_{0})$. On Figure~\ref{fig:meanos}, we represent  $m_{x_{0},R}(x_{0})$ depending on $n$, for various values of $\beta$ in $I_{2}$. Notice that now, $x_{0}$ is fixed equal to 1. The Monte-Carlo estimation of the mean number of offspring behaves badly. We see higher and higher peaks occuring for random values of $n$ when $\beta$ increases. It becomes harder and harder to estimate correctly this expectation with Monte Carlo estimation, because very rare trajectories like on Figure~\ref{fig:sixth} have a very high contribution. This is a typical heavy-tailed situation, where the mean number of offspring can be very high due to very rare events with a huge amount of births, which highlights a classical weakness of Monte-Carlo estimates. Simulations for others value of $x_{0}$ lead to a similar behavior, so we conjecture that we observe this unstability because $m_{x_{0},R}(x_{0}) = + \infty$ for every $x_{0}$ for these values of $\beta$.

\newpage

\subsection{Sufficient and non sufficient conditions to have supercriticality in the $I_{1}$ case ($\beta=\alpha-1$)}
\label{subsec:proveit}

\begin{figure}[h!]
\centering
\begin{subfigure}{0.4\textwidth}
    \includegraphics[width=\textwidth]{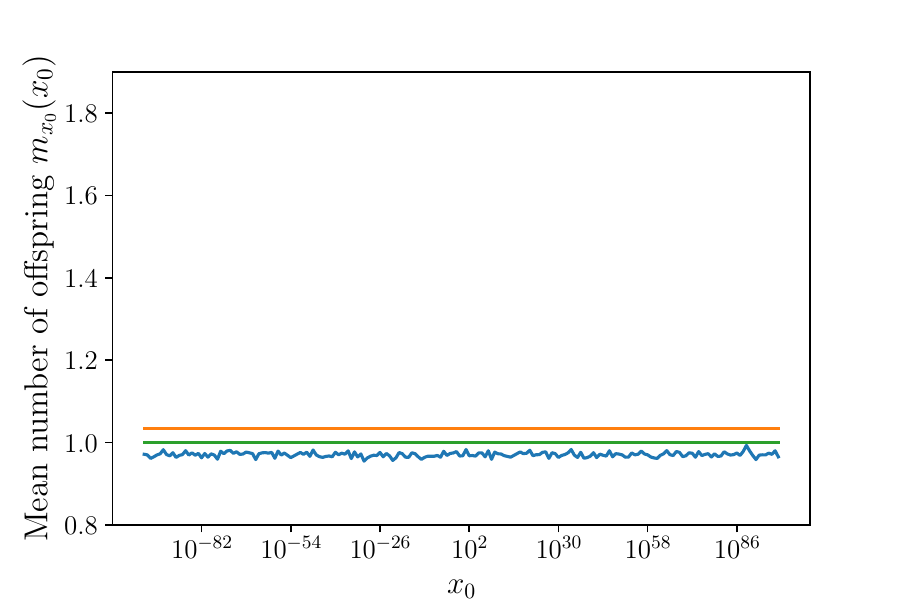} 
    \caption{$C_{\beta}=0.31$}
    \label{fig:firstxo}
\end{subfigure}
\hspace{0.01 cm}
\begin{subfigure}{0.4\textwidth}
   \includegraphics[width=\textwidth]{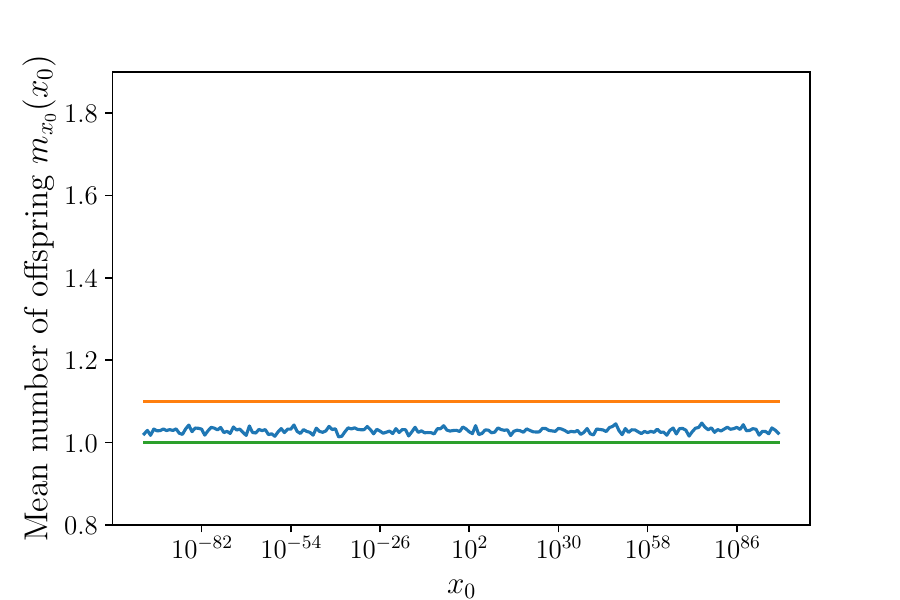}  
    \caption{$C_{\beta}=0.33$}
    \label{fig:secondxos}
\end{subfigure}
\hspace{0.01 cm}
\begin{subfigure}{0.4\textwidth}
   \includegraphics[width=\textwidth]{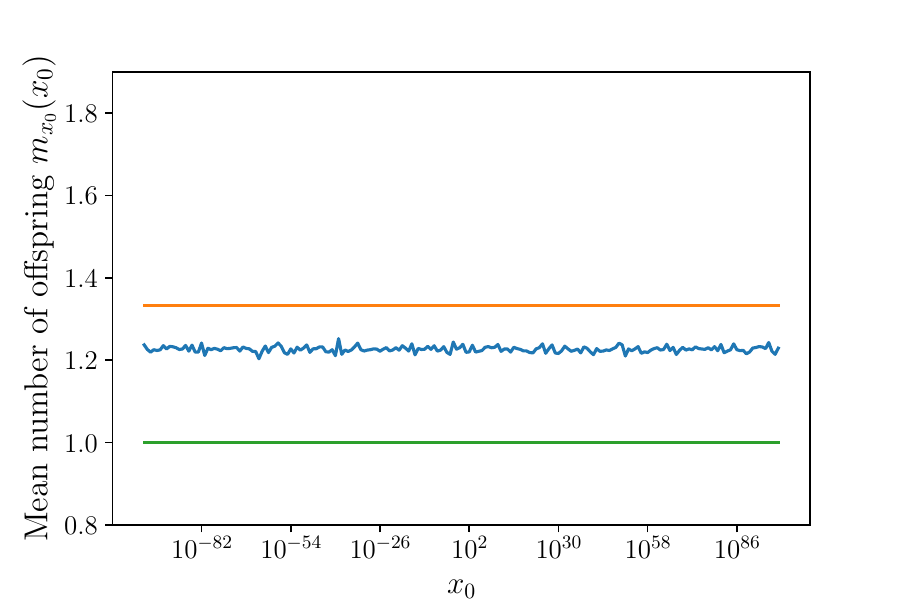}
    \caption{$C_{\beta}=0.4$}
    \label{fig:thirdxo}
\end{subfigure}
\begin{subfigure}{0.4\textwidth}
  \includegraphics[width=\textwidth]{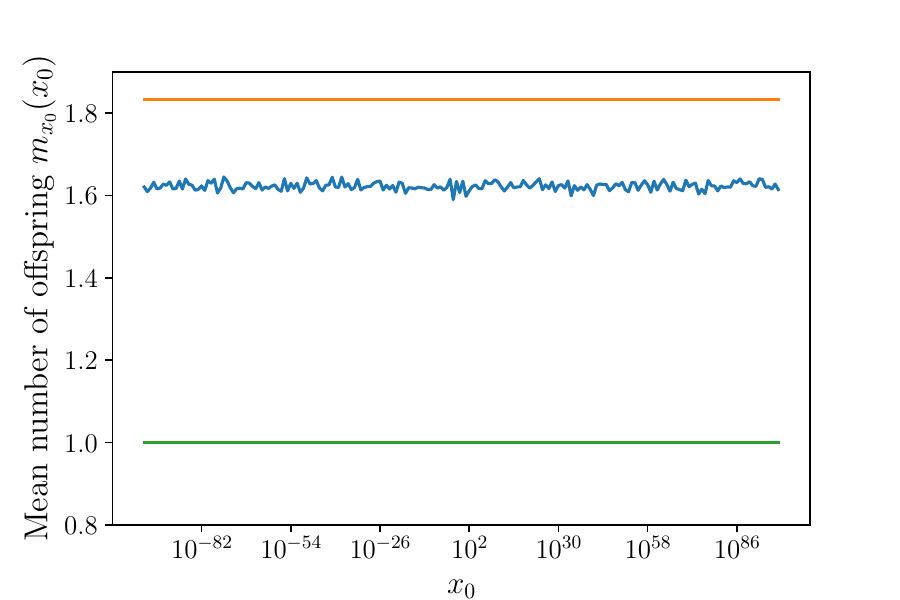}
    \caption{$C_{\beta}=0.55$}
    \label{fig:fourthxo}
\end{subfigure}
\caption{Monte Carlo estimation of the average number of offspring $m_{x_{0},R}(x_{0})$ for $x_{0} \in [10^{-100},10^{100}]$, plotted on a log-scale, and different values of $C_{\beta}$, with parameters $\beta=\alpha-1$, $C_{\delta}=0.3$. For every value of $x_{0}$, we simulated $n :=50000$ Monte-Carlo samples to obtain an estimation of $m_{x_{0},R}(x_{0})$, represented in blue. The orange line is constant equal to $C_{\beta}/C_{\delta}$, and the green line constant equal to 1. If at some $x_{0}$, the blue line is above the green line, the population process is supercritical for the corresponding value of $x_{0}$, and if it is below, the population process is subcritical.}
\label{fig:xobetaegal}
\end{figure}

We illustrate that condition~\eqref{eq:recipro} is sufficient to obtain Assumption~\ref{ass:supercritical}, and that~\eqref{eq:iuncase} is not sufficient (see Section~\ref{subsec:reciprocal} for theoretical proofs).
\\\\
On Figure~\ref{fig:xobetaegal}, we represent the average number of offspring $m_{x_{0},R}(x_{0})$ depending on $x_{0}$. 
We illustrate that $m_{x_{0},R}(x_{0})<C_\beta/C_\delta$ (the blue line is always below the orange one), which is proven in Section~\ref{subsec:reciprocal} (see~\eqref{eq:touslesk}). 
This is a major difference between the $I_{1}$ and the $I_{2}$ case: in the $I_{1}$ case, we prove that $m_{x_{0},R}(x_{0})$ is always finite, and we conjecture that in the $I_{2}$ case, it is always infinite. The scale of the $x$-axis is extremely large comparatively to the scale of the $y$-axis. Hence, we can consider that the variations of the blue curve are essentially due to the intrinsic randomness of the simulations, and that $m_{x_{0},R}(x_{0})$ is independent from $x_{0}$, which is really different from the situation depicted on Figure~\ref{fig:beta}.  This illustrates the independence of our model dynamics along the energy at birth $x_{0}$ when $\beta=\delta=\alpha-1$. It was possible to foresee this result, because in this particular case, the ratio $d/b_{x_{0}}$ is constant on $]x_{0},+\infty[$. A simple time change shows that the birth and death dynamics of the processes $(\xi_{x_{0},R,x_{0}})_{x_{0}>0}$ are the same (up to a time scaling), no matter the value of $x_{0}$. On Figure~\ref{fig:firstxo}, we choose parameters that verify~\eqref{eq:iuncase} (and not~\eqref{eq:recipro}), but Assumption~\ref{ass:supercritical} is not verified as $m_{x_{0},R}(x_{0})<1$, so that \eqref{eq:iuncase} is not sufficient for Assumption~\ref{ass:supercritical}. On Figure~\ref{fig:fourthxo}, we choose parameters that verify~\eqref{eq:recipro} and Assumption~\ref{ass:supercritical} is verified. Notice that on Figures~\ref{fig:secondxos} and \ref{fig:thirdxo}, Assumption~\ref{ass:supercritical} is verified but not~\eqref{eq:recipro}, so that this condition is not necessary.

\subsection{Necessary and sufficient condition for supercriticality in the $I_{1}$ case}

\begin{figure}[h!]
\centering
\begin{subfigure}{0.45\textwidth}
  \includegraphics[width=\textwidth]{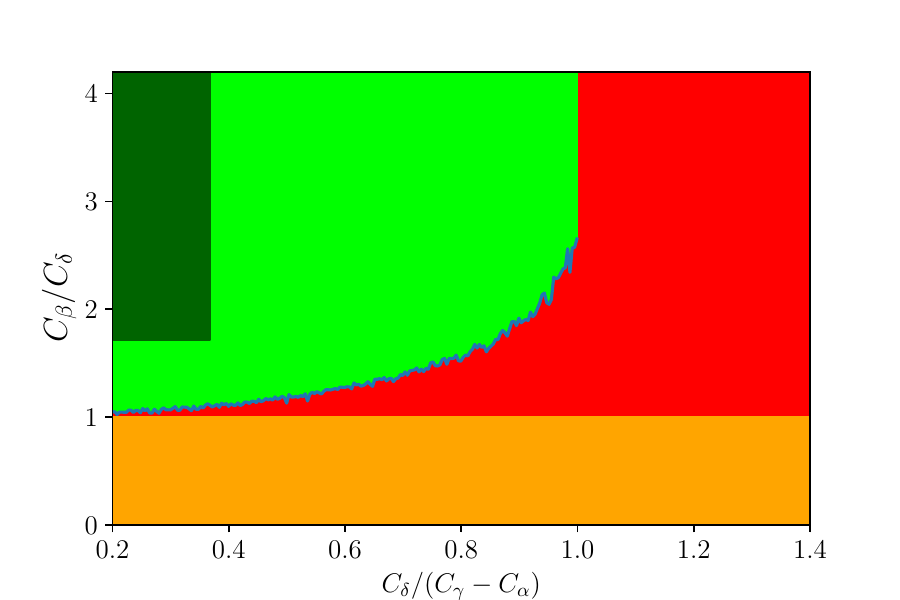}   
    \caption{$C_{\gamma}-C_{\alpha}=0.1$}
    \label{fig:firstdeux}
\end{subfigure}
\hspace{0.01 cm}
\begin{subfigure}{0.45\textwidth}
     \includegraphics[width=\textwidth]{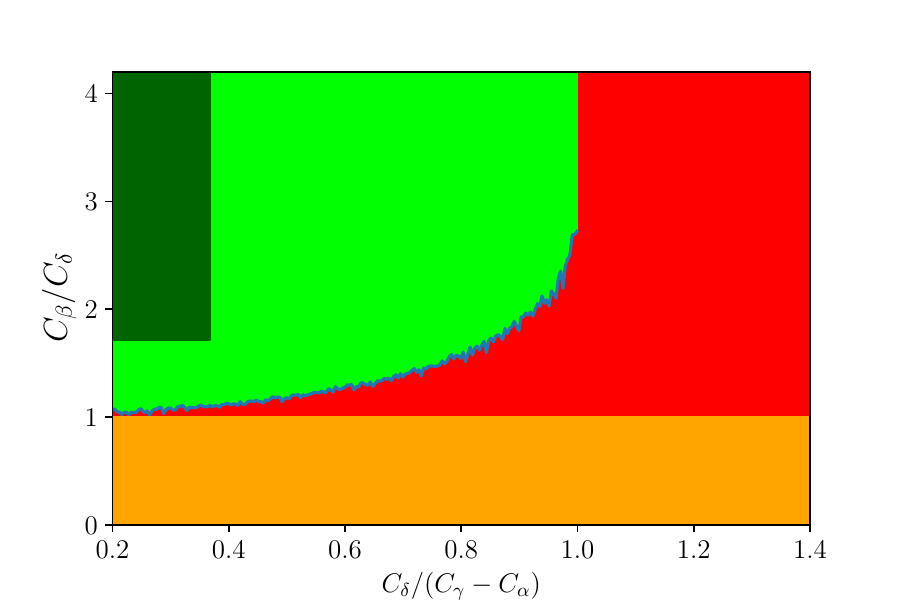}
    \caption{$C_{\gamma}-C_{\alpha}=1$}
    \label{fig:seconddeux}
\end{subfigure}

\caption{Supercriticality of the population process with respect to $C_{\beta}/C_{\delta}$ and $C_{\delta}/(C_{\gamma}-C_{\alpha})$, for $x_{0}=1$, $\beta=\alpha-1$ and two values of $C_\gamma-C_\alpha$. The dark and light green areas correspond to parameters leading to a supercritical regime, whereas the orange and red areas lead to critical or subcritical regimes. The set of parameters for which the necessary condition for supercriticality \eqref{eq:iuncase} is not verified is the orange area, and the set of parameters verifying the sufficient condition for supercriticality \eqref{eq:recipro} are in the dark green area.
For $ 0.2 \leq C_{\delta}/(C_{\gamma}-C_{\alpha}) < 1$, the blue curve describes the value of transition from (sub)critical to supercritical regime. For $C_{\delta}/(C_{\gamma}-C_{\alpha}) \geq 1$, our simulations led only to (sub)critical regime.}
\label{fig:endeuxd}
\end{figure}

We have verified numerically in Section~\ref{subsec:proveit} (and we prove it in Section~\ref{subsec:reciprocal}) that for the case $\beta=\delta=\alpha-1$, the necessary condition~\eqref{eq:iuncase} and the sufficient condition~\eqref{eq:recipro}, respectively highlighted in Theorem~\ref{theo:short} and Proposition~\ref{prop:theoreciprocal}, are not sharp. This is why, for a given value of $C_{\gamma}-C_{\alpha}$, we numerically search for a necessary and sufficient condition on the parameters to have a real equivalence in Theorem~\ref{theo:short} (this gives rise to Conjecture~\ref{conj}). Precisely, we numerically find a function $\Xi : C_{\delta}/(C_{\gamma}-C_{\alpha}) \in ]0,1[ \mapsto \Xi(C_{\delta}/(C_{\gamma}-C_{\alpha})) \in ]1, + \infty[$ such that Assumptions~\ref{hyp:probamort},~\ref{hyp:tempsinf} and~\ref{ass:supercritical} are verified, if and only if $C_{\delta} < C_{\gamma}-C_{\alpha}$ and $\frac{C_{\beta}}{C_{\delta}} > \Xi(C_{\delta}/(C_{\gamma}-C_{\alpha}))$.
\\\\
On Figure~\ref{fig:endeuxd}, we represent in green the region of couples $(C_{\beta}/C_{\delta}, C_{\delta}/(C_{\gamma}-C_{\alpha}))$ for which the regime is supercritical. We observe that if $C_{\delta}/(C_{\gamma}-C_{\alpha}) \geq 1$ then the population process is subcritical. Moreover, the boundary between subcritical and supercritical regimes (see the blue curve) draws a convex increasing curve, leading to the function $\Xi$ in our Conjecture~\ref{conj}, defined on $]0,1[$. Finally, Figures~\ref{fig:firstdeux} and \ref{fig:seconddeux} are strikingly similar (and it is a general observation for other values of $C_{\gamma}-C_{\alpha}$), which is why we formulate Conjecture~\ref{conj} with $\Xi$ depending on the quotient $C_{\delta}/(C_{\gamma}-C_{\alpha})$, rather than on the couple $(C_{\delta},C_{\gamma}-C_{\alpha})$.

\appendix
\appendixpage

\section{Construction of the population process}
\label{app:construc}
\subsection{Construction of the measure-valued process $\mu$}
\label{app:mesure}
We give a formal construction of the population process $\mu$, using the construction of the individual process $\xi$ in Section~\ref{subsec:construc}. In the following, we fix $x_{0}>0$, $R \geq 0$, and we work under the general setting of Section~\ref{subsec:gensett}. We define $\mu_{t}$ for $t \in [0,\overline{\Theta}[$, where $\overline{\Theta}< + \infty$ if there is an accumulation of jump times (and in that case, $\overline{\Theta}$ is the supremum of the jump times), and $\overline{\Theta}= + \infty$ otherwise. Then, we prove in Proposition~\ref{prop:final}, that under Assumptions~\ref{hyp:probamort} and \ref{hyp:tempsinf}, then $\overline{\Theta} = + \infty$ almost surely for every $\mu_{0} \in \mathcal{M}_{P}$. Recall from Section~\ref{subsec:whole} that we adopt the Ulam-Harris-Neveu notation to index individuals in the population. We define
$$ \mathcal{U} := \bigcup_{n \in \mathbb{N}}(\mathbb{N}^{*})^{n+1}.$$
Over time, every individual will have an index of the form $ u := u_{0}...u_{n}$ with some $n \geq 0$, and some positive integers $u_{0},..., u_{n}$, \textit{i.e.} some $u \in \mathcal{U}$. The \textit{generation} of $u$ is $|u| := n$. For every $u := u_{0}...u_{n}$ and $k \geq 1$, we define $uk := u_{0}...u_{n}k$.
\\\\
At time $t=0$, we pick a random variable $\mu_{0} \in \mathcal{M}_{P}$. It means that there exists a random variable $\mathcal{C}_{0} \in \mathbb{N}$ and a random vector $(\xi_{0}^{1},..., \xi_{0}^{\mathcal{C}_{0}}) \in (\mathbb{R}^{*}_{+})^{\mathcal{C}_{0}}$ such that
\begin{center}
        $\mu_{0} = \sum\limits_{u \in V_{0}}{\delta_{\xi_{0}^{u}}},$
\end{center}
with $V_{0} := \{1,...,\mathcal{C}_{0}\}$. Then, we define a family of processes $(\xi^{u})_{u \in \mathcal{U}}$, independent conditionally to $\mu_{0}$, such that for every $u \in \mathcal{U}$, $\xi^{u}$ is distributed as the process $\xi$ of Section~\ref{subsec:construc}, started from $\xi_{0}^{u}$ if $u \in V_{0}$, and started from $x_{0}$ if $u \in \mathcal{U} \setminus V_{0}$ (which corresponds to the fixed amount of energy transferred to offspring). We use the construction of Section~\ref{subsec:construc} with independent exponential random variables $(E_{i}^{u})_{u \in \mathcal{U},i \geq 1}$ and independent uniform random variables $(U_{i}^{u})_{u \in \mathcal{U}, i \geq 1}$. We assume that these sequences are independent from each other, and both independent from $\mu_{0}$. For $u \in \mathcal{U}$, we write $(J_{k}^{u})_{k \geq 1}$ for the jump times of $\xi^{u}$, the time of death is $T_{d}^{u}$ and the time when the process possibly reaches $\flat$ is $T_{0}^{u} \wedge T_{\infty}^{u}$. Also, for $u \in \mathcal{U}$ and $\tau \geq 0$, we define the shifted process $\xi^{u, \tau}$ on $[\tau,+ \infty[$ with $\xi^{u,\tau}_{t} := \xi^{u}_{t-\tau}$ for every $t \geq \tau$.  Over time, we will describe $\mu_{t}$ as a point measure given by 
\begin{align}
\mu_{t} := \sum\limits_{u \in V_{t}}{\delta_{\xi_{t}^{u,\tau_{u}}}},
\label{muinitial}
\end{align}
which means that at time $t$, alive individuals have indices $u$ in some $V_{t} \subseteq \mathcal{U}$. Recall that by definition of alive individuals, $V_{t}$ does not contain individuals absorbed in $\{\partial, \flat \}$. Individual $u$ is born at time $\tau_{u}$ and follows the process $\xi^{u,\tau_{u}}$. In particular, we set $\tau_{u}=0$ for $u \in V_{0}$.
\\\\
At time $t=0$, the initial condition of the population process is 
\begin{center}
        $\mu_{0} := \sum\limits_{u \in V_{0}}{\delta_{\xi_{0}^{u,0}}}.$
\end{center}
We now define the sequence $(\Theta_{n})_{n \in \mathbb{N}}$ of successive times of jump of the population process. First, we set $\Theta_{0} := 0$, and then suppose that our process is described until some time $\Theta_{n}$ with $n \geq 0$. If $\Theta_{n}= + \infty$, the process is already well-defined for every $t \geq 0$ and we set $\Theta_{n+1}= + \infty$. Else, at time $\Theta_{n}$, the process is of the form~\eqref{muinitial} with some finite $V_{\Theta_{n}} \subseteq \mathcal{U}$. The next time of jump for the whole population is then 
$$ \Theta_{n+1} := \inf_{u \in V_{\Theta_{n}}} \{(J_{k}^{u}\wedge T_{d}^{u} \wedge T_{0}^{u} \wedge T_{\infty}^{u})  + \tau_{u}, \quad  k \in \mathbb{N}^{*}, J_{k}^{u} + \tau_{u} \geq \Theta_{n} \},$$  
with the convention $\inf(\varnothing) = + \infty$. For $s \in [\Theta_{n},\Theta_{n+1}[$, we set $$\mu_{s} = \sum\limits_{u \in V_{\Theta_{n}}}{\delta_{\xi_{s}^{u,\tau_{u}}}}.$$
If $\Theta_{n+1} = + \infty$, we have defined the process $\mu_{t}$ for every $t \geq 0$. Else, almost surely, the infimum in the definition of $\Theta_{n+1}$ is reached at a single element $u \in V_{\Theta_{n}}$. First, if $\Theta_{n+1}$ is of the form $J_{k}^{u} + \tau_{u}$ for some $k \in \mathbb{N}^{*}, u \in V_{\Theta_{n}}$ and $J^{u}_{k} \neq T_{d}^{u}$, it means that this jump is the $k$-th birth for the individual $u$. In that case, we set $$\mu_{\Theta_{n+1}} = \mu_{\Theta_{n+1}-} + \delta_{\xi_{\Theta_{n+1}}^{uk,\Theta_{n+1}}}$$
and $\tau_{uk} := \Theta_{n+1}$, $V_{\Theta_{n+1}} = V_{\Theta_{n}} \cup \{ uk \}$.
\\
Then, if $\Theta_{n+1}$ is of the form $T_{d}^{u} + \tau_{u}$ for some $u \in V_{\Theta_{n}}$, it means that this jump is the death of individual $u$. In that case, we set 
$$\mu_{\Theta_{n+1}} = \mu_{\Theta_{n+1}-} - \delta_{\xi_{\Theta_{n+1}-}^{u,\tau_{u}}}$$
and $V_{\Theta_{n+1}} = V_{\Theta_{n}} \setminus \{ u \}$.
\\
Finally, if $\Theta_{n+1}$ is of the form $(T_{0}^{u} \wedge T_{\infty}^{u}) + \tau_{u}$, it means that at this jump, the energy of the individual $u$ reaches $\{ 0, + \infty \}$. In that case, we also remove the individual from $V_{\Theta_{n+1}}$:
$$\mu_{\Theta_{n+1}} = \mu_{\Theta_{n+1}-} - \delta_{\xi_{\Theta_{n+1}-}^{u,\tau_{u}}}$$
and $V_{\Theta_{n+1}} = V_{\Theta_{n}} \setminus \{ u \}$. Finally, we define $\overline{\Theta} := \sup\limits_{n \in \mathbb{N}} \Theta_{n}$. If for all $t \in [0, \overline{\Theta}[$, $u \notin V_{t}$, we set $\tau_{u} := + \infty$. 
\\\\
Remark that, until now, $\mu$ is well-defined only on $[0,\overline{\Theta}[$, with $\overline{\Theta}$ possibly finite or infinite. We work now under Assumptions~\ref{hyp:probamort} and \ref{hyp:tempsinf}. Under these assumptions, any individual process $\xi^{u}$ is almost surely biologically relevant by Theorem~\ref{theo:cns}. Thus, we can construct the individual process $\xi^{u}$ as in Section~\ref{subsec:construpta}, using a Poisson point measure $\mathcal{N}_{u}$ and an exponential random variable $E_{u}$ for the death. It is equivalent in law to define another independent Poisson point measure $\mathcal{N}'_{u}$ with intensity $\mathrm{d}s \times \mathrm{d}h$, and to work with $\xi^{u}$ defined for $t \geq 0$ by
\begin{align*} 
\xi^{u}_{t} := &  A_{\xi^{u}_{0}}(t) + \displaystyle{\int_{0}^{t}\int_{\mathbb{R}^{+}}}\mathbb{1}_{\{ h \leq b_{x_{0}}(\xi^{u}_{s-})\}} \bigg(A_{\xi^{u}_{s-}-x_{0}}(t-s) - A_{\xi^{u}_{s-}}(t-s)\bigg) \mathcal{N}_{u}(\mathrm{d}s, \mathrm{d}h) \\
& + \displaystyle{\int_{0}^{t}\int_{\mathbb{R}^{+}}}\mathbb{1}_{\{ h \leq d(\xi^{u}_{s-})\}} \bigg(\partial - A_{\xi^{u}_{s-}}(t-s)\bigg) \mathcal{N}'_{u}(\mathrm{d}s, \mathrm{d}h),
\end{align*}
with the convention $b_{x_{0}}(\partial) = d(\partial) = 0$. Then, from the construction of $\mu$ and as the individual processes $(\xi^{u}_{t})_{t \geq 0}$ are independent from each other, we have for $t \in [0,\overline{\Theta}[$,
\begin{multline}
    \hspace{-0.3cm} \mu_{t} =  \sum\limits_{u \in V_{0}} \delta_{A_{\xi^{u}_{0}}(t)} + \sum\limits_{u \in \mathcal{U}}\bigg[ \displaystyle{\int_{0}^{t}\int_{\mathbb{R}^{+}}} \mathbb{1}_{\{ u \in V_{s-} \}} \mathbb{1}_{\{ h \leq b_{x_{0}}(\xi^{u,\tau_{u}}_{s-})\}} \bigg(\delta_{A_{x_{0}}(t-s)} +\delta_{A_{\xi^{u,\tau_{u}}_{s-}-x_{0}}(t-s)} - \delta_{A_{\xi^{u,\tau_{u}}_{s-}}(t-s)} \bigg) \mathcal{N}_{u}(\mathrm{d}s, \mathrm{d}h)\\ - \displaystyle{\int_{0}^{t}\int_{\mathbb{R}^{+}}} \mathbb{1}_{\{ u \in V_{s-} \}}\mathbb{1}_{\{ h \leq d(\xi^{u,\tau_{u}}_{s-})\}} \delta_{A_{\xi^{u,\tau_{u}}_{s-}}(t-s)} \mathcal{N}'_{u}(\mathrm{d}s, \mathrm{d}h) \bigg].
    \label{eq:muencoremu}
\end{multline} 
Notice that this definition is valid almost surely on $[0, \overline{\Theta}[$, because individual trajectories are biologically relevant. This is now very similar to the classical framework for the construction of measure-valued processes \cite{four_04, coranicofab, marguet2016}, but the birth rate $b_{x_{0}}$ is not necessarily bounded on $\mathbb{R}^{*}_{+}$.
\begin{prop}
Under the general setting of Section~\ref{subsec:gensett}, under Assumptions~\ref{hyp:probamort} and \ref{hyp:tempsinf}, for every random variable $\mu_{0} \in \mathcal{M}_{P}$, we have
$$ \forall x_{0}>0, \forall R \geq 0, \quad \mathbb{Q}_{\mu_{0},x_{0},R}(\overline{\Theta}=+ \infty)=1 .$$
\label{prop:final}
\end{prop}
\begin{proof}
If we write $\Lambda$ for the law of $\mu_{0}$, then
$$ \mathbb{Q}_{\mu_{0},x_{0},R}(\overline{\Theta} < + \infty) = \displaystyle{\int_{\mathcal{M}_{P}} \mathbb{Q}_{\mu,x_{0},R}(\overline{\Theta} < + \infty) \mathrm{d}\Lambda(\mu)}. $$
Thus, it is sufficient to show that for any deterministic initial condition $\mu_{0} := \sum\limits_{u \in V_{0}} \xi^{u}_{0} \in \mathcal{M}_{P}$, for every $T>0$, $\mathbb{Q}_{\mu_{0},x_{0},R}(\overline{\Theta} < T)=0$. We fix such $\mu_{0}$ and $T$ in the following, then the initial energies $(\xi^{u}_{0})_{u \in V_{0}}$ are in a compact depending on $\mu_{0}$, and $\langle \mu_{0},1 \rangle := \int_{\mathbb{R}^{*}_{+}} \mu_{0}(\mathrm{d}x)= \mathrm{Card}(V_{0})$ is finite. For $M >0$, we define the stopping time $\tau_{M} := \inf \{ t \geq 0, \langle \mu_{t},1 \rangle \geq M \}$. We have from \eqref{eq:muencoremu} that for $t \leq \tau_{M} \wedge T$,
\begin{align}
    \langle \mu_{t}, 1 \rangle \leq  \langle \mu_{0}, 1 \rangle + \sum\limits_{u \in \mathcal{U}} \displaystyle{\int_{0}^{t}\int_{\mathbb{R}^{+}}} \mathbb{1}_{\{ u \in V_{s-} \}} \mathbb{1}_{\{ h \leq b_{x_{0}}(\xi^{u,\tau_{u}}_{s-})\}} \mathcal{N}_{u}(\mathrm{d}s, \mathrm{d}h).
 \label{eq:controlemu}
\end{align}
The initial condition $\mu_{0}$ is compactly supported, so individual trajectories are almost surely in a compact $\mathfrak{C}_{\mu_{0},T}$ until time $t \leq T$. Also, $b_{x_{0}}$ is bounded by a constant $\overline{b}_{\mu_{0},T}$ on $\mathfrak{C}_{\mu_{0},T}$, and $\mathrm{Card}(V_{s}) \leq M$ until time $t \leq \tau_{M} \wedge T$. Hence, the integrand in the previous upper bound is almost surely bounded, and classical arguments using \eqref{eq:controlemu} allow us to prove first that $\mathbb{E}\left( \sup_{t \leq \tau_{M} \wedge T} \langle \mu_{t},1 \rangle\right) < C'_{T} $, where $C'_{T}< + \infty$ does not depend on $M$ (see the proof of Corollary 4.3. in \cite{campillo2015weak}). This entails that almost surely, $\tau_{M}$ goes to $+ \infty$ when $M$ goes to $+ \infty$.
\\\\
Finally, let us suppose by contradiction that $\overline{\Theta} < T$ with positive probability and work on this event. In particular, this accumulation of jump times means that there exists a random subsequence of $(\Theta_{n})_{n \in \mathbb{N}}$ corresponding to an infinite number of birth times, and denoted as $(\mathcal{B}_{n})_{n \in \mathbb{N}}$ in the following. Also, we fix $M >0$ such that $\tau_{M} > T$, and remark that if 
\begin{align}
\exists n \in \mathbb{N}, \tau_{M} \leq \Theta_{n},
\label{eq:coralie}
\end{align}
then there is a contradiction because $\overline{\Theta} < T < \tau_{M} \leq \Theta_{n} \leq \overline{\Theta}$. Hence, to conclude, it suffices to show that on the event $\{\overline{\Theta} < T\}$, \eqref{eq:coralie} is  verified almost surely. We assume by contradiction that $\overline{\Theta} < T$, but \eqref{eq:coralie} is not verified with positive probability. Then, we would have $\langle \mu_{t}, 1 \rangle < M$ for every $t \in [0, \overline{\Theta}[$. In that case, the subsequence $(\mathcal{B}_{n})_{n \in \mathbb{N}}$ can be constructed in law as a subsequence of a Poisson point process with intensity $M\overline{b}_{\mu_{0},T} $. The only accumulation point of $(\mathcal{B}_{n})_{n \in \mathbb{N}}$ would then almost surely be $+ \infty$, which contradicts the fact that for all $n \in \mathbb{N}$, $ 0 \leq \mathcal{B}_{n} \leq \overline{\Theta} < T$ and concludes the proof.
\end{proof}

Remark that in the inductive definition of $\mu$, we assumed that the processes $(\xi^{u})_{u \in \mathcal{U}}$ are independent. This is what we call a \textit{branching property} for the population process, and will be illustrated in the next section.

\subsection{Proof of Proposition~\ref{prop:jesuisungalton} and Proposition~\ref{prop:galtonwatson} in Section~\ref{subsec:generation}}
\label{appendix:embedding}
First, we give the proof of Proposition~\ref{prop:jesuisungalton}, based on Lemma~\ref{petitlemme} below. Recall that for $n \in \mathbb{N}$, $$G_{n} := \{ u \in \mathcal{U}, |u|=n, \exists t \geq 0, u \in V_{t} \} $$
contains all the individuals of the $n$-th generation and $ \Upsilon_{n} := \mathrm{Card}(G_{n}).$
\begin{lemme}
Under Assumption~\ref{ass:pasinfini}, 
$$\forall \mu_{0} \in \mathcal{M}_{P}, \forall x_{0}>0, \forall R \geq 0, \quad \mathbb{Q}_{\mu_{0},x_{0},R}(\forall n \in \mathbb{N}, \Upsilon_{n} < + \infty)=1.$$
\label{petitlemme}
\end{lemme}
\begin{proof}
We fix $\mu_{0},x_{0},R$ and use induction on $n$. For the 0-generation, because $\mathcal{C}_{0} = \langle \mu_{0},1 \rangle$ is a random variable taking values in $\mathbb{N}$, we almost surely have 
$$ \mathrm{Card}(G_{0}) = \mathrm{Card}(V_{0}) < + \infty. $$
Now, let us suppose that our result holds for some $n \geq 0$, that is almost surely,
\begin{align}
\mathrm{Card}(G_{n}) < + \infty.
\label{eq:supgn}
\end{align}
By Assumption~\ref{ass:pasinfini}, we have almost surely
\begin{align}
\forall u \in G_{n}, \quad N^{u} <+ \infty,
\label{eq:nfni}
\end{align}   
where $N^{u}$ is the number of direct offspring of individual $u$ defined in Section~\ref{subsec:construc}. Combining~\eqref{eq:supgn} and~\eqref{eq:nfni} entails that almost surely,
$$ \mathrm{Card}(G_{n+1}) = \sum\limits_{u \in G_{n}} N^{u} < + \infty.$$
\end{proof}
By Lemma~\ref{petitlemme}, $\Upsilon_{1}$ is almost surely finite, and for all $n \geq 0$,
$$ \Upsilon_{n+1} := \sum_{u \in G_{n}}N^{u}.  $$
From the construction of Appendix~\ref{app:mesure}, the state at birth of every offspring after generation $G_{1}$ is the same, starting from $x_{0}$, and the resource remains constant equal to $R$, so the law of every $N^{u}$ is $\nu_{x_{0},R,x_{0}}$. Also, all these individual trajectories are independent, so are the $N^{u}$. This is exactly the setting defining a Galton-Watson process with offspring distribution $\nu_{x_{0},R,x_{0}}$, which ends the proof of Proposition~\ref{prop:jesuisungalton}. Now, we give the proof of Proposition~\ref{prop:galtonwatson}.
\begin{lemme}
\label{lemm:generations}
Under Assumptions~\ref{hyp:probamort} and \ref{hyp:tempsinf}, for all $n \in \mathbb{N}$, there exists almost surely a random time $\sigma_{n} \in [0, + \infty]$ such that $\sigma_{n}=+ \infty$, if and only if $G_{n} = \varnothing$; and if $\sigma_{n} < + \infty$, then
\begin{align*}
(t \geq \sigma_{n}) \Leftrightarrow (\forall s \geq t, G_{n} \cap V_{s} = \varnothing).
\label{eq:taunnnn}
\end{align*} 
In addition, $\overline{\sigma}_{n} := \sup\limits_{0 \leq m \leq n} \sigma_{m} \xrightarrow[n \rightarrow + \infty]{} + \infty$.
\end{lemme}
\begin{proof}
First, under Assumptions~\ref{hyp:probamort} and \ref{hyp:tempsinf}, by Proposition~\ref{prop:final}, we work on the almost sure event $\{ \overline{\Theta} = + \infty \}$. By Corollary~\ref{corr:swagosss} in Section~\ref{subsec:teun}, Assumption~\ref{hyp:tempsinf} implies Assumption~\ref{ass:pasinfini}, so we use Lemma~\ref{petitlemme} and work on the event $\{\forall n \in \mathbb{N}, \Upsilon_{n}<+\infty \}$. For $n \in \mathbb{N}$, we define $\sigma_{n} := \mathrm{sup}\{\tau_{u} + T_{d}^{u}, u \in G_{n}\}$, with the convention $\sup(\varnothing)=+ \infty$.  Under Assumptions~\ref{hyp:probamort} and \ref{hyp:tempsinf}, thanks to Theorem~\ref{theo:cns}, for every $t \geq 0$, $u \in V_{t}$, $\tau_{u} + T_{d}^{u} < + \infty$ almost surely. Also, for every $n \in \mathbb{N}$, $\Upsilon_{n} = \mathrm{Card}(G_{n})<+\infty$. Hence, $\sigma_{n} = + \infty$, if and only if $G_{n}= \varnothing$.  Also, it is clear by definition that if $\sigma_{n}< + \infty$, then
\[(t \geq \sigma_{n}) \Leftrightarrow (\forall s \geq t, G_{n} \cap V_{s} = \varnothing).\]
The sequence $(\overline{\sigma}_{n})_{n \in \mathbb{N}}$ is increasing, so converges to some $\tilde{\Theta} \leq + \infty$. If by contradiction, $\tilde{\Theta} < + \infty$, we would have
$$ \forall t \geq \tilde{\Theta}, \forall n \in \mathbb{N}, G_{n} \cap V_{t} = \varnothing, $$
so no random jumps can occur after time $\tilde{\Theta}$, which contradicts the fact that $\overline{\Theta}= + \infty$.
\end{proof}
We take now $x_{0} >0$, and we work under Assumptions~\ref{hyp:probamort} and \ref{hyp:tempsinf}, so by Proposition~\ref{prop:final}, we have $\overline{\Theta}=+ \infty$ almost surely and we work on this event. By Corollary~\ref{corr:swagosss} in Section~\ref{subsec:teun}, Proposition~\ref{prop:jesuisungalton} and \cite[Theorem 1, p.7]{athreya:ney}, for all $R \geq 0$, for all $\mu_{0} \in \mathcal{M}_{P,x_{0},R}$ (so $\mathfrak{B} = \{G_{1} \neq \varnothing \}$ occurs with positive probability), we have $$ m_{x_{0},R}(x_{0}) >1 \Leftrightarrow \mathbb{Q}_{\mu_{0},x_{0},R}( \forall n \in \mathbb{N}, G_{n} \neq \varnothing )>0.$$
Thus, to prove Proposition~\ref{prop:galtonwatson}, it suffices to show that on the event $\{ \overline{\Theta} = + \infty\}$, $$ \{ \forall t \geq 0, V_{t} \neq \varnothing \} = \{ \forall n \in \mathbb{N}, G_{n} \neq \varnothing \}.$$
We begin with the inclusion $\subseteq$. If there exists $n \in \mathbb{N}$ such that $G_{n} = \varnothing$, then we can define $m := \min \{ n \in \mathbb{N}, G_{n} = \varnothing \} $. Remark that $\mu_{0} \in \mathcal{M}_{P,x_{0},R}$ implies that $\mu_{0} \neq 0$, so $m \geq 1$. Then, by Lemma~\ref{lemm:generations}, $\sigma_{m-1}<+ \infty$ and $\sigma_{n}=+ \infty$ for every $n \geq m$. Hence, $V_{t} = \varnothing$ for every $t \in [\sigma_{m-1},+ \infty[$. Then for the inclusion $\supseteq$, suppose that: $\forall n \in \mathbb{N}, G_{n} \neq \varnothing$. With Lemma~\ref{lemm:generations}, we have $\overline{\sigma}_{n} \xrightarrow[n \rightarrow + \infty]{} + \infty$. Let $t \geq 0$, then there exists $n \in \mathbb{N}$ with $\overline{\sigma}_{n} >t$. By definition, there exists $s \in [t, \overline{\sigma}_{n}[$ and $m \leq n$ such that $G_{m} \cap V_{s} \neq \varnothing$, so necessarily $V_{t} \neq \varnothing$. This reasoning means exactly
$$ \{ \forall t \geq 0, V_{t} \neq \varnothing \} \supseteq \{ \forall n \in \mathbb{N}, G_{n} \neq \varnothing \}.$$
\section{Additional content for Section~\ref{sec:pta}}
\label{appendix:pta}
\subsection{Proof of Proposition~\ref{prop:martingaleprobx}}
\label{app:propsept}
Let $G$ be a measurable function. We have, as an immediate consequence of the telescopic expression of $X$ given by \eqref{eq:constructionofx},
\begin{align} 
G(X_{t}) =  G(A_{\xi_{0}}(t))
+ \displaystyle{\int_{0}^{t}\int_{\mathbb{R}^{+}}}\mathbb{1}_{\{ h \leq b_{x_{0}}(X_{s-})\}} \bigg(G(A_{X_{s-}-x_{0}}(t-s)) - G(A_{X_{s-}}(t-s))\bigg) \mathcal{N}(\mathrm{d}s, \mathrm{d}h).
\label{eq:etoile}
\end{align}
For any fixed $t \geq 0$, we apply~\eqref{eq:etoile} to $x \mapsto F(t,x)$:
\begin{align}
F(t,X_{t}) = & \hspace{0.1 cm} F(t,A_{\xi_{0}}(t)) + \displaystyle{\int_{0}^{t}\int_{\mathbb{R}^{+}}}\mathbb{1}_{\{ h \leq b_{x_{0}}(X_{s-})\}} \bigg(F(t,A_{X_{s-}-x_{0}}(t-s)) - F(t,A_{X_{s-}}(t-s))\bigg) \mathcal{N}(\mathrm{d}s, \mathrm{d}h).
\label{eq:tomate}
\end{align} 
Take $z > 0$ and $0 \leq s \leq t$. Since $F$ and $u \in [s,t] \mapsto A_{z}(u-s)$ are respectively $\mathcal{C}^{1,1}$ and $\mathcal{C}^{1}$, the fundamental theorem of calculus gives:
\begin{align*}
F(t,A_{z}(t-s)) &  = F(s,z) + \displaystyle{\int_{s}^{t} \bigg(\partial_{1} F(u,A_{z}(u-s)) + \partial_{2} F(u,A_{z}(u-s))g(A_{z}(u-s),R) \bigg) \mathrm{d}u }. \\
\end{align*}
We get
\begin{align*}
F(t,X_{t}) = & \hspace{0.1 cm} F(0,\xi_{0}) + F_{0} + \displaystyle{\int_{0}^{t}\int_{\mathbb{R}^{+}}}\mathbb{1}_{\{ h \leq b_{x_{0}}(X_{s-})\}} \bigg(F(s,X_{s-}-x_{0})- F(s,X_{s-})\bigg) \mathcal{N}(\mathrm{d}s, \mathrm{d}h) + F_{1},
\end{align*}
with 
\begin{align*}
F_{0} := & \displaystyle{\int_{0}^{t} \bigg(\partial_{1} F(u,A_{\xi_{0}}(u)) + \partial_{2} F(u,A_{\xi_{0}}(u))g(A_{\xi_{0}}(u),R) \bigg) \mathrm{d}u }, \\
F_{1} := & \displaystyle{\int_{0}^{t}\int_{\mathbb{R}^{+}}}\mathbb{1}_{\{ h \leq b_{x_{0}}(X_{s-})\}} \displaystyle{\int_{s}^{t}} \bigg[ \bigg(\partial_{1} F(u,A_{X_{s-}-x_{0}}(u-s)) + \partial_{2} F(u,A_{X_{s-}-x_{0}}(u-s))g(A_{X_{s-}-x_{0}}(u-s),R) \bigg) \\
& \hspace{2 cm} - \bigg( \partial_{1} F(u,A_{X_{s-}}(u-s)) + \partial_{2} F(u,A_{X_{s-}}(u-s))g(A_{X_{s-}}(u-s),R) \bigg) \bigg]  \mathrm{d}u \hspace{0.1 cm} \mathcal{N}(\mathrm{d}s, \mathrm{d}h).
\end{align*}
Recall that the first two integrals in the definition of $F_{1}$ refer to the Poisson point measure $\mathcal{N}$, so are only a formal writing of an almost surely finite random sum, meaning that we can swap the integrals and obtain
\begin{align*}
F_{1} = & \displaystyle{\int_{0}^{t}\int_{0}^{u}\int_{\mathbb{R}^{+}}}\mathbb{1}_{\{ h \leq b_{x_{0}}(X_{s-})\}} \bigg[ \bigg(\partial_{1} F(u,A_{X_{s-}-x_{0}}(u-s)) + \partial_{2} F(u,A_{X_{s-}-x_{0}}(u-s))g(A_{X_{s-}-x_{0}}(u-s),R) \bigg) \\
& \hspace{2 cm} - \bigg( \partial_{1} F(u,A_{X_{s-}}(u-s)) + \partial_{2} F(u,A_{X_{s-}}(u-s))g(A_{X_{s-}}(u-s),R) \bigg) \bigg]  \mathcal{N}(\mathrm{d}s, \mathrm{d}h) \hspace{0.1 cm} \mathrm{d}u.
\end{align*}
If we look back at~\eqref{eq:tomate}, we realize that
\[ F_{0}+F_{1} = \displaystyle{\int_{0}^{t} \bigg( \partial_{1} F(u,X_{u}) + \partial_{2} F(u,X_{u}) g(X_{u},R) \bigg) \mathrm{d}u}, \]
which, associated to the decomposition $\mathcal{N}(\mathrm{d}s,\mathrm{d}h) = \mathcal{N}_{C}(\mathrm{d}s,\mathrm{d}h) + \mathrm{d}s\mathrm{d}h$, leads to the expected expression of $F(t,X_{t})$.
\\\\
At this step, $M_{F}$ is only a local martingale. Let $T>0$, and $t \in [0,T]$. From a fixed initial condition $\xi_{0}>0$, as $A_{\xi_{0}}(.)$ is defined on $\mathbb{R}^{+}$ and $\mathcal{C}^{1}$, $X_{t} \leq A_{\xi_{0}}(t)$ is bounded by some constant $C_{T}$ for $t \in [0,T]$. Also, $b_{x_{0}} \equiv 0$ on $]0,x_{0}]$ and is continuous on $]x_{0},C_{T}]$, and $F$ is $\mathcal{C}^{1,1}$. Hence, $b_{x_{0}}(X_{t})$ and quantities of the form $F(t,X_{t})$ are also bounded for $t \in [0,T]$. Finally, $M_{F,t}$ is $L^{2}$ for all $t \in [0,T]$, so $M_{F}$ is a true martingale.
To obtain the expression of $\langle M_{F} \rangle_{t}$, we give two different semi-martingale decompositions of $F^{2}(t,X_{t})$. This is the technique used to prove Proposition 3.4. in \cite{four_04}. We also refer to \cite{rudiger2006ito}. First, we apply the previous result to $F^{2}$ and obtain:
\begin{align*} 
F^{2}(t,X_{t}) = & \hspace{0.1 cm} F^{2}(0,\xi_{0}) + \displaystyle{\int_{0}^{t}2F(s,X_{s})\bigg(\partial_{1} F(s,X_{s}) + \partial_{2} F(s,X_{s})g(X_{s},R) \bigg) \mathrm{d}s} \\
& + \displaystyle{\int_{0}^{t}} b_{x_{0}}(X_{s}) \bigg(F^{2}(s,X_{s}-x_{0})- F^{2}(s,X_{s})\bigg) \mathrm{d}s + M_{F^{2},t}.
\end{align*}
Then, we apply Itô's formula to $F(t,X_{t})$ and the function $x \mapsto x^{2}$ to get:
\begin{align*}
F^{2}(t,X_{t}) = &  F^{2}(0,\xi_{0})  + \displaystyle{\int_{0}^{t}2F(s,X_{s})\bigg(\partial_{1} F(s,X_{s}) +\partial_{2} F(s,X_{s})g(X_{s},R) \bigg) \mathrm{d}s} \\
& + \displaystyle{\int_{0}^{t}} 2F(s,X_{s}) b_{x_{0}}(X_{s})\bigg(F(s,X_{s}-x_{0})- F(s,X_{s})\bigg) \mathrm{d}s \\
& + \displaystyle{\int_{0}^{t}} 2F(s,X_{s}) \mathrm{d}M_{F,s} + \langle M \rangle_{F,t}.
\end{align*}
By the Doob-Meyer theorem, we can identify the martingale parts and the finite variation parts of these two semi-martingale decompositions of $F^{2}(X_{t})$, and deduce the expected expression of $\langle M \rangle_{F,t}$.
\subsection{Proof of Lemma~\ref{corr:martzt}}
\label{app:lemmamr}
We simply use Lemma~\ref{lemm:yt} to affirm that
\begin{align*} 
F(Z_{t}) = & \hspace{0.1cm} F(Y_{0}) + \displaystyle{\int_{0}^{\pi(t)}F'(Y_{s})\dfrac{g(A_{\xi_{0}}(s),R)}{A_{\xi_{0}}(s)}\bigg(\dfrac{g(A_{\xi_{0}}(s)Y_{s},R)}{g(A_{\xi_{0}}(s),R)}-Y_{s} \bigg) \mathrm{d}s} \\
& + \displaystyle{\int_{0}^{\pi(t)}\int_{\mathbb{R}^{+}}}\mathbb{1}_{\{ h \leq b_{x_{0}}(A_{\xi_{0}}(s)Y_{s-})\}} \left(F\bigg(Y_{s-}-\dfrac{x_{0}}{A_{\xi_{0}}(s)}\bigg) - F(Y_{s-})\right) \mathcal{N}(\mathrm{d}s, \mathrm{d}h).
\end{align*}
Then, we perform the change of variables $u=\pi^{-1}(s)$ in the first integral, and also obtain by definition of $\widehat{\mathcal{N}}$:
\begin{align*} 
F(Z_{t})= &  \hspace{0.1cm} F(Z_{0}) + \displaystyle{\int_{0}^{t}F'(Z_{u})\pi'(u)\dfrac{g(A_{\xi_{0}}(\pi(u)),R)}{A_{\xi_{0}}(\pi(u))}\left( \dfrac{g(A_{\xi_{0}}(\pi(u))Z_{u},R)}{g(A_{\xi_{0}}(\pi(u)),R)} - Z_{u} \right) \mathrm{d}u} \\
& + \displaystyle{\int_{0}^{t}\int_{\mathbb{R}^{+}}}\mathbb{1}_{\{ k \leq \pi'(u) b_{x_{0}}\left(A_{\xi_{0}}(\pi(u))Z_{u-}\right) \}} \left(F\bigg(Z_{u-}-\dfrac{x_{0}}{A_{\xi_{0}}(\pi(u))}\bigg) - F(Z_{u-})\right) \widehat{\mathcal{N}}(\mathrm{d}u, \mathrm{d}k).
\end{align*}
We conclude by definition of $\pi'$, and splitting $\widehat{\mathcal{N}}(\mathrm{d}u, \mathrm{d}k) = \widehat{\mathcal{N}}_{C}(\mathrm{d}u, \mathrm{d}k) + \mathrm{d}u\mathrm{d}k$ thanks to Lemma~\ref{lemm:equapoi}. By the same arguments as in the proof of Proposition~\ref{prop:martingaleprobx}, the highlighted local martingale is a true $L^{2}$-martingale, and we compute $\langle \mathcal{M}_{F} \rangle_{t}$.

\section*{Acknowledgements}
This work was partially supported by the Chaire ``Mod\'elisation Math\'ematique et Biodiversit\'e'' of VEOLIA Environment, \'Ecole Polytechnique, Mus\'eum National d'Histoire Naturelle and Fondation X, and by the European Union (ERC, SINGER, 101054787). Views and opinions expressed are however those of the author(s) only and do not necessarily reflect those of the European Union or the European Research Council. Neither the European Union nor the granting authority can be held responsible for them.

\bibliographystyle{alpha}
\bibliography{biblioarticle1}

\end{document}